\documentclass[hidelinks]{siamart251216}

\usepackage[utf8]{inputenc}
\usepackage[T1]{fontenc}
\usepackage{mathtools}
\usepackage{amssymb}
\usepackage{bm}
\usepackage{mathrsfs}
\usepackage{microtype}
\newcommand{\IC}{{\mathbb C}}
\newcommand{\IR}{{\mathbb R}}
\newcommand{\D}{\mathrm{D}}
\newcommand{\br}{{\normalfont{\boldsymbol{r}}}}
\newcommand{\x}{{\bf{x}}}
\newcommand{\y}{{\bf{y}}}
\newcommand{\z}{{\bf{z}}}

\newcommand{\dotp}[2]{\left ( #1,#2 \right )}
\newcommand{\dual}[2]{\left\langle#1,#2\right\rangle}
\newcommand{\norm}[1]{\left \lVert #1 \right \rVert}
\newcommand{\snorm}[1]{\left \lvert #1 \right \rvert}
\newcommand{\dd}{\,{\rm d}}
\renewcommand{\d}{\!\!\operatorname{d}}
\DeclareMathOperator*{\essinf}{ess\,inf}

\newsiamthm{problem}{Problem}
\newsiamthm{assumption}{Assumption}
\newsiamremark{example}{Example}
\newsiamremark{remark}{Remark}

\headers{Shape Holomorphy of the Helmholtz Calder\'on Projector}
{F. Henr\'iquez, J. D\"olz, and J. Pinto}
\hypersetup{
  pdftitle={Shape Holomorphy of the Calderon Projector on Lipschitz Surfaces: Frequency-Explicit Bounds for the Helmholtz Operator},
  pdfauthor={Fernando Henriquez, Jurgen Dolz, Jose Pinto}
}

\title{Shape Holomorphy of the Calder\'on Projector
on Lipschitz Surfaces: Frequency-Explicit Bounds for the Helmholtz Operator\thanks{\textbf{Funding:} FH's work is funded by the Deutsche Forschungsgemeinschaft (DFG, German
Research Foundation), Project-ID 258734477, SFB 1173, and the Austrian
Science Fund (FWF) under the project I6667-N. Funding was also received from
the European Research Council (ERC) under the Horizon
Europe research and innovation programme of the European Union (Grant agreement No.~101125225,
Project Name: ``New Frontiers in Optimal Adaptivity'', PI: Prof.\ Dr.\
Michael Feischl). JP's work is funded by the ANID grant Fondecyt
Iniciaci\'on No.~11230248. JD gratefully acknowledges funding by the
Deutsche Forschungsgemeinschaft (DFG, German Research Foundation) under
Germany's Excellence Strategy, EXC-2047/1, Project-ID 390685813, and
through the CRC 1720, Project-ID 53930965.}}

\author{Fernando Henr\'iquez\thanks{Institute for Analysis and Scientific Computing,
Technische Universit\"{a}t Wien,
Wiedner Hauptstrasse 8-10, A-1040 Wien, Austria
(\email{fernando.henriquez@asc.tuwien.ac.at}).}
\and
J\"{u}rgen D\"{o}lz\thanks{Institute for Numerical Simulation, University of Bonn,
Friedrich-Hirzebruch-Allee 7, Bonn 53115, Germany
(\email{doelz@ins.uni-bonn.de}).}
\and
Jos\'e Pinto\thanks{Faculty of Engineering and Sciences,
Universidad Adolfo Ib\'añez,
Diagonal Las Torres 2700, Santiago 7941169, Chile
(\email{jose.pinto@uai.cl}).}}

\begin{document}
\maketitle

\begin{abstract}
We establish the holomorphic dependence of the Calder\'on projector for the
Helmholtz operator upon shape boundary deformations of a family of bounded
Lipschitz domains in two and three spatial dimensions, with bounds explicit
in the wavenumber and the size of the holomorphic extension. More precisely,
after transporting the four boundary integral operators contained in the Calder\'on projector
to a fixed reference
boundary, we prove that the resulting parameter-to-operator maps admit bounded
holomorphic extensions in their natural fractional Sobolev spaces.
Unlike approaches that require boundary regularity of class
$\mathscr C^{1,1}$ or higher, our argument traces the holomorphic
dependence on the shape by transferring elliptic estimates to the boundary
through the fundamental solution and the Dirichlet and conormal trace
operators. This construction is the key step that allows us to treat
Lipschitz-regular boundary deformations.

Furthermore, we show that the domain of holomorphy can be chosen independently of the
wavenumber $\kappa$, while the operator bounds in the corresponding fractional Sobolev 
spaces contain algebraic factors
in $\kappa$ multiplied by $\exp(C\kappa\delta)$. Here, $\delta>0$
quantifies the size of the domain of complex holomorphic extension, and
$C>0$ is independent of $\kappa$ and $\delta$. Choosing
$\delta$ proportional to $\kappa^{-1}$ controls the exponential factor
and retains the algebraic growth of the real-domain bounds for the corresponding
operator norm in the natural fractional Sobolev spaces.
In particular, for real domains, corresponding to $\delta=0$, our single
layer estimate recovers, to our knowledge, the best available frequency
dependence on general Lipschitz boundaries. For the double layer, adjoint
double layer, and hypersingular operators, we obtain novel frequency-explicit
bounds in their natural fractional Sobolev spaces under the same
Lipschitz assumptions.

The results proved in this work provide the analytic foundation for
frequency-explicit uncertainty quantification, model order reduction, and high-dimensional
approximation of boundary integral operators, their Galerkin discretizations,
and the resulting boundary densities and quantities of interest on
families of parametrically varying boundaries.
\end{abstract}



\section{Introduction}

We consider the Calder\'on projector for the Helmholtz
operator on the boundaries of a collection of bounded Lipschitz domains
in two and three spatial dimensions. This projector comprises four
boundary integral operators (BIOs), namely the single layer, double
layer, adjoint double layer, and hypersingular operator. We establish
the holomorphic dependence of these operators on domain deformations,
a property referred to as \emph{shape holomorphy}.
This work extends our previous work on shape holomorphy of the Calder\'on projector
for the Laplace operator \cite{HDPS26Laplace} to the time-harmonic scalar wave propagation setting. 
In addition, we quantify how the bounds depend
on the wavenumber $\kappa$ and the size of the domain of holomorphic
extension.

Our interest in the shape holomorphy property of BIOs stems from the
approximation of parametric maps arising from the boundary integral
formulation of boundary value problems on parametrically varying domains, in particular, 
wave scattering problems in unbounded domains.
The dependence of the BIOs and the corresponding boundary densities on
the geometric parameters gives rise to parameter-to-operator and
parameter-to-solution maps. In forward and inverse uncertainty
quantification (UQ), these maps may involve high-dimensional, possibly
countably infinite, parametric inputs.
Approximating these maps efficiently becomes challenging when the number of
parameters is large. For tensor product methods, for instance, the computational
effort grows exponentially with the number of parameters, leading to the
so-called \emph{curse of dimensionality} in the parameter space.

A strategy to overcome this difficulty relies on constructing bounded
holomorphic extensions of the parametric maps to suitable products of
complex domains, as described in \cite{CCS15}. The sizes of these domains reflect the
anisotropic dependence on the individual parametric variables. Together with
appropriate summability assumptions on the parametrization, this
regularity enables sparse polynomial approximation with algebraic
convergence rates independent of the parametric dimension, including
the case of countably infinitely many parameters.

A variety of computational methods rely on this holomorphic dependence
for their convergence analysis, including sparse grid interpolation and
quadrature \cite{ZS17,SS13,NTW2008,HHPS2018} and higher-order quasi-Monte
Carlo integration \cite{DKL14,DLC16}. The same regularity
underlies approximation results for artificial neural network surrogates
\cite{SZ19,DH2026,HS26NeuralShape} and sparsity results for parametric
posterior densities in Bayesian inverse problems \cite{CS16_2}.
These applications require holomorphic extensions and bounds for the
relevant parametric maps. For Helmholtz problems, the
dependence of these bounds on the wavenumber is essential to make the
approximation results explicit in the problem's leading frequency.

\subsection{Previous Work}

The holomorphic dependence of solution maps on shape deformations has
been studied for a variety of boundary value problems. In the context
of variational formulations in the domain, this property has been
established for elliptic diffusion \cite{CNT2016,HPS2016,HS2022},
electromagnetic scattering \cite{JSZ2017}, stationary Navier--Stokes
flow \cite{CSZ18}, and acoustic scattering \cite{HSSS2018}. For Helmholtz
problems, the dependence of the region of parametric holomorphy of
solution maps on the wavenumber is investigated in \cite{SW23}, while
\cite{HSE25} establishes frequency-explicit shape holomorphy for acoustic
scattering. In the context of volume integral formulations for the scattering by an inhomogeneous
obstacle, \cite{henriquez2024domain} proves shape holomorphy of the
Lippmann--Schwinger operator under Lipschitz domain transformations.

For BIOs, earlier work established real shape differentiability
for acoustic, elastic, and electromagnetic scattering
\cite{Pot94,Charam95,Pot96,Pot96_3,LC121,LC122}.
Even when available to arbitrary order, these differentiability results
do not by themselves imply analytic or holomorphic dependence on shape
deformations.

Shape holomorphy of BIOs was established in \cite{HS21,henriquez2021shape}
for the Laplace and Helmholtz Calder\'on projectors on boundaries of
class $\mathscr C^2$ in two spatial dimensions. The corresponding
question for weakly singular and
hypersingular BIOs on multiple open arcs is addressed in \cite{PHJ23}.
In addition, the real analytic dependence of layer potentials has also been studied in
H\"older spaces in \cite{DR04,LR08,dalla2022multi,DLM22}, with an extension
to heat potentials in \cite{DL23}. In fractional Sobolev spaces,
\cite{Henriquez2026C11} proves shape holomorphy of the Helmholtz
Calder\'on projector on closed $\mathscr C^{1,1}$ surfaces and of the
Maxwell EFIE on closed surfaces and screens of the same regularity.
Still, these results require boundary smoothness beyond Lipschitz regularity.

For Lipschitz boundaries,
\cite{dolz2024parametric} establishes parametric holomorphy under
affine boundary deformations in three spatial dimensions following an analysis
that is carried out in an $L^2$-based setting.
Recently, for the Laplace operator, shape holomorphy in fractional Sobolev trace spaces
is established in \cite{HDPS26Laplace} by expressing the BIOs through volume potentials
and trace operators.
In the present work, we extend our analysis of the Laplace operator to
the Helmholtz one with the additional objective of establishing frequency-explicit
bounds for the domain of holomorphy and the corresponding operator norms.

\subsection{Contributions}
The main result of this work (Theorem~\ref{thm:Calderon_Projector})
extends the Laplace analysis of \cite{HDPS26Laplace} to the Helmholtz
Calder\'on projector under Lipschitz domain and boundary deformations
in two and three spatial dimensions. We establish shape holomorphy
of all four BIOs in their natural fractional Sobolev spaces and give
bounds explicit in the wavenumber and the size of the holomorphic extension.

On Lipschitz boundaries, for $\kappa\ge\kappa_0>0$, our bounds in the natural fractional Sobolev spaces have orders
$(1+\kappa)^{\frac{d-1}{2}}$ for the single layer operator,
$(1+\kappa)^{\frac32}$ for the double layer operator and its adjoint,
and $(1+\kappa)^{\frac{d+1}{2}}$ for the hypersingular operator.
The single layer estimate recovers
\cite[Corollary~1.2, item~(i)]{GrahamLohndorfMelenkSpence2015}, to our knowledge
the best available frequency dependence on general Lipschitz boundaries.
The latter three bounds in the natural fractional Sobolev spaces are
new to our knowledge under the same Lipschitz assumptions.

We also quantify the extension of these estimates to complex shape
deformations. The domain of holomorphy can be chosen independently of
$\kappa$, however, the bounds contain algebraic factors in $\kappa$
multiplied by $\exp(C\kappa\delta)$, with $C>0$ independent of $\kappa$
and $\delta$. Choosing $\delta$ proportional to $\kappa^{-1}$ keeps the
exponential factor uniformly bounded and retains the algebraic growth
of the real-domain bounds. To our knowledge, this work and
\cite{SW23,HSE25} are the only works addressing the wavenumber dependence
of holomorphic extensions, and the present work is the first to do so
for BIOs.

We hasten to add that this quantitative control provides the analytic foundation for
frequency-explicit uncertainty quantification, model order reduction,
and high-dimensional approximation on parametrically varying domains.
Indeed, it supplies the operator estimates needed to study Galerkin
discretizations, boundary densities, and quantities of interest.

\subsection{Outline}
\label{sec:Outl}
This work is organized as follows. Section~\ref{sec:preliminaries}
introduces the domain transformations, pullback operators, and trace
spaces used in the analysis, and recalls the Helmholtz potentials and
BIOs on Lipschitz boundaries.
Section~\ref{sec:BiosResults} states the main shape holomorphy theorem
for the Calder\'on projector, with bounds explicit in the wavenumber and
the size of the holomorphic extension. Its proof is developed by first
treating the Newton and layer potentials, and then the four BIOs in
their natural fractional Sobolev spaces.
Finally, in Section~\ref{sec:concluding_remarks}, we provide concluding remarks,
outline directions of future research, and refer to the comparison
with existing frequency bounds. Technical results are
deferred to the appendix.

\subsection{Notation}
\label{sec:Notat}
For (complex) Banach spaces $X$ and $Y$,
we denote by $\mathscr{L}(X,Y)$ the space of bounded linear operators 
from $X$ into $Y$ and by $\mathscr{L}_{\text{iso}}(X,Y)$ the (open)
subset of isomorphisms, i.e.,~bounded linear operators with a bounded inverse. 
Recall that $\mathscr{L}(X,Y)$
is a (complex) Banach space equipped with the standard operator norm (see, e.g., \cite[Theorem~III.2]{RS80_Vol1}).
We also denote by $X'$ the dual space of $X$, i.e., $X' \coloneq \mathscr{L}(X, \IC)$, 
and by $X^\star$ the space of all bounded \emph{anti-linear} functionals acting on $X$.
Dual pairings are denoted by $\dual{\cdot}{\cdot}_{X' \times X}$
and $\dual{\cdot}{\cdot}_{X^{\star}\times X}$.
We denote by
$\mathscr{L}_{\normalfont\text{bl}}(X\times Y,Z)$ the space of bounded
bilinear maps from $X\times Y$ into $Z$, equipped with the standard norm.

The complex Hilbert spaces used below are complexifications of real
Hilbert spaces. We write $\overline{x}$ for the canonical conjugation
and denote by $(\cdot,\cdot)_X$ the complex bilinear extension of the
real inner product. Thus, the Hilbert inner product is
$(x,\overline{y})_X$, and $\|x\|_X^2=(x,\overline{x})_X$, where the explicit reference to $X$ in the norm is dropped whenever $X = \IC^d$ or $X= \IR^d$.
With this convention, the map
$\mathcal{R}_X: X \rightarrow X':x \mapsto \dotp{\cdot}{x}_X$
is a linear isometric isomorphism, obtained by composing the usual
Riesz map to the linear dual with conjugation. For each $x'\in X'$, we have
\begin{equation*}
	\dual{x' }{x}_{X' \times X} = (x, \mathcal{R}_X^{-1} x' )_X, \quad \forall x \in X.
\end{equation*}
For a given $T \in \mathscr{L}(X,Y)$, we define its adjoint 
$T' \in \mathscr{L}(Y',X')$ as the unique operator satisfying
\begin{equation*}
	\dual{y'}{T x}_{Y' \times Y}
	=
	\dual{T' y'}{x}_{X' \times X},
	\quad
	\forall x \in X, \;\forall y' \in Y'.
\end{equation*}

All Lebesgue and Sobolev spaces are understood over $\IC$, unless stated otherwise.
Let $\D \subset \IR^{d}$, $d\geq 2$, be a Lipschitz domain (\cite[Chap.~1.2]{Grisvard}).
Denote by $L^p(\D)$, $p\in[1,\infty]$, the usual Lebesgue spaces on $\D$,
and by $H^s(\D)$, $s \geq 0$, the standard Sobolev spaces in $\D$. We further define
$\widetilde{H}^s(\D)$, $s \geq 0$, as the closure of the smooth compactly supported
functions $\mathscr{C}_0^\infty(\D)$ in the norm of $H^s(\IR^d)$. The spaces $H^s(\D)$ and $\widetilde{H}^s(\D)$, for $s <0$, are defined as the duals of $\widetilde{H}^{-s}(\D)$ and $H^{-s}(\D)$, respectively.
In addition, let $H^s(\partial \D)$, for $s\in[0,1]$, {be the Sobolev}
space of traces on $\partial \D$ (\cite[Section~2.4]{SS10} and \cite[Section~4.2]{HW08}).
As is customary, we identify $H^0(\partial \D)$ with $L^2({\partial \D})$ and, 
for $s\in [0,1]$, define $H^{-s}(\partial \D)$ as the dual space of $H^s(\partial \D)$.

For $\x\in\IR^d$ and $R>0$, we write
$B(\x,R)=\{\y\in\IR^d:\norm{\y-\x}<R\}$ for the open Euclidean ball
centered at $\x$. We use the Euclidean vector norm and its induced matrix norm.
We equip $W^{1,\infty}(\D;\IC^d)$ with the norm
\begin{equation*}
\|\boldsymbol u\|_{W^{1,\infty}(\D;\IC^d)}=
\|\boldsymbol u\|_{L^\infty(\D;\IC^d)}
+\|D\boldsymbol u\|_{L^\infty(\D;\IC^{d\times d})},
\end{equation*}
where $D\boldsymbol u$ denotes the Jacobian matrix. The corresponding real-valued spaces are equipped with the same norms. We also denote by $\operatorname{Lip}(\br; \D)$ the Lipschitz constant of
$\br \in W^{1,\infty}(\D; \IC^d)$.

\section{Preliminaries}
\label{sec:preliminaries}
\subsection{Domain Transformations}
\label{sec:domTransforms}
We consider domains in $\mathbb{R}^d$ which are images of one common,
bounded Lipschitz \emph{nominal} domain $\widehat{\D}\subset\IR^d$ under a
bijective and bi-Lipschitz map. 
That is, we consider a collection $\{\D_\br\}_{\br\in\mathfrak{T}}$ of parametrized domains of
the form $\D_\br\coloneqq\br(\widehat{\D})$, where $\mathfrak{T}$
denotes the set of admissible domain transformations.
We write
$\widehat{\Gamma}\coloneqq\partial\widehat{\D}$ and, for each
$\br\in\mathfrak{T}$, $\Gamma_\br\coloneqq\partial\D_\br=\br(\widehat{\Gamma})$.
We assume that the set of admissible domain transformations satisfies the
following conditions, cf.~\cite[Assumption~2.1]{HDPS26Laplace}.

\begin{assumption}\label{assump:transformations_T}
Let $\widehat{\D}\subset\IR^d$,  $d\geq 2$,  
be a bounded Lipschitz domain. The set 
$\mathfrak{T}\subset W^{1,\infty}(\widehat\D;\IR^d)$ of admissible domain transformations is 
such that:
\begin{itemize}
	\item[\normalfont{(i)}] 
	$\mathfrak{T}$ is compact as a subset of $W^{1,\infty}(\widehat\D;\IR^d)$.
	\item[\normalfont{(ii)}] 
	Each $\br\in\mathfrak{T}$ is bijective and bi-Lipschitz 
	when considered as a map $\br\colon\widehat{\D}\to\D_\br$.
	\item[\normalfont{(iii)}] 
	For each $\br\in\mathfrak{T}$, the domain $\br(\widehat{\D})$ is a bounded Lipschitz domain in $\mathbb{R}^d$. 
	\item[\normalfont{(iv)}] 
	For each $\br \in \mathfrak{T}$, it holds that
\begin{equation*}\label{eq:LwrBd}
		\essinf_{\widehat\x \in \widehat{\D}}
		J_\br(\widehat\x) 
		> 0.
	\end{equation*}
\end{itemize}
Here $J_\br \coloneqq \det D\br$, for each $\br \in \mathfrak{T}$,
where $D\br$ denotes the Jacobian matrix of $\br \in \mathfrak{T}$\footnote{
Observe that $D\br$ denotes the Jacobian matrix of $\br$, 
whereas $\widehat{\text{D}}$ refers to the reference domain.}.
\end{assumption}

Establishing shape holomorphy requires complex neighborhoods of the
domain transformations $\mathfrak{T}$. To this end, for
$\mathfrak{T}$ satisfying Assumption~\ref{assump:transformations_T}
and for $\varepsilon>0$, we define
\begin{equation*}
	\mathfrak{T}_{\varepsilon}
	\coloneqq
	\left\{
		\br 
		\in 
		W^{1,\infty}(\widehat\D;\IC^d)
		:
		\exists\,
		\widetilde{\br} \in \mathfrak{T}
		\text{ such that }
		\norm{\br - \widetilde{\br}}_{W^{1,\infty}(\widehat\D;\IC^d)}
		<
		\varepsilon
	\right\}.
\end{equation*}
The following two technical results will be used repeatedly in this work.
\begin{lemma}[{\cite[Lemma~A.2]{HDPS26Laplace}}]
\label{lem:lipschitz_preservation}
Let $\mathfrak T$ be a set of domain transformations satisfying
Assumption~\ref{assump:transformations_T}.
Then there exist $\eta>0$ and $\delta>0$ such that
\begin{equation*}
\inf_{\br\in\mathfrak T_\delta}
\inf_{\substack{\widehat\x,\widehat\y\in\widehat\Gamma\\
\widehat\x\ne\widehat\y}}
\operatorname{Re}\left\{
\frac{(\br(\widehat\x)-\br(\widehat\y))\cdot
(\br(\widehat\x)-\br(\widehat\y))}
{\norm{\widehat\x-\widehat\y}^2}\right\}
\geq\eta>0,
\end{equation*}
and
\begin{equation*}
\inf_{\br\in\mathfrak T_\delta}
\inf_{\substack{\widehat\x,\widehat\y\in\widehat\D\\
\widehat\x\ne\widehat\y}}
\operatorname{Re}\left\{
\frac{(\br(\widehat\x)-\br(\widehat\y))\cdot
(\br(\widehat\x)-\br(\widehat\y))}
{\norm{\widehat\x-\widehat\y}^2}\right\}
\geq\eta>0.
\end{equation*}
\end{lemma}

\begin{lemma}[{\cite[Lemma~C.1]{HDPS26Laplace}}]
\label{lem:shape_hol_bilinear_forms}
Let $\mathfrak T$ be a set of domain transformations satisfying
Assumption~\ref{assump:transformations_T}.
For each $\br\in\mathfrak T$, we introduce the bilinear forms
$\mathsf a_\br:H^1(\widehat\D)\times H^1(\widehat\D)\rightarrow\IC$
and
$\mathsf b_\br:L^2(\widehat\D)\times L^2(\widehat\D)\rightarrow\IC$, which are defined as follows
\begin{equation*}
\begin{aligned}
\mathsf a_\br(\widehat u,\widehat v)
&\coloneqq\int_{\widehat\D}
(D\br)^{-\top}\widehat\nabla\widehat u
\cdot(D\br)^{-\top}\widehat\nabla\widehat v
J_\br\,\dd\widehat\x,\\
\mathsf b_\br(\widehat u,\widehat v)
&\coloneqq\int_{\widehat\D}
\widehat u\widehat vJ_\br\,\dd\widehat\x.
\end{aligned}
\end{equation*}
There exists $\delta>0$ such that the maps
\begin{equation*}
\begin{aligned}
\mathcal A:&\ \mathfrak T\rightarrow
\mathscr L_{\normalfont\text{bl}}
(H^1(\widehat\D)\times H^1(\widehat\D),\IC):
\br\mapsto\mathsf a_\br,\\
\mathcal B:&\ \mathfrak T\rightarrow
\mathscr L_{\normalfont\text{bl}}
(L^2(\widehat\D)\times L^2(\widehat\D),\IC):
\br\mapsto\mathsf b_\br,
\end{aligned}
\end{equation*}
admit holomorphic extensions to $\mathfrak T_\delta$, which
are bounded according to
\begin{equation*}
\begin{aligned}
\sup_{\br\in\mathfrak T_\delta}
\norm{\mathsf a_\br}_{\mathscr L_{\normalfont\text{bl}}
(H^1(\widehat\D)\times H^1(\widehat\D),\IC)}
&\leq C_{\mathsf a}(\mathfrak T,\delta)
\coloneqq\sup_{\br\in\mathfrak T_\delta}
\norm{(D\br)^{-\top}}_{L^\infty(\widehat\D;\IC^{d\times d})}^2
\norm{J_\br}_{L^\infty(\widehat\D)},\\
\sup_{\br\in\mathfrak T_\delta}
\norm{\mathsf b_\br}_{\mathscr L_{\normalfont\text{bl}}
(L^2(\widehat\D)\times L^2(\widehat\D),\IC)}
&\leq C_{\mathsf b}(\mathfrak T,\delta)
\coloneqq\sup_{\br\in\mathfrak T_\delta}
\norm{J_\br}_{L^\infty(\widehat\D)}.
\end{aligned}
\end{equation*}
\end{lemma}

\subsection{Pullback Operator}
\label{sec:Pullbck}
Let $\mathfrak{T}$ be a set of domain transformations satisfying Assumption~\ref{assump:transformations_T}.
Following \cite[Section~2.4]{HDPS26Laplace}, for each $\br \in \mathfrak{T}$, we consider Riesz maps
associated with the complex Sobolev spaces $H^1(\D_\br)$ and $H^1(\widehat{\D})$, i.e.,
\begin{align}\label{eq:riesz_maps_volume}
    \mathcal{R}_{\br}: H^1(\D_\br) \rightarrow \widetilde{H}^{-1}(\D_\br)
    \quad
    \text{and}
    \quad
    \widehat{\mathcal{R}}: H^1(\widehat{\D}) \rightarrow \widetilde{H}^{-1}(\widehat{\D}),
\end{align}
respectively. 
Furthermore, when dealing with boundary densities we consider the Riesz maps
corresponding to the complex Sobolev spaces $H^{\frac{1}{2}}(\Gamma_\br)$ and $H^\frac{1}{2}(\widehat{\Gamma})$, i.e.,
\begin{align}\label{eq:riesz_maps_boundary}
    \mathcal{R}_\br: H^{\frac{1}{2}}(\Gamma_\br) \rightarrow {H}^{-\frac{1}{2}}(\Gamma_\br)
    \quad
    \text{and}
    \quad
    \widehat{\mathcal{R}}: H^\frac{1}{2}(\widehat{\Gamma}) \rightarrow {H}^{-\frac{1}{2}}(\widehat{\Gamma}).
\end{align}
The precise meaning will be clear from the context.

Next, for each $\br \in \mathfrak{T}$, we define the \emph{pullback} operator 
\begin{equation}\label{eq:pullback_operator}
    \tau_{\br} u
    \coloneqq
    \left\{
    \begin{array}{cc}
        u\circ \br, & u\in H^1(\D_\br)\ \text{or}\ u\in L^2(\D_\br), \\
        \widehat{\mathcal{R}}
        \tau_\br
        \mathcal{R}^{-1}_{\br}
        u, & u\in \widetilde{H}^{-1}(\D_\br),
    \end{array}
    \right.
\end{equation}
where $\mathcal{R}_\br$ and $\widehat{\mathcal{R}}$ are the Riesz maps in
\eqref{eq:riesz_maps_volume}. Its boundary counterpart is
\begin{equation}
\label{eq:pullbackBoundary}
    \tau_{\br} \varphi
    \coloneqq
    \left\{
    \begin{array}{cc}
        \varphi\circ \br, & \varphi\in H^{\frac{1}{2}}(\Gamma_\br)\ \text{or}\ \varphi\in L^2(\Gamma_\br), \\
        \widehat{\mathcal{R}}
        \tau_\br 
        \mathcal{R}^{-1}_{\br}
        \varphi, & \varphi \in H^{-\frac{1}{2}}(\Gamma_\br),
    \end{array}
    \right.
\end{equation}
where $\mathcal{R}_\br$ and $\widehat{\mathcal{R}}$ are the Riesz maps in
\eqref{eq:riesz_maps_boundary}.
The following result states the mapping properties of the pullback operator.

\begin{proposition}[{\cite[Proposition~2.2]{HDPS26Laplace}}]\label{prop:iso_pullback}
Let $\mathfrak{T}$ be a set of domain transformations satisfying Assumption~\ref{assump:transformations_T}. 
Then, for each $\br\in\mathfrak{T}$,
\begin{equation*}
\begin{aligned}
\tau_\br&\in
\mathscr{L}_{\normalfont\text{iso}}(H^1(\D_\br),H^1(\widehat{\D})),
&\quad
\tau_\br&\in
\mathscr{L}_{\normalfont\text{iso}}(\widetilde{H}^{-1}(\D_\br),
\widetilde{H}^{-1}(\widehat{\D})),\\
\tau_\br&\in
\mathscr{L}_{\normalfont\text{iso}}(H^{\frac{1}{2}}(\Gamma_\br),
H^{\frac{1}{2}}(\widehat{\Gamma})),
&
\tau_\br&\in
\mathscr{L}_{\normalfont\text{iso}}(H^{-\frac{1}{2}}(\Gamma_\br),
H^{-\frac{1}{2}}(\widehat{\Gamma})).
\end{aligned}
\end{equation*}
Furthermore, the maps $\br \mapsto \tau_\br$ and $\br \mapsto \tau^{-1}_\br$,
as elements of the space of bounded linear operators with the above mapping properties,
are uniformly bounded in $\mathfrak{T}$.
\end{proposition}

Following \cite[Section~2.4]{HDPS26Laplace}, for each $\br \in \mathfrak{T}$, we define the \emph{bilinear} form
$\dual{\cdot}{\cdot}_{\widetilde{H}^{-1}(\widehat\D) \times H^1(\widehat\D),\br}:
\widetilde{H}^{-1}(\widehat\D) \times H^1(\widehat\D) \rightarrow \mathbb{C}$
for any $\widehat{u} \in \widetilde{H}^{-1}(\widehat{\D})$ and any
$\widehat{v}\in H^1(\widehat{\D})$ as follows:
\begin{equation}\label{eq:bilnfrmH1_2}
\begin{aligned}
  	\dual{
		\widehat{u}
	}{
		\widehat{v}
	}_{\widetilde{H}^{-1}(\widehat\D) \times H^1(\widehat\D),\br}
    &\coloneqq
    \int_{\widehat\D}
    \left(D\br\right)^{-\top}
    \widehat{\nabla} ({\widehat{\mathcal{R}}}^{-1}\widehat{u})
    \cdot
    \left(D\br\right)^{-\top}
    \widehat{\nabla}\widehat{v}
    J_{\br}
    \mathrm{d}\widehat{\x}
    \\
    &\quad+
    \int_{\widehat\D}
    ({\widehat{\mathcal{R}}}^{-1}\widehat{u})
    \widehat{v}
    J_{\br}
    \mathrm{d}\widehat{\x}.
\end{aligned}
\end{equation}
Moreover, for each $\br\in\mathfrak{T}$, $u\in\widetilde H^{-1}(\D_\br)$,
and $v\in H^1(\D_\br)$, we have that
\begin{equation*}\label{eq:transport_domain}
	\dual{u}{{v}}_{\widetilde{H}^{-1}(\D_{\br}) \times H^1(\D_{\br})}
	=
    	\dual{
    		\tau_\br {u}
	}{
		\tau_\br {v}
	}_{\widetilde{H}^{-1}(\widehat\D) \times H^1(\widehat\D),\br},
\end{equation*}
with the pullback operator $\tau_\br$
as in \eqref{eq:pullback_operator}.

Next, we extend the definition of the bilinear form
$\dual{\cdot}{\cdot}_{\widetilde{H}^{-1}(\widehat\D) \times H^1(\widehat\D),\br}:
\widetilde{H}^{-1}(\widehat\D) \times H^1(\widehat\D) \rightarrow \mathbb{C}$
to $\br\in\mathfrak{T}_\delta$, for some $\delta>0$, and recall its main properties.

\begin{lemma}[{\cite[Lemma~2.5]{HDPS26Laplace}}]\label{lem:bilinear_ext_D}
Let $\mathfrak{T}$ be a set of domain transformations
satisfying Assumption~\ref{assump:transformations_T}.
\begin{itemize}
\item[\normalfont{(i)}]
There exists $\delta>0$ such that the bilinear form
$ \dual{\cdot}{\cdot}_{\widetilde{H}^{-1}(\widehat\D) \times H^1(\widehat\D),\br}:
\widetilde{H}^{-1}(\widehat\D) \times H^1(\widehat\D) \rightarrow \mathbb{C}$, $\br \in \mathfrak{T}$,
introduced in \eqref{eq:bilnfrmH1_2} admits a bounded extension to $\mathfrak{T}_\delta$.
More precisely, for each $\br \in \mathfrak{T}_\delta$
and for all  $\widehat{u} \in \widetilde{H}^{-1}(\widehat{\D})$, $\widehat{v}\in H^1(\widehat{\D})$, it holds that
\begin{equation}\label{eq:C_1_domain}
	\snorm{
		\dual{\widehat{u}}{\widehat{v}}_{\widetilde{H}^{-1}(\widehat\D) \times H^1(\widehat\D),\br}
	}
	\leq
	C_1(\mathfrak{T},\delta)
    \norm{\widehat{u}}_{\widetilde{H}^{-1}(\widehat\D)}
    \norm{\widehat{v}}_{H^1(\widehat\D)},
\end{equation}
where $C_1(\mathfrak{T},\delta)>0$ depends only on $\mathfrak{T}$
and $\delta>0$.
\item[\normalfont{(ii)}]
There exists $\delta>0$ such that the map
\begin{equation*}
    \mathfrak{T} \rightarrow 
    \mathscr{L}_{\normalfont\text{bl}}
      \left(
        \widetilde{H}^{-1}(\widehat{\D})
        \times
        H^{1}(\widehat{\D}),
        \mathbb{C}
    \right):
    \br
    \mapsto
    \left(
        (\widehat{u},\widehat{v})
        \mapsto
        \dual{
        		\widehat{u}
		}{
			\widehat{v}
		}_{\widetilde{H}^{-1}(\widehat{\D}) \times H^{1}(\widehat{\D}),\br}
    \right)
\end{equation*}
admits a bounded holomorphic extension to $\mathfrak{T}_\delta$.
\item[\normalfont{(iii)}]
There exist $\delta>0$ and $C_2(\mathfrak{T},\delta)>0$, 
depending only on $\mathfrak{T}$ and $\delta$, 
such that, for each $\br \in \mathfrak{T}_\delta$,
\begin{equation}\label{eq:inf_sup_sesq_r_D}
\begin{aligned}
    \inf_{\widehat{u}\in \widetilde{H}^{-1}(\widehat{\D})\setminus\{0\}}
    \sup_{\widehat{v}\in H^{1}(\widehat{\D})\setminus\{0\}}
    \frac{
        \snorm{\dual{\widehat{u}}{\widehat{v}}_{\widetilde{H}^{-1}(\widehat{\D}) \times H^{1}(\widehat{\D}),\br}}
    }{
        \norm{\widehat{u}}_{\widetilde{H}^{-1}(\widehat{\D})}
        \norm{\widehat{v}}_{H^{1}(\widehat{\D})}
    }
    &
    \geq C_2(\mathfrak{T},\delta)
    \quad
    \text{and}
    \\
    \inf_{\widehat{v}\in H^{1}(\widehat{\D})\setminus\{0\}}
    \sup_{\widehat{u}\in \widetilde{H}^{-1}(\widehat{\D})\setminus\{0\}}
    \frac{
        \snorm{\dual{\widehat{u}}{\widehat{v}}_{\widetilde{H}^{-1}(\widehat{\D}) \times H^{1}(\widehat{\D}),\br}}
    }{
        \norm{\widehat{u}}_{\widetilde{H}^{-1}(\widehat{\D})}
        \norm{\widehat{v}}_{H^{1}(\widehat{\D})}
    }
    &
    \geq C_2(\mathfrak{T},\delta).
\end{aligned}
\end{equation}
\end{itemize}
\end{lemma}

For each $\br\in\mathfrak T$, let $J_{\mathrm S,\br}$ denote the surface
Jacobian, therefore for any $ g\in L^1(\Gamma_\br)$ it holds that
\begin{equation*}
\int_{\Gamma_\br}g\,\dd s
=
\int_{\widehat\Gamma}(g\circ\br)J_{\mathrm S,\br}\,\dd s,
\end{equation*}
The surface Jacobian is defined through tangential derivatives of local
Lipschitz parametrizations.
For its definition and, in particular, the construction of its bounded holomorphic extension we refer to
\cite[Lemma~2.4, item~(i)]{HDPS26Laplace}.

Furthermore, for each $\br \in \mathfrak{T}$, 
we define the \emph{bilinear} form
$ \dual{\cdot}{\cdot}_{{H}^{-\frac{1}{2}}(\widehat\Gamma) \times H^{\frac{1}{2}}(\widehat\Gamma),\br}:
{H}^{-\frac{1}{2}}(\widehat\Gamma) \times H^{\frac{1}{2}}(\widehat\Gamma) \rightarrow \mathbb{C}$, for $\widehat{\phi} \in {H}^{-\frac{1}{2}}(\widehat\Gamma)$
and $\widehat{\varphi}\in H^\frac{1}{2}(\widehat\Gamma)$, as follows:
\begin{equation}\label{eq:bilinear_form_boundary}
\begin{aligned}
  &\dual{
		\widehat{\phi} 
	}{
		\widehat{\varphi} 
	}_{{H}^{-\frac{1}{2}}(\widehat\Gamma) \times H^{\frac{1}{2}}(\widehat\Gamma),\br}
    \\
    &\coloneqq
    \int_{\widehat\Gamma} 
    (\widehat{\mathcal{R}}^{-1} \widehat{\phi})
    {\widehat{\varphi}}
    J_{\text{S},\br}
    \dd s_{\widehat{\x}}
    \\
    &\quad +
    \int_{\widehat\Gamma} 
    \int_{\widehat\Gamma} 
    \frac{
        ((\widehat{\mathcal{R}}^{-1} \widehat{\phi})(\widehat{\x})-(\widehat{\mathcal{R}}^{-1} \widehat{\phi})(\widehat{\y}))
        {(\widehat{\varphi}(\widehat{\x})-\widehat{\varphi}(\widehat{\y}))}
    }{
        \norm{\br(\widehat{\x})-\br(\widehat{\y})}^d
    }
    J_{\text{S},\br}(\widehat{\x})
    J_{\text{S},\br}(\widehat{\y})
    \dd s_{\widehat{\x}}
    \dd s_{\widehat{\y}}.
\end{aligned}
\end{equation}

Observe that for each $\br\in\mathfrak{T}$ and for all
$\phi\in H^{-\frac{1}{2}}(\Gamma_\br)$ and
$\varphi\in H^{\frac{1}{2}}(\Gamma_\br)$, it holds that
\begin{equation}\label{eq:transport_boundary}
\dual{\phi}{\varphi}_{H^{-\frac{1}{2}}(\Gamma_\br)
\times H^{\frac{1}{2}}(\Gamma_\br)}
=
\dual{\tau_\br\phi}{\tau_\br\varphi}_{
H^{-\frac{1}{2}}(\widehat\Gamma)
\times H^{\frac{1}{2}}(\widehat\Gamma),\br},
\end{equation}
where the pullback operator is defined in \eqref{eq:pullbackBoundary}.
In the following, we recall the equivalent result to \Cref{lem:bilinear_ext_D} for bilinear forms $\dual{\cdot}{\cdot}_{{H}^{-\frac{1}{2}}(\widehat\Gamma) \times H^{\frac{1}{2}}(\widehat\Gamma),\br}$,
$\br \in \mathfrak{T}$.

\begin{lemma}[{\cite[Lemma~2.7]{HDPS26Laplace}}]
\label{lem:bilinear_Gamma}
Let $\mathfrak{T}$ be a set of domain transformations
satisfying Assumption \ref{assump:transformations_T}.
\begin{itemize}
\item[\normalfont{(i)}]
There exists $\delta>0$ such that 
the bilinear form 
$ \dual{\cdot}{\cdot}_{{H}^{-\frac{1}{2}}(\widehat\Gamma) \times H^\frac{1}{2}(\widehat\Gamma),\br}:
{H}^{-\frac{1}{2}}(\widehat\Gamma) \times H^\frac{1}{2}(\widehat\Gamma) \rightarrow \mathbb{C}$ 
introduced in \eqref{eq:bilinear_form_boundary}
admits a bounded extension to $\mathfrak{T}_\delta$.
More precisely, for each $\br \in \mathfrak{T}_\delta$ 
and for all $\widehat{\phi} \in {H}^{-\frac{1}{2}}(\widehat\Gamma)$,
$\widehat{\varphi}  \in H^{\frac{1}{2}}(\widehat\Gamma)$, it holds that
\begin{equation}\label{eq:C_1_boundary}
	\snorm{
  		\dual{
			\widehat{\phi} 
		}{
			\widehat{\varphi} 
		}_{{H}^{-\frac{1}{2}}(\widehat\Gamma) \times H^{\frac{1}{2}}(\widehat\Gamma),\br}
	}
	\leq
	C_1(\mathfrak{T},\delta)
    \norm{\widehat{\phi} }_{{H}^{-\frac{1}{2}}(\widehat\Gamma)}
    \norm{\widehat{\varphi} }_{{H}^{\frac{1}{2}}(\widehat\Gamma)},
\end{equation}
where $C_1(\mathfrak{T},\delta)>0$ depends only on $\mathfrak{T}$
and $\delta>0$.
\item[\normalfont{(ii)}]
There exists $\delta>0$ such that the map
\begin{equation*}
    \mathfrak{T} \rightarrow 
    \mathscr{L}_{\normalfont\text{bl}}
    \left(
    H^{-\frac{1}{2}}(\widehat\Gamma)  \times H^{\frac{1}{2}}(\widehat\Gamma)
    ,
    \mathbb{C}
    \right):
    \br
    \mapsto
    \left(
        (\widehat\phi,\widehat\psi)
        \mapsto
        \dual{\widehat\phi}{\widehat\psi}_{H^{-\frac{1}{2}}(\widehat\Gamma)  \times H^{\frac{1}{2}}(\widehat\Gamma),\br}
    \right)
\end{equation*}
admits a bounded holomorphic extension to $\mathfrak{T}_\delta$.
\item[\normalfont{(iii)}]
There exist $\delta>0$ and $C_2(\mathfrak{T},\delta)>0$, depending only
on $\mathfrak{T}$ and $\delta$, such that
for each $\br \in \mathfrak{T}_\delta$
\begin{equation}\label{eq:inf_sup_sesq_r_boundary}
\begin{aligned}
    \inf_{\widehat{\phi}\in {H}^{-\frac{1}{2}}(\widehat{\Gamma})\setminus\{0\}}
    \sup_{\widehat{\psi}\in H^{\frac{1}{2}}(\widehat{\Gamma})\setminus\{0\}}
    \frac{
        \snorm{\dual{\widehat\phi}{\widehat\psi}_{H^{-\frac{1}{2}}(\widehat\Gamma)  \times H^{\frac{1}{2}}(\widehat\Gamma),\br}}
    }{
        \norm{\widehat{\phi}}_{{H}^{-\frac{1}{2}}(\widehat{\Gamma})}
        \norm{\widehat{\psi}}_{{H}^{\frac{1}{2}}(\widehat{\Gamma})}
    }
    &
    \geq C_2(\mathfrak{T},\delta)
    \quad
    \text{and}
    \\
    \inf_{\widehat{\psi}\in H^{\frac{1}{2}}(\widehat{\Gamma})\setminus\{0\}}
    \sup_{\widehat{\phi}\in {H}^{-\frac{1}{2}}(\widehat{\Gamma})\setminus\{0\}}
    \frac{
        \snorm{\dual{\widehat\phi}{\widehat\psi}_{H^{-\frac{1}{2}}(\widehat\Gamma)  \times H^{\frac{1}{2}}(\widehat\Gamma),\br}}
    }{
        \norm{\widehat{\phi}}_{{H}^{-\frac{1}{2}}(\widehat{\Gamma})}
        \norm{\widehat{\psi}}_{{H}^{\frac{1}{2}}(\widehat{\Gamma})}
    }
    &
    \geq C_2(\mathfrak{T},\delta).
\end{aligned}
\end{equation}
\end{itemize}
\end{lemma}

\subsection{Boundary Integral Operators}
\label{sec:bios_Lipschitz}
In this subsection, we recall the definitions and main properties of
Helmholtz BIOs acting between Sobolev spaces on the boundaries of
Lipschitz domains.
We refer the reader to \cite[Chapters~6, 7, and~9]{mclean2000strongly}
for further details.

Let $d\in\{2,3\}$ and $\kappa>0$, and set
$\mathcal{L}^{(\kappa)}\coloneqq-\Delta-\kappa^2$.
For $\br\in\mathfrak T$ and $\varphi\in L^2(\D_\br)$,
we define the Newton potential in $\D_\br$ as follows:
\begin{equation}\label{eq:newtonpot}
(\mathsf N^{(\kappa)}_\br\varphi)(\x)
\coloneqq\int_{\D_\br}G^{(\kappa)}(\x-\y)\varphi(\y)\,\dd\y,
\quad \x\in\D_\br,
\end{equation}
where $G^{(\kappa)}:\IR^d\setminus\{\boldsymbol0\}\rightarrow\IC$
is the outgoing fundamental solution of the Helmholtz operator,
\begin{equation}\label{eq:green_G_omega}
G^{(\kappa)}(\z)\coloneqq
\begin{cases}
\displaystyle\frac{\imath}{4}\mathrm{H}^{(1)}_0(\kappa\norm{\z}),&d=2,\\[6pt]
\displaystyle\frac{e^{\imath\kappa\norm{\z}}}{4\pi\norm{\z}},&d=3,
\end{cases}
\quad \z\in\mathbb R^d\setminus\{\boldsymbol0\}.
\end{equation}
Here $\mathrm{H}^{(1)}_0$ is the Hankel function of the first kind and
order zero. For each $\br\in\mathfrak T$, $\mathsf N^{(\kappa)}_\br$
is linear and bounded according to
\begin{equation*}
\mathsf N^{(\kappa)}_\br\in
\mathscr L(\widetilde H^{-1}(\D_\br),H^1(\D_\br)),
\end{equation*}
see, e.g., \cite[Theorem~6.1]{mclean2000strongly}. Furthermore, we also have that
$\mathcal L^{(\kappa)}\mathsf N^{(\kappa)}_\br\varphi=\varphi$
for all $\varphi\in L^2(\D_\br)$.

Next, we briefly recall the definition of trace operators on $\Gamma_\br$,
$\br\in\mathfrak T$. For each $\br\in\mathfrak T$, the Dirichlet trace is a bounded,
surjective, linear operator
$\gamma_{\Gamma_\br}:H^1(\D_\br)\rightarrow H^{\frac{1}{2}}(\Gamma_\br)$. We denote its right inverse, being a bounded linear operator itself,
by $Z_{\D_\br}:H^{\frac{1}{2}}(\Gamma_\br)\rightarrow H^1(\D_\br)$.
For further technical details regarding the Dirichlet trace, we refer to, e.g.,
\cite[Theorem~3.37]{mclean2000strongly}.

For each $\br\in\mathfrak T$, we denote by $\gamma_{1,\Gamma_\br}$
the Neumann trace, which, as usual, is defined for
$u\in H^1(\D_\br;\Delta)\coloneqq
\{u\in H^1(\D_\br):\Delta u\in L^2(\D_\br)\}$
through integration by parts as follows
\begin{equation}\label{eq:defneumanntrace}
\dual{\gamma_{1,\Gamma_\br}u}{\varphi}_
{H^{-\frac{1}{2}}(\Gamma_\br)\times H^{\frac{1}{2}}(\Gamma_\br)}
=\int_{\D_\br}\bigl(\Delta u\,Z_{\D_\br}\varphi
+\nabla u\cdot\nabla Z_{\D_\br}\varphi\bigr)\,\dd\x,
\end{equation}
for all $\varphi\in H^{\frac{1}{2}}(\Gamma_\br)$.
Moreover, for each $\br\in\mathfrak T$ and any
$u\in H^1(\D_\br;\Delta)$
\begin{equation*}
\norm{\gamma_{1,\Gamma_\br}u}_{H^{-\frac{1}{2}}(\Gamma_\br)}
\leq\norm{Z_{\D_\br}}_{\mathscr L(H^{\frac{1}{2}}(\Gamma_\br),H^1(\D_\br))}
\bigl(\norm{u}_{H^1(\D_\br)}+\norm{\Delta u}_{L^2(\D_\br)}\bigr),
\end{equation*}
see, e.g., \cite[Lemma~4.3]{mclean2000strongly}.
The definition of $\gamma_{1,\Gamma_\br}$ is independent of the
choice of $Z_{\D_\br}$; see, e.g., \cite[Remark~2.2]{HDPS26Laplace}.

The single- and double-layer potentials are defined, for each
$\br\in\mathfrak T$ and smooth densities, as
\begin{align}
(\mathsf S^{(\kappa)}_{\Gamma_\br}\varphi)(\x)
&\coloneqq\int_{\Gamma_\br}G^{(\kappa)}(\x-\y)
\varphi(\y)\,\dd s_\y,
&&\x\in\D_\br,
\label{eq:single_layer}\\
(\mathsf D^{(\kappa)}_{\Gamma_\br}\phi)(\x)
&\coloneqq\int_{\Gamma_\br}\bm\nu_\br(\y)\cdot
\nabla_\y G^{(\kappa)}(\x-\y)\phi(\y)\,\dd s_\y,
&&\x\in\D_\br,
\label{eq:double_layer}
\end{align}
where $\bm\nu_\br$ denotes the outward unit normal to $\D_\br$.
Standard analysis (see, e.g., \cite{mclean2000strongly,HW08,SS10})
of these volume potential operators shows that
\begin{equation*}
\mathsf S^{(\kappa)}_{\Gamma_\br}\in
\mathscr L(H^{-\frac{1}{2}}(\Gamma_\br),H^1(\D_\br))
\quad\text{and}\quad
\mathsf D^{(\kappa)}_{\Gamma_\br}\in
\mathscr L(H^{\frac{1}{2}}(\Gamma_\br),H^1(\D_\br)).
\end{equation*}

Moreover, for each $\br \in \mathfrak{T}$, a direct computation shows that
\begin{equation}\label{eq:single_double_sol}
\mathcal L^{(\kappa)}\mathsf S^{(\kappa)}_{\Gamma_\br}\varphi=0,
\quad
\mathcal L^{(\kappa)}\mathsf D^{(\kappa)}_{\Gamma_\br}\phi=0,
\end{equation}
for all $\varphi\in H^{-\frac{1}{2}}(\Gamma_\br)$ and
$\phi\in H^{\frac{1}{2}}(\Gamma_\br)$, respectively.
From Green's identities, it is possible to obtain the following representation formula:
\begin{equation}\label{eq:repformula}
u=\mathsf N^{(\kappa)}_\br\mathcal L^{(\kappa)}u
+\mathsf S^{(\kappa)}_{\Gamma_\br}\gamma_{1,\Gamma_\br}u
-\mathsf D^{(\kappa)}_{\Gamma_\br}\gamma_{\Gamma_\br}u,
\quad u\in H^1(\D_\br;\Delta),
\end{equation}
in $\D_\br$.
The same integrals define the exterior layer potentials, which belong to
$H^1_{\mathrm{loc}}(\IR^d\setminus\overline{\D_\br})$, solve the homogeneous
Helmholtz equation, and satisfy the Sommerfeld radiation condition.
We denote the exterior Dirichlet and Neumann traces by
$\gamma^c_{\Gamma_\br}$ and $\gamma^c_{1,\Gamma_\br}$, respectively,
using the same normal $\bm\nu_\br$ as for the interior traces.
For a radiating solution of the homogeneous exterior problem,
Green's representation formula gives
\begin{equation*}\label{eq:repformula_exterior}
u=\mathsf D^{(\kappa)}_{\Gamma_\br}\gamma^c_{\Gamma_\br}u
-\mathsf S^{(\kappa)}_{\Gamma_\br}\gamma^c_{1,\Gamma_\br}u
\quad\text{in }\mathbb R^d\setminus\overline{\D_\br}.
\end{equation*}

For the interior and exterior Cauchy traces, respectively, we write
\begin{equation*}
\boldsymbol\gamma u\coloneqq
(\gamma_{\Gamma_\br}u,\gamma_{1,\Gamma_\br}u)^\top,
\quad
\boldsymbol\gamma^c u\coloneqq
(\gamma^c_{\Gamma_\br}u,\gamma^c_{1,\Gamma_\br}u)^\top.
\end{equation*}
The jump and average are defined by
\begin{equation}\label{eq:jumps}
[\boldsymbol\gamma u]\coloneqq
\boldsymbol\gamma^c u-\boldsymbol\gamma u,
\quad
\{\boldsymbol\gamma u\}\coloneqq
\frac{1}{2}(\boldsymbol\gamma^c u+\boldsymbol\gamma u),
\end{equation}
whenever both traces exist. We use the same notation for their components.

Upon these preparations, 
we recall the four BIOs in the Calder\'{o}n-projector on $\Gamma_\br$: 
\begin{itemize}
    \item {\bf Single Layer BIO:}
    \begin{equation}\label{eq:def_single_layer}
    	\mathsf{V}^{(\kappa)}_{\Gamma_\br} 
	\coloneqq \gamma_{\Gamma_\br} \mathsf{S}^{(\kappa)}_{\Gamma_\br} \in \mathscr{L}({H}^{-\frac{1}{2}}(\Gamma_\br), {H}^{\frac{1}{2}}(\Gamma_\br)).
	\end{equation}
    \item {\bf Double Layer BIO:}
    \begin{equation}\label{eq:def_double_layer}
    \mathsf{K}^{(\kappa)}_{\Gamma_\br} \coloneqq \gamma_{\Gamma_\br} \mathsf{D}^{(\kappa)}_{\Gamma_\br}+\frac{1}{2} \mathsf{I}_{\Gamma_\br} \in \mathscr{L}({H}^{\frac{1}{2}}(\Gamma_\br), {H}^{\frac{1}{2}}(\Gamma_\br)),
    \end{equation}
    where $\mathsf{I}_{\Gamma_\br}$ denotes the identity operator on $H^{\pm\frac{1}{2}}(\Gamma_\br)$. 
    \item 
    {\bf Adjoint Double Layer BIO:} 
    \begin{equation}\label{eq:def_adjoint_double_layer}
    \mathsf{K}^{(\kappa)\prime}_{\Gamma_\br} \coloneqq \gamma_{1,\Gamma_\br} \mathsf{S}^{(\kappa)}_{\Gamma_\br}-\frac{1}{2} \mathsf{I}_{\Gamma_\br}\in \mathscr{L}({H}^{-\frac{1}{2}}(\Gamma_\br), {H}^{-\frac{1}{2}}(\Gamma_\br)).
    \end{equation}
    \item 
    {\bf Hypersingular BIO:}
    \begin{equation}\label{eq:def_hypersingular}
    \mathsf{W}^{(\kappa)}_{\Gamma_\br} \coloneqq -\gamma_{1,\Gamma_\br} \mathsf{D}^{(\kappa)}_{\Gamma_\br}\in \mathscr{L}({H}^{\frac{1}{2}}(\Gamma_\br), {H}^{-\frac{1}{2}}(\Gamma_\br)).
    \end{equation}
\end{itemize}

The single layer potential has continuous Dirichlet trace and Neumann
jump equal to the negative of its density; the double layer potential
has Dirichlet jump equal to its density and continuous Neumann trace,
cf.~\cite[Chapter~6]{mclean2000strongly}.

\section{Shape Holomorphy of the Calder\'{o}n Projector}
\label{sec:BiosResults}
\label{sec:main_results}

In this section, we prove the main result of this work,
i.e., we establish wavenumber-explicit shape holomorphy of
the Calder\'{o}n projector for the Helmholtz operator, which is stated in
Theorem~\ref{thm:Calderon_Projector}.
Following \cite{HDPS26Laplace}, our analysis is based on the
abstract operator construction of the four BIOs in \eqref{eq:def_single_layer},
\eqref{eq:def_double_layer}, \eqref{eq:def_adjoint_double_layer}, and
\eqref{eq:def_hypersingular}, rather than on their explicit
representations, as in \cite{henriquez2021shape,HS21,DL23,PHJ23,dolz2024parametric}.

In the following, we consider $d\in\{2,3\}$ and fix a wavenumber $\kappa>0$.
For each $\br\in\mathfrak T$, we define the Calder\'on projector on
$\Gamma_\br$ as
\begin{equation*}
    \mathsf{C}^{(\kappa)}_{\Gamma_\br}
    \coloneqq
    \begin{pmatrix}
        -\mathsf{K}^{(\kappa)}_{\Gamma_\br}
        +\frac{1}{2}\mathsf{I}_{\Gamma_\br}
        & \mathsf{V}^{(\kappa)}_{\Gamma_\br}
        \\
        \mathsf{W}^{(\kappa)}_{\Gamma_\br}
        & \mathsf{K}^{(\kappa)\prime}_{\Gamma_\br}
        +\frac{1}{2}\mathsf{I}_{\Gamma_\br}
    \end{pmatrix}.
\end{equation*}
For $\boldsymbol v=(v_D,v_N)^\top\in
H^{\frac{1}{2}}(\Gamma_\br)\times H^{-\frac{1}{2}}(\Gamma_\br)$,
using the BIO definitions and the jump relations from
Subsection~\ref{sec:bios_Lipschitz}, we obtain
\begin{equation*}
\boldsymbol\gamma
(\mathsf S^{(\kappa)}_{\Gamma_\br}v_N
-\mathsf D^{(\kappa)}_{\Gamma_\br}v_D)
=\mathsf C^{(\kappa)}_{\Gamma_\br}\boldsymbol v,
\quad
\boldsymbol\gamma^c
(\mathsf D^{(\kappa)}_{\Gamma_\br}v_D
-\mathsf S^{(\kappa)}_{\Gamma_\br}v_N)
=(\mathbb I_{\Gamma_\br}-\mathsf C^{(\kappa)}_{\Gamma_\br})
\boldsymbol v,
\end{equation*}
where $\mathbb I_{\Gamma_\br}$ denotes the identity on
$H^{\frac{1}{2}}(\Gamma_\br)\times H^{-\frac{1}{2}}(\Gamma_\br)$.
For each $\br \in \mathfrak{T}$, $\mathsf{C}^{(\kappa)}_{\Gamma_\br}$ is a bounded projector on
$H^{\frac{1}{2}}(\Gamma_\br)\times H^{-\frac{1}{2}}(\Gamma_\br)$.
Indeed, the representation formula and the jump relations imply
$\left(\mathsf{C}^{(\kappa)}_{\Gamma_\br}\right)^{2}=\mathsf{C}^{(\kappa)}_{\Gamma_\br}$.
Its range consists of the Cauchy data of solutions of
$\mathcal L^{(\kappa)}u=0$ in $\D_\br$.
For more details, we refer to \cite[Section~3.6, Proposition~3.6.2]{SS10}.

Using the pullback operator $\tau_\br$ introduced in
Subsection~\ref{sec:Pullbck}, and recalling the mapping properties stated in
Proposition~\ref{prop:iso_pullback}, we define the four transported BIOs
to the reference boundary $\widehat{\Gamma}$ as follows:
\begin{equation}\label{eq:transported_calderon_bios}
\begin{aligned}
    \widehat{\mathsf{V}}^{(\kappa)}_\br
    &\coloneqq \tau_\br\mathsf{V}^{(\kappa)}_{\Gamma_\br}\tau_\br^{-1}
    \in\mathscr{L}\left(H^{-\frac{1}{2}}(\widehat{\Gamma}),
    H^{\frac{1}{2}}(\widehat{\Gamma})\right),
    \\
    \widehat{\mathsf{K}}^{(\kappa)}_\br
    &\coloneqq \tau_\br\mathsf{K}^{(\kappa)}_{\Gamma_\br}\tau_\br^{-1}
    \in\mathscr{L}\left(H^{\frac{1}{2}}(\widehat{\Gamma}),
    H^{\frac{1}{2}}(\widehat{\Gamma})\right),
    \\
    \widehat{\mathsf{K}}^{(\kappa)\prime}_\br
    &\coloneqq \tau_\br\mathsf{K}^{(\kappa)\prime}_{\Gamma_\br}\tau_\br^{-1}
    \in\mathscr{L}\left(H^{-\frac{1}{2}}(\widehat{\Gamma}),
    H^{-\frac{1}{2}}(\widehat{\Gamma})\right),
    \\
    \widehat{\mathsf{W}}^{(\kappa)}_\br
    &\coloneqq \tau_\br\mathsf{W}^{(\kappa)}_{\Gamma_\br}\tau_\br^{-1}
    \in\mathscr{L}\left(H^{\frac{1}{2}}(\widehat{\Gamma}),
    H^{-\frac{1}{2}}(\widehat{\Gamma})\right).
\end{aligned}
\end{equation}
For each $\br \in \mathfrak{T}$, the transported Calder\'on projector is
\begin{equation*}
    \widehat{\mathsf{C}}^{(\kappa)}_\br
    \coloneqq
    \begin{pmatrix}
        -\widehat{\mathsf{K}}^{(\kappa)}_\br
        +\frac{1}{2}\mathsf{I}_{\widehat{\Gamma}}
        & \widehat{\mathsf{V}}^{(\kappa)}_\br
        \\
        \widehat{\mathsf{W}}^{(\kappa)}_\br
        & \widehat{\mathsf{K}}^{(\kappa)\prime}_\br
        +\frac{1}{2}\mathsf{I}_{\widehat{\Gamma}}
    \end{pmatrix}.
\end{equation*}
Observe that $\left(\widehat{\mathsf{C}}^{(\kappa)}_\br\right)^{2}=\widehat{\mathsf{C}}^{(\kappa)}_\br$
for each $\br \in \mathfrak{T}$.

Assumption~\ref{assump:domain_T_2}, introduced ahead in
Subsection~\ref{sec:shape_hol_Newton}, refers to the existence of
extensions of elements of $\mathfrak T$ to a common hold-all domain.
This is required to prove shape holomorphy of the Newton potential
(Theorem~\ref{thm:shape_hol_Newton}).
The following result follows from Theorems~\ref{thm:Vhlmrp}, \ref{thm:KHlm}, \ref{thm:KpHlm}, and~\ref{thm:WHlm}.

\begin{theorem}\label{thm:Calderon_Projector}
Let $\mathfrak T$ be a set of domain transformations
satisfying Assumptions~\ref{assump:transformations_T}
and~\ref{assump:domain_T_2}.
Then there exists $\delta_0>0$, depending only on $\mathfrak T$, such that,
for every $\kappa>0$ and $0<\delta<\delta_0$, the map
\begin{equation*}\label{eq:map_hol_C}
	\mathcal{C}^{(\kappa)}:\mathfrak{T}
	\rightarrow
	\mathscr{L}\left(
	H^{\frac{1}{2}}(\widehat{\Gamma})\times
	H^{-\frac{1}{2}}(\widehat{\Gamma}),
	H^{\frac{1}{2}}(\widehat{\Gamma})\times
	H^{-\frac{1}{2}}(\widehat{\Gamma})
	\right):\br
    \mapsto \widehat{\mathsf{C}}^{(\kappa)}_\br
\end{equation*}
admits a bounded holomorphic extension to $\mathfrak T_\delta$.
For each $0<R<\delta_0$ and $\kappa_0>0$, there exist
$C_1=C_1(\mathfrak T,\delta_0,\kappa_0)>0$ and $C_2=C_2(\mathfrak T,\delta_0)>0$ such that, for all
$\kappa\geq\kappa_0$ and $0<\delta<\min\{R,\delta_0-R\}$,
\begin{equation*}
\sup_{\br\in\mathfrak T_\delta}
\norm{\widehat{\mathsf C}^{(\kappa)}_\br}_{\mathscr L(H^{\frac{1}{2}}(\widehat\Gamma)\times H^{-\frac{1}{2}}(\widehat\Gamma))}
\leq C_1(1+\kappa)^{\frac{d+1}{2}}
e^{\frac{4\delta}{\pi R}\left(\kappa\delta_0C_2+2\log(1+\kappa)\right)}.
\end{equation*}
\end{theorem}

\subsection{Shape Holomorphy of the Newton Potential}\label{sec:shape_hol_Newton}
In this subsection, we prove shape holomorphy of the Newton potential
introduced in \eqref{eq:newtonpot}. More precisely, for each $\br\in\mathfrak T$, we define
\begin{equation*}
	\widehat{\mathsf{N}}^{(\kappa)}_{\br}
	\coloneqq
	\tau_{\br}
	\mathsf{N}^{(\kappa)}_{\br}
	\tau^{-1}_{\br}
	\in 
	\mathscr{L}
	\left(	
		\widetilde{H}^{-1}(\widehat{\text{D}})
		,
		H^{1}(\widehat{\text{D}})
	\right),
\end{equation*}
and show shape holomorphy of the \emph{shape-to-operator} map
\begin{equation*}
	\mathfrak{T}
	\ni
	\br 
	\mapsto 
	\widehat{\mathsf{N}}^{(\kappa)}_{\br}
	\in
	\mathscr{L}
	\left(	
		\widetilde{H}^{-1}(\widehat{\text{D}})
		,
		H^{1}(\widehat{\text{D}})
	\right).
\end{equation*}
For $\widehat u\in L^2(\widehat\D)$ and for each $\br\in\mathfrak T$, this operator admits the explicit representation
\begin{equation}
\label{eq:Nhatexplicit}
	\left(
		\widehat{\mathsf{N}}^{(\kappa)}_{\br}
		\widehat{u}
	\right)
	(\widehat\x)
	=
	\int_{\widehat{\D}}
		\mathsf{g}^{(\kappa)}_{\br}
	\left(
		\widehat\x,\widehat\y
	\right)
	\widehat{u}(\widehat\y)
	\dd
	\widehat\y,
	\quad
	\widehat\x\in\widehat{\D},
\end{equation}
where
\begin{equation*}
	\mathsf{g}^{(\kappa)}_{\br}
	\left(
		\widehat\x,\widehat\y
	\right)
	\coloneqq
	G^{(\kappa)}(\br(\widehat\x)-\br(\widehat\y))
	J_{\br}(\widehat\y),
	\quad
	(\widehat\x,\widehat\y)
	\in
	(\widehat{\D} \times \widehat{\D})^\star,
\end{equation*}
with $G^{(\kappa)}$ as in \eqref{eq:green_G_omega}, and
$(\widehat\D\times\widehat\D)^\star
=\{(\widehat\x,\widehat\y)\in\widehat\D\times\widehat\D:
\widehat\x\ne\widehat\y\}$.

As previously explained, in addition to Assumption~\ref{assump:transformations_T}, our analysis of the Newton potential requires the following assumption that refers to the
existence of extensions of each element of $\mathfrak T$ to one common
hold-all domain.

\begin{assumption}\label{assump:domain_T_2}
There exists a bounded Lipschitz domain $\mathcal{O} \subset \mathbb{R}^d$, $d \in \{2,3\}$,
such that:
\begin{itemize}
    \item[\normalfont{(i)}]
    $\overline{\widehat{\D}},\overline{\D_\br} \subset \mathcal{O}$ for each
    $\br \in \mathfrak{T}$.
    \item[\normalfont{(ii)}]
    Each $\br\in \mathfrak{T}$ admits an extension, also denoted by
    $\br$, which is a bi-Lipschitz bijection from $\mathcal{O}$ onto itself
    and whose restriction to $\widehat{\D}$ coincides with the original
    transformation. The collection of these extensions, denoted by
    $\mathfrak{T}_{\mathcal{O}}$, is a compact subset of
    $W^{1,\infty}(\mathcal{O};\mathbb{R}^d)$.
    \item[\normalfont{(iii)}]
    The extensions from item (ii) become the identity in a vicinity of
    $\partial\mathcal{O}$, i.e., there exists $\varepsilon>0$ such that, for
    each $\br\in\mathfrak{T}_\mathcal{O}$, it holds that
\begin{equation*}
        \br(\x)=\x
        \quad
        \text{for each }
        \x\in\mathcal{O}\cap
        \bigcup_{\y \in \partial\mathcal{O}}B(\y,\varepsilon).
    \end{equation*}
\end{itemize}
\end{assumption}

For each $\br \in \mathfrak{T}_{\mathcal{O}}$, we define the
pullback operator in $\mathcal{O}$ as follows
\begin{equation}\label{eq:pullback_O}
    \tau_{\br,\mathcal{O}} u
    \coloneqq
    \left\{
    \begin{array}{cl}
        u\circ \br, & u\in H^1(\mathcal{O})\ \text{or}\ u\in L^2(\mathcal{O}), \\
        \mathcal{R}_{\mathcal{O}}
        \tau_{\br,\mathcal{O}}
        \mathcal{R}^{-1}_{\mathcal{O}} u, & u\in \widetilde{H}^{-1}(\mathcal{O}),
    \end{array}
    \right.
\end{equation}
with $\mathcal{R}_{\mathcal{O}}: H^1(\mathcal{O}) \rightarrow \widetilde{H}^{-1}(\mathcal{O})$
being the corresponding Riesz map.

Working under Assumption~\ref{assump:domain_T_2}
in addition to Assumption~\ref{assump:transformations_T}, for each
$\br \in \mathfrak{T}_\mathcal{O}$, we define the \emph{bilinear} form
$\dual{\cdot}{\cdot}_{\widetilde{H}^{-1}(\mathcal{O}) \times H^1(\mathcal{O}),\br}:
\widetilde{H}^{-1}(\mathcal{O}) \times H^1(\mathcal{O}) \rightarrow \mathbb{C}$
as follows:
\begin{equation}\label{eq:bilnfrmH1_2_O}
\begin{aligned}
  	\dual{
		\widehat{u}
	}{
		\widehat{v}
	}_{\widetilde{H}^{-1}(\mathcal{O}) \times H^1(\mathcal{O}),\br}
    \coloneqq
    &
    \int_{\mathcal{O}}
    \left(D\br\right)^{-\top}
    \widehat{\nabla} ({{\mathcal{R}}}^{-1}_\mathcal{O}\widehat{u})
    \cdot
    \left(D\br\right)^{-\top}
    \widehat{\nabla}\widehat{v}
    J_{\br}
    \mathrm{d}\widehat{\x}
    \\
    &
    +
    \int_{\mathcal{O}}
    ({{\mathcal{R}}}^{-1}_\mathcal{O}\widehat{u})
    \widehat{v}
    J_{\br}
    \mathrm{d}\widehat{\x},
    \quad
    \forall \widehat{u} \in \widetilde{H}^{-1}(\mathcal{O}), \; \forall \widehat{v}\in H^1(\mathcal{O}).
\end{aligned}
\end{equation}

The results introduced in Subsection~\ref{sec:Pullbck}
also apply to $\mathcal O$, $\partial\mathcal O$, and
$\mathfrak T_{\mathcal O}$ in place of $\widehat\D$,
$\widehat\Gamma$, and $\mathfrak T$, respectively.
In particular, for each $\br \in \mathfrak{T}_\mathcal{O}$ it holds that
for all $u \in \widetilde{H}^{-1}(\mathcal{O})$ and all $v\in H^1(\mathcal{O})$
\begin{equation*}\label{eq:transport_domainO}
	\dual{u}{{v}}_{\widetilde{H}^{-1}(\mathcal{O}) \times H^1(\mathcal{O})}
	=
    	\dual{
    		\tau_{\br,\mathcal{O}} {u}
	}{
		\tau_{\br,\mathcal{O}} {v}
	}_{\widetilde{H}^{-1}(\mathcal{O}) \times H^1(\mathcal{O}),\br},
\end{equation*}
where $\tau_{\br,\mathcal{O}}$ is defined in \eqref{eq:pullback_O}.

Given $\varepsilon>0$, we define
\begin{equation*}\label{eq:TepsO}
	\mathfrak{T}_{\varepsilon,\mathcal{O}}
	\coloneqq
	\left\{
		\br 
		\in 
		W^{1,\infty}(\mathcal{O};\IC^d)
		:
		\exists\,
		\widetilde{\br} \in \mathfrak{T}_\mathcal{O}
		\text{ such that }
		\norm{\br - \widetilde{\br}}_{W^{1,\infty}(\mathcal{O};\IC^d)}
		<
		\varepsilon
	\right\}.
\end{equation*}
Applying Lemma~\ref{lem:bilinear_ext_D} to $\mathcal O$ and
$\mathfrak T_{\mathcal O}$ gives a bounded holomorphic extension of
the bilinear form in \eqref{eq:bilnfrmH1_2_O} to
$\mathfrak T_{\delta,\mathcal O}$ for some $\delta>0$.%

For each $\varphi\in L^2(\mathcal{O})$, we define the
Newton potential in $\mathcal{O}$ by
\begin{equation}\label{eq:newton_potential_O}
    \big(\mathsf{N}^{(\kappa)}_{\mathcal{O}}\varphi\big)(\x)
    \coloneqq
    \int_{\mathcal{O}}G^{(\kappa)}(\x-\y)\varphi(\y)\,\mathrm{d}\y,
    \quad \x\in\mathcal{O}.
\end{equation}

We first establish shape holomorphy for the Newton potential acting
from $L^2(\widehat\D)$ to $H^1(\widehat\D)$, and then consider
its extension to $\widetilde H^{-1}(\widehat\D)$ in
Theorem~\ref{thm:shape_hol_Newton}.
The complete proof of the following result is given in Appendix~\ref{app:proof_shape_hol_Newton_L2}.

\begin{theorem}\label{thm:shape_hol_Newton_L2}
Let $\mathfrak T$ be a set of domain transformations satisfying
Assumption~\ref{assump:transformations_T}.
Then there exists $\delta_0 = \delta_0(\mathfrak{T})>0$, depending only on $\mathfrak T$, such that,
for every $\kappa>0$ and $0<\delta<\delta_0$, the map
\begin{equation*}\label{eq:map_hol_N}
\mathcal N^{(\kappa)}:\mathfrak T
\rightarrow\mathscr L(L^2(\widehat\D),H^1(\widehat\D)):\br
\mapsto\widehat{\mathsf N}^{(\kappa)}_\br
\end{equation*}
admits a bounded holomorphic extension to $\mathfrak T_\delta$.
For each $0<R<\delta_0$ and $\kappa_0>0$, there exists
$C_1=C_1(\mathfrak T,\delta_0,\kappa_0)>0$ such that, for all
$\kappa\geq\kappa_0$ and $0<\delta<\min\{R,\delta_0-R\}$,
\begin{equation}\label{eq:newton_reference_improved_complex_bound}
\sup_{\br\in\mathfrak T_\delta}
\norm{\widehat{\mathsf N}^{(\kappa)}_\br}_{
\mathscr L(L^2(\widehat\D),H^1(\widehat\D))}
\leq C_1
e^{\frac{4\delta}{\pi R}
\left(\kappa\delta_0C_2
+\frac{d-1}{2}\log(1+\kappa)\right)}.
\end{equation}
Here, $C_2=C_2(\mathfrak T,\delta_0)$, where for $0<\delta\leq\delta_0$ we set
\begin{equation}\label{eq:kernel_exponential_constant}
C_2(\mathfrak T,\delta)\coloneqq
\frac{1}{\sqrt{\eta}}
\sup_{\br\in\mathfrak T_\delta}
\sup_{\widetilde\br\in\mathfrak T}
\left(\operatorname{Lip}(\br;\widehat\D)
+\operatorname{Lip}(\widetilde\br;\widehat\D)\right)<\infty,
\end{equation}
where $\eta>0$ is as in Lemma~\ref{lem:lipschitz_preservation}.
\end{theorem}

\begin{theorem}\label{thm:shape_hol_Newton}
Let $\mathfrak T$ be a set of domain transformations satisfying
Assumptions~\ref{assump:transformations_T} and~\ref{assump:domain_T_2}.
Then there exists $\delta_0=\delta_0(\mathfrak T)>0$ such that,
for every $\kappa>0$ and $0<\delta<\delta_0$, the map
\begin{equation*}\label{eq:map_hol_N_H_-1_H1}
\mathcal N^{(\kappa)}:\mathfrak T
\rightarrow\mathscr L(\widetilde H^{-1}(\widehat\D),H^1(\widehat\D)):\br
\mapsto\widehat{\mathsf N}^{(\kappa)}_\br
\end{equation*}
admits a bounded holomorphic extension to $\mathfrak T_\delta$.
Moreover, there exists $C_2=C_2(\mathfrak T,\delta_0)>0$ such that, for every $\kappa_0>0$, there exists $C_1=C_1(\mathfrak T,\delta_0,\kappa_0)>0$
satisfying
\begin{equation}\label{eq:Newton_explicit_refined_frequency}
\sup_{\br\in\mathfrak T_\delta}
\norm{\widehat{\mathsf N}^{(\kappa)}_\br}_
{\mathscr L(\widetilde H^{-1}(\widehat\D),H^1(\widehat\D))}
\leq C_1(1+\kappa)
e^{\frac{4\delta}{\pi R}\left(\kappa\delta_0C_2+\frac{d+1}{2}\log(1+\kappa)\right)}
\end{equation}
for all $\kappa\geq\kappa_0$, $0<R<\delta_0$, and
$0<\delta<\min\{R,\delta_0-R\}$.
\end{theorem}
\begin{proof}
We proceed in four steps as follows.

\textbf{Step 1: Representation using weak formulation.}
In this step we cover the first three steps of the proof of \cite[Theorem~3.3]{HDPS26Laplace}.
We make the computations explicit to keep track of additional terms involving $\kappa$.
In the following, we do not distinguish in terms of notation between functions
defined on $\mathcal O$ and their restriction to $\widehat\D$.
Conversely, densities supported in $\widehat\D$ or $\D_\br$ considered on
$\mathcal O$ are identified with their extensions by zero.

For each $\br\in\mathfrak T_{\mathcal O}$, consider $\widehat{\varphi}\in \mathscr{C}^\infty_0(\widehat{\D})$ 
and define $\ell_{\widehat{\varphi},\br}\in\widetilde H^{-1}(\D_\br)$ as
\begin{equation*}
    \dual{\ell_{\widehat\varphi,\br}}{v}_{\widetilde H^{-1}(\D_\br)\times H^1(\D_\br)}
    =
    \int_{\D_\br}
    \tau^{-1}_\br\widehat{\varphi} v\,\dd\x,
    \quad\forall v\in H^1(\D_\br).
\end{equation*}
Set $\widehat\ell_{\widehat{\varphi},\br}\coloneqq\tau_\br\ell_{\widehat{\varphi},\br}$, thus
\begin{equation*}
    \dual{\widehat\ell_{\widehat{\varphi},\br}}{\widehat v}_
    {\widetilde H^{-1}(\widehat\D)\times H^1(\widehat\D),\br}
    =
    \int_{\widehat\D}\widehat\varphi\widehat vJ_\br\,\dd\widehat\x,
    \quad\forall\widehat v\in H^1(\widehat\D).
\end{equation*}
Moreover, Lemma~\ref{lem:bilinear_ext_D}, item~\normalfont{(i)}, yields
\begin{equation}\label{eq:bound_dual_J_r_new}
    \snorm{
        \int_{\widehat\D}\widehat\varphi\widehat vJ_\br\,\dd\widehat\x
    }
    \leq
    C_1(\mathfrak{T},\delta)
    \norm{\widehat\ell_{\widehat{\varphi},\br}}_{\widetilde H^{-1}(\widehat\D)}
    \norm{\widehat{v}}_{H^1(\widehat\D)},
    \quad\forall\widehat v\in H^1(\widehat\D),
\end{equation}
with $C_1(\mathfrak{T},\delta)$ as in \eqref{eq:C_1_domain}. 
The definition of the Newton potential by duality, cf.~\cite[Section~6.1]{STEI08}, yields
\begin{equation*}
    \bigl(\mathsf N^{(\kappa)}_\br\ell_{\widehat{\varphi},\br}\bigr)(\x)
    =
    \int_{\D_\br}G^{(\kappa)}(\x-\y) (\tau^{-1}_\br\widehat\varphi)(\y)\,\dd\y
\quad\text{for almost every }\x\in\D_\br.
\end{equation*}
Therefore, for each $\br \in \mathfrak{T}$, we obtain, for $\widehat\x\in\widehat\D$,
\begin{equation}\label{eq:defExt}
\begin{aligned}
    \bigl(\widehat{\mathsf N}^{(\kappa)}_\br\widehat\ell_{\widehat\varphi,\br} \bigr)(\widehat\x)
    &=\int_{\D_\br}G^{(\kappa)}(\br(\widehat\x)-\y)
    (\tau^{-1}_\br\widehat\varphi)(\y)\,\dd\y\\
    &=\int_{\widehat\D}G^{(\kappa)}(\br(\widehat\x)-\br(\widehat\y))
    \widehat\varphi(\widehat\y)J_\br(\widehat\y)\,\dd\widehat\y \\
    & = 
    \bigl(\widehat{\mathsf{N}}^{(\kappa)}_{\br,\mathcal O} \widehat\varphi\bigr)(\widehat\x),
\end{aligned}
\end{equation}
where $\widehat{\mathsf{N}}^{(\kappa)}_{\br,\mathcal O}
\coloneqq\tau_{\br,\mathcal O}\mathsf N^{(\kappa)}_{\mathcal O}\tau^{-1}_{\br,\mathcal O}$,
and $\widehat\varphi$ is extended by zero to $\mathcal O$.

For each $\br\in\mathfrak T_{\mathcal O}$, we set
$u_\br\coloneqq\mathsf N^{(\kappa)}_{\mathcal O}
\tau_\br^{-1}\widehat\varphi\in H^1(\mathcal O)$.
As in Step~1 of the proof of \cite[Theorem~3.3]{HDPS26Laplace},
integration by parts yields, for all $v \in H^1(\mathcal{O})$,
\begin{equation}\label{eq:integrationbyparts}
\begin{aligned}
   \int_{\mathcal{O}}
   \left(
   \nabla  u_\br
   \cdot
   \nabla  {v}
   -\kappa^2u_\br v
   \right)
   \mathrm{d}{\x}
   =
   \int_{\D_\br}
   \tau^{-1}_{\br}\widehat{\varphi}
   {v}
   \mathrm{d}{\x}
   +
   \int_{\partial \mathcal{O}}
   (\boldsymbol n_{\mathcal O} \cdot \nabla  u_\br) \gamma_{\partial\mathcal{O}} {v}
   \mathrm{ds}_{\x}.
\end{aligned}
\end{equation}
Here, $\boldsymbol n_{\mathcal O}$ signifies the outward unit normal to $\mathcal O$,
and $\gamma_{\partial\mathcal O}$ denotes the Dirichlet trace operator onto $\partial\mathcal O$.
Recalling the explicit representation of the Newton potential
given in \eqref{eq:newton_potential_O}, and using
the Fubini--Tonelli theorem, we obtain
\begin{equation}\label{eq:reformulatetoB}
\int_{\partial\mathcal O}
(\boldsymbol n_{\mathcal O}\cdot\nabla u_\br)\gamma_{\partial\mathcal O}v
\,\dd s_\x
=\int_{\widehat\D}\widehat\varphi
\bigl(\mathsf B^{(\kappa)}_\br\gamma_{\partial\mathcal O}v\bigr)
J_\br\,\dd\widehat\y,
\end{equation}
where $\mathsf B^{(\kappa)}_\br:H^{\frac{1}{2}}(\partial\mathcal O)
\rightarrow H^1(\widehat\D)$ is defined, for
$\br\in\mathfrak T_{\mathcal O}$, as
\begin{equation}\label{eq:B_operator}
    \left(\mathsf B_\br^{(\kappa)}\phi\right)(\widehat\y)
    \coloneqq
    \int_{\partial\mathcal O}
    \bigl(\boldsymbol n_{\mathcal O}(\x)\cdot\nabla_{\x}
    G^{(\kappa)}(\x-\br(\widehat\y))\bigr)
    \phi(\x)\,\dd s_{\x},
    \quad\widehat\y\in\widehat\D.
\end{equation}

Inserting \eqref{eq:reformulatetoB} into \eqref{eq:integrationbyparts}
and transporting to $\mathcal O$ and $\widehat\D$ as in Step~2 of the proof of
\cite[Theorem~3.3]{HDPS26Laplace},
one has that, for each $\br \in \mathfrak{T}_{\mathcal{O}}$ and all $\widehat{v} \in H^1(\mathcal{O})$,
\begin{equation}\label{eq:pulledpackpartialintegration2}
\begin{aligned}
	\int_{\mathcal{O}}
	\left(
	(D\br)^{-\top}
	\widehat{\nabla} \widehat{u}_\br
	\cdot
	(D\br)^{-\top}
	\widehat{\nabla}  {\widehat{v}}
	-\kappa^2\widehat{u}_\br\widehat{v}
	\right)
	J_{\br}
	\mathrm{d}{\widehat{\x}}
    =
    &
	\int_{\widehat{\D}}
	\widehat{\varphi}
	{\widehat{v}}
	J_{\br}
	\mathrm{d}\widehat{\x}
	\\
	&
	+
	\int_{\widehat{\D}}
	\widehat{\varphi}
	\left(\mathsf{B}^{(\kappa)}_\br \gamma_{\partial\mathcal{O}} {\widehat{v}} \right)
	J_{\br}
	\dd \widehat{\y},
\end{aligned}
\end{equation}
where we have set $\widehat{u}_\br
\coloneqq  \tau_{\br,\mathcal{O}} u_\br = \tau_{\br,\mathcal{O}} \mathsf{N}^{(\kappa)}_{\mathcal{O}} \tau^{-1}_\br\widehat{\varphi}\in H^1(\mathcal{O})$ and we have used that $\gamma_{\partial\mathcal{O}} {v} = 
\gamma_{\partial\mathcal{O}} \widehat{v}$ for all $\widehat{v} \in H^1(\mathcal{O})$,
with $v=\tau^{-1}_{\br,\mathcal O}\widehat v$, due to item (iii) in
Assumption~\ref{assump:domain_T_2}.

For $\br\in\mathfrak T_{\mathcal O}$, consider the maps
\begin{equation*}
\begin{aligned}
\mathcal M_1 :&\mathfrak T_{\mathcal O}
\rightarrow\mathscr L(H^1(\mathcal O),\IC):
\br\mapsto\left(\widehat v\mapsto\int_{\mathcal O}
(D\br)^{-\top}\widehat\nabla\widehat u_\br
\cdot(D\br)^{-\top}\widehat\nabla\widehat v
J_\br\,\dd\widehat\x\right),\\
\mathcal M_2 :&\mathfrak T_{\mathcal O}
\rightarrow\mathscr L(H^1(\mathcal O),\IC):
\br\mapsto\left(\widehat v\mapsto\kappa^2\int_{\mathcal O}
\widehat u_\br\widehat vJ_\br\,\dd\widehat\x\right),\\
\mathcal M_3 :&\mathfrak T_{\mathcal O}
\rightarrow\mathscr L(H^1(\mathcal O),\IC):
\br\mapsto\left(\widehat v\mapsto\int_{\widehat\D}
\widehat\varphi\widehat vJ_\br\,\dd\widehat\x\right),\\
\mathcal M_4 :&\mathfrak T_{\mathcal O}
\rightarrow\mathscr L(H^1(\mathcal O),\IC):
\br\mapsto\left(\widehat v\mapsto\int_{\widehat\D}\widehat\varphi
\left(\mathsf B^{(\kappa)}_\br
\gamma_{\partial\mathcal O}\widehat v\right)
J_\br\,\dd\widehat\y\right).
\end{aligned}
\end{equation*}
It follows from Lemma~\ref{lem:shape_hol_bilinear_forms} and
Theorem~\ref{thm:shape_hol_Newton_L2}, both applied to $\mathcal O$,
that the maps $\mathcal M_1$ and $\mathcal M_2$ admit bounded holomorphic
extensions to $\mathfrak T_{\delta,\mathcal O}$, for some $\delta>0$.
Note that for $\mathcal M_2$ the domain of holomorphy is independent of $\kappa$.
It follows from Lemma~\ref{lem:shape_hol_bilinear_forms} and the
continuity of the restriction from $H^1(\mathcal O)$ to
$L^2(\widehat\D)$ that $\mathcal M_3$ admits a bounded holomorphic extension to
$\mathfrak T_{\delta,\mathcal O}$, for a possibly different $\delta>0$.
Finally, Proposition~\ref{prop:shape_hol_B_r_O} provides a bounded
holomorphic extension of the map
$
\mathfrak T_{\mathcal O}\ni\br\mapsto\mathsf B^{(\kappa)}_\br
\in\mathscr L(H^{\frac{1}{2}}(\partial\mathcal O),H^1(\widehat\D))
$
to $\mathfrak T_{\delta,\mathcal O}$, implying that $\mathcal M_4$ admits a bounded holomorphic
extension to $\mathfrak T_{\delta,\mathcal O}$, for some $\delta>0$.
In what follows, we choose $\delta>0$ sufficiently small and independent of $\kappa$ such that $\mathcal{M}_1$, $\mathcal{M}_2$, $\mathcal{M}_3$ and $\mathcal{M}_4$ admit bounded holomorphic extensions to $\mathfrak{T}_{\delta,\mathcal{O}}$ (still denoted by the same symbols).

Next, we prove that the identity
$\mathcal M_1(\br)-\mathcal M_2(\br)-\mathcal M_3(\br)-\mathcal M_4(\br)=0$
holds for each real-valued domain transformation
$\br\in V_{\widetilde\br}\subset\mathfrak T_{\delta,\mathcal O}$,
where, for a fixed $\widetilde\br\in\mathfrak T_{\mathcal O}$, we define
\begin{equation*}
U_{\widetilde\br}\coloneqq
B_{W^{1,\infty}(\mathcal O;\IC^d)}(\widetilde\br,\delta),
\quad
V_{\widetilde\br}\coloneqq
U_{\widetilde\br}\cap W^{1,\infty}(\mathcal O;\IR^d).
\end{equation*}
Let $\delta>0$ be sufficiently small, independently of $\kappa$,
so that \cite[Lemma~B.1, item~(ii)]{HDPS26Laplace} and
Proposition~\ref{prop:shape_hol_B_r_O} apply.
The set $U_{\widetilde\br}$ is connected and open in
$W^{1,\infty}(\mathcal O;\IC^d)$, whereas $V_{\widetilde\br}$ is
nonempty and relatively open in $W^{1,\infty}(\mathcal O;\IR^d)$.
For each $\br\in V_{\widetilde\br}$, set
\begin{equation*}
    \mathcal{O}_\br
    \coloneqq
    \br(\mathcal{O})
    \quad
    \text{and}
    \quad
    \D_\br
    \coloneqq
    \br(\widehat{\D}).
\end{equation*}
By \cite[Lemma~B.1, item~(ii)]{HDPS26Laplace}, $\br:\mathcal{O}\rightarrow\mathcal{O}_\br$ is bi-Lipschitz and $\mathcal{O}_\br$ is a bounded Lipschitz domain.
Moreover, $\overline{\D_\br}$ is compactly contained in $\mathcal{O}_\br$, since $\overline{\widehat{\D}}$ is compactly contained in $\mathcal{O}$.
Integration by parts on $\mathcal O_\br$, change of variables,
and \eqref{eq:B_operator_real_O_r}, as in Step~3 of
the proof of \cite[Theorem~3.3]{HDPS26Laplace}, give
\begin{equation*}
\mathcal M_1(\br)-\mathcal M_2(\br)-\mathcal M_3(\br)-\mathcal M_4(\br)=0
\end{equation*}
for each $\br\in V_{\widetilde\br}$.

Consequently,
$\mathcal M\coloneqq\mathcal M_1-\mathcal M_2-\mathcal M_3-\mathcal M_4$
is holomorphic on $U_{\widetilde\br}$ and vanishes on
$V_{\widetilde\br}$.
It follows from \cite[Theorem~2.2]{HDPS26Laplace} that
$\mathcal M(\br)=0$ for each $\br\in U_{\widetilde\br}$.
Finally, recalling that
$\mathfrak T_{\delta,\mathcal O}
=\bigcup_{\widetilde\br\in\mathfrak T_{\mathcal O}}
U_{\widetilde\br}$, we conclude that
\eqref{eq:pulledpackpartialintegration2} holds for each
$\br\in\mathfrak T_{\delta,\mathcal O}$.

\begin{samepage}
\textbf{Step 2: Continuity of $\widehat{\mathsf N}^{(\kappa)}_\br:
\widetilde H^{-1}(\widehat\D)\rightarrow H^1(\widehat\D)$.}
For $\br\in\mathfrak T_{\delta,\mathcal O}$, let
$\widehat\ell_{\widehat\varphi,\br}\in\widetilde H^{-1}(\widehat\D)$
be the unique functional satisfying
\begin{equation*}
\dual{\widehat\ell_{\widehat\varphi,\br}}{\widehat v}_
{\widetilde H^{-1}(\widehat\D)\times H^1(\widehat\D),\br}
=
\int_{\widehat\D}\widehat\varphi\widehat vJ_\br\,\dd\widehat\x,
\quad \widehat v\in H^1(\widehat\D),
\end{equation*}
whose existence, uniqueness, and bound \eqref{eq:bound_dual_J_r_new}
follow from Lemma~\ref{lem:bilinear_ext_D}.
\end{samepage}
Let $C_{\mathrm{is}}>0$ be the inf-sup constant from
Lemma~\ref{lem:bilinear_ext_D}, item~\normalfont{(iii)},
applied on $\mathcal O$, and $C_{\mathrm{p}}>0$ the constant from
\eqref{eq:bound_dual_J_r_new}.
For each $\br\in\mathfrak T_{\delta,\mathcal O}$, set
$\widehat u_\br=\widehat{\mathsf N}^{(\kappa)}_{\br,\mathcal O}\widehat\varphi$,
where the operator is the holomorphic extension supplied by
Theorem~\ref{thm:shape_hol_Newton_L2} applied to $\mathcal O$.
We have that
$\|\widehat u_\br\|_{H^1(\widehat\D)}
\leq\|\widehat u_\br\|_{H^1(\mathcal O)}$, and the inf-sup estimate yields
\begin{equation*}
\|\widehat u_\br\|_{H^1(\mathcal O)}
\leq\frac{1}{C_{\mathrm{is}}}
\sup_{\widehat w\in\widetilde H^{-1}(\mathcal O)\setminus\{0\}}
\frac{\left|\dual{\widehat w}{\widehat u_\br}_
{\widetilde H^{-1}(\mathcal O)\times H^1(\mathcal O),\br}\right|}
{\|\widehat w\|_{\widetilde H^{-1}(\mathcal O)}}.
\end{equation*}
Next, it follows from the definition of the bilinear form in
\eqref{eq:bilnfrmH1_2_O} that
\begin{equation*}
\begin{aligned}
\dual{\widehat w}{\widehat u_\br}_
{\widetilde H^{-1}(\mathcal O)\times H^1(\mathcal O),\br}
&=\int_{\mathcal O}
(D\br)^{-\top}\widehat\nabla(\mathcal R_{\mathcal O}^{-1}\widehat w)
\cdot(D\br)^{-\top}\widehat\nabla\widehat u_\br
J_\br\,\dd\widehat\x\\
&\quad+\int_{\mathcal O}
(\mathcal R_{\mathcal O}^{-1}\widehat w)\widehat u_\br
J_\br\,\dd\widehat\x.
\end{aligned}
\end{equation*}
Recalling that $\mathcal R_{\mathcal O}:H^1(\mathcal O)
\rightarrow\widetilde H^{-1}(\mathcal O)$ is an isometric isomorphism,
and that \eqref{eq:pulledpackpartialintegration2} holds for each
$\br\in\mathfrak T_{\delta,\mathcal O}$ by \textbf{Step~1}, we obtain
\begingroup
\small
\begin{equation*}
\begin{aligned}
\norm{\widehat u_\br}_{H^1(\mathcal O)}&\\
&\hspace{-1.3cm}\leq\frac{1}{C_{\mathrm{is}}}
\sup_{\widehat v\in H^1(\mathcal O)\setminus\{0\}}
\frac{
\left|\displaystyle\int_{\widehat\D}
\widehat\varphi\widehat vJ_\br\,\dd\widehat\x\right|
+\left|\displaystyle\int_{\widehat\D}\widehat\varphi
\mathsf B^{(\kappa)}_\br(\gamma_{\partial\mathcal O}\widehat v)
J_\br\,\dd\widehat\y\right|
+(1+\kappa^2)\left|\displaystyle\int_{\mathcal O}
\widehat u_\br\widehat vJ_\br\,\dd\widehat\x\right|
}{\norm{\widehat v}_{H^1(\mathcal O)}}.
\end{aligned}
\end{equation*}
\endgroup
Furthermore, recalling the explicit representation of the Newton potential from
\eqref{eq:Nhatexplicit}, and using symmetry of the fundamental solution and the
Fubini--Tonelli theorem, we obtain, for each $\widehat v\in H^1(\mathcal O)$,
\begin{equation*}
\begin{aligned}
\int_{\mathcal O}\widehat u_\br\widehat vJ_\br\,\dd\widehat\x
&=\int_{\widehat\D}\widehat\varphi(\widehat\y)
\left(\int_{\mathcal O}
\mathsf g^{(\kappa)}_{\br,\mathbb C}(\widehat\y,\widehat\x)
\widehat v(\widehat\x)\,\dd\widehat\x\right)
J_\br(\widehat\y)\,\dd\widehat\y\\
&=\int_{\widehat\D}\widehat\varphi
\left(\widehat{\mathsf N}^{(\kappa)}_{\br,\mathcal O}\widehat v\right)
J_\br\,\dd\widehat\y.
\end{aligned}
\end{equation*}
It follows from \eqref{eq:bound_dual_J_r_new} and
Theorem~\ref{thm:shape_hol_Newton_L2},
applied to $\mathcal O$, that
\begin{equation*}
\begin{aligned}
    \left|\int_{\widehat\D}\widehat\varphi
    \left(\widehat{\mathsf N}^{(\kappa)}_{\br,\mathcal O}\widehat v\right)
    J_\br\,\dd\widehat\y\right|
    &\leq C_{\mathrm{p}}\norm{\widehat{\ell}_{\widehat\varphi,\br}}_{\widetilde H^{-1}(\widehat\D)}
    \norm{\widehat{\mathsf N}^{(\kappa)}_{\br,\mathcal O}\widehat v}_{H^1(\widehat\D)}\\
    &\leq C_{\mathrm{p}}\norm{\widehat{\ell}_{\widehat\varphi,\br}}_{\widetilde H^{-1}(\widehat\D)}
    \norm{\widehat{\mathsf N}^{(\kappa)}_{\br,\mathcal O}}_{
    \mathscr L(L^2(\mathcal O),H^1(\mathcal O))}
    \norm{\widehat v}_{H^1(\mathcal O)},
\end{aligned}
\end{equation*}
where we have used
$\norm{\widehat v}_{L^2(\mathcal O)}\leq\norm{\widehat v}_{H^1(\mathcal O)}$.
Next, recalling \eqref{eq:bound_dual_J_r_new}, we obtain
\begin{equation*}
    \left|\int_{\widehat\D}\widehat\varphi\widehat vJ_\br\,\dd\widehat\x\right|
    \leq 
    C_{\mathrm{p}}\norm{\widehat{\ell}_{\widehat\varphi,\br}}_{\widetilde H^{-1}(\widehat\D)}
    \norm{\widehat v}_{H^1(\widehat\D)}
    \leq 
    C_{\mathrm{p}}\norm{\widehat{\ell}_{\widehat\varphi,\br}}_{\widetilde H^{-1}(\widehat\D)}
    \norm{\widehat v}_{H^1(\mathcal O)}.
\end{equation*}
Using similar arguments, continuity of $\mathsf B^{(\kappa)}_\br$,
continuity of the Dirichlet trace operator on $\partial\mathcal O$,
and \eqref{eq:bound_dual_J_r_new}, we obtain
\begin{equation*}
\begin{aligned}
    \left|\int_{\widehat\D}\widehat\varphi
    \left(\mathsf B^{(\kappa)}_\br\gamma_{\partial\mathcal O}\widehat v\right)
    J_\br\,\dd\widehat\y\right|
    \\
    &
    \hspace{-2.85cm}
    \leq C_{\mathrm{p}}\norm{\widehat{\ell}_{\widehat\varphi,\br}}_{\widetilde H^{-1}(\widehat\D)}
    \norm{\mathsf B^{(\kappa)}_\br\gamma_{\partial\mathcal O}\widehat v}_{H^1(\widehat\D)}\\
    &
    \hspace{-2.85cm}
    \leq C_{\mathrm{p}}\norm{\widehat{\ell}_{\widehat\varphi,\br}}_{\widetilde H^{-1}(\widehat\D)}
    \norm{\mathsf B^{(\kappa)}_\br}_{
    \mathscr L(H^{\frac{1}{2}}(\partial\mathcal O),H^1(\widehat\D))}
    \norm{\gamma_{\partial\mathcal O}}_{
    \mathscr L(H^1(\mathcal O),H^{\frac{1}{2}}(\partial\mathcal O))}
    \norm{\widehat v}_{H^1(\mathcal O)}.
\end{aligned}
\end{equation*}
Combining these estimates, we obtain
\begin{equation}\label{eq:NewtonBound}
\begin{aligned}
\norm{\widehat u_\br}_{H^1( \widehat \D)} & \leq\norm{\widehat u_\br}_{H^1(\mathcal O)}
\leq\frac{C_{\mathrm{p}}}{C_{\mathrm{is}}}\Bigg(
1+(1+\kappa^2)
\norm{\widehat{\mathsf N}^{(\kappa)}_{\br,\mathcal O}}_{
\mathscr L(L^2(\mathcal O),H^1(\mathcal O))}
\\
&\;+
\norm{\mathsf B^{(\kappa)}_\br}_{\mathscr L(H^{\frac{1}{2}}(\partial\mathcal O),H^1(\widehat\D))}
\norm{\gamma_{\partial\mathcal O}}_{\mathscr L(H^1(\mathcal O),H^{\frac{1}{2}}(\partial\mathcal O))}
\Bigg)\norm{\widehat{\ell}_{\widehat\varphi,\br}}_{\widetilde H^{-1}(\widehat\D)}\\
&\leq C(\mathfrak T_{\mathcal O},\delta)
e^{\kappa\delta C_2(\mathfrak T_{\mathcal O},\delta)}
\norm{\widehat{\ell}_{\widehat\varphi,\br}}_{\widetilde H^{-1}(\widehat\D)}
\begin{cases}
1+\kappa^{\frac{5}{2}}+\kappa^{-\frac{1}{2}},&d=2,\\
1+\kappa^3,&d=3
\end{cases}.
\end{aligned}
\end{equation}
Here $\kappa>0$, $\br\in\mathfrak T_{\delta,\mathcal O}$, and
$C(\mathfrak T_{\mathcal O},\delta)>0$ is independent of $\kappa$ and $\br$.
The last inequality follows from \eqref{eq:newton_reference_preliminary_bound} and \eqref{eq:BopNorm}.

For each fixed $\br\in\mathfrak T_{\delta,\mathcal O}$, the functionals
$\widehat\ell_{\widehat\varphi,\br}$, with
$\widehat\varphi\in\mathscr C^\infty_0(\widehat\D)$, are dense in
$\widetilde H^{-1}(\widehat\D)$.
Indeed, by Lemma~\ref{lem:bilinear_ext_D} and the Hahn--Banach theorem,
it suffices to show that a function $\widehat v\in H^1(\widehat\D)$ satisfying
\begin{equation*}
    \dual{\widehat\ell_{\widehat\varphi,\br}}{\widehat v}_
    {\widetilde H^{-1}(\widehat\D)\times H^1(\widehat\D),\br}
    =\int_{\widehat\D}\widehat\varphi\widehat vJ_\br\,\dd\widehat\x
    =0
\end{equation*}
for every $\widehat\varphi\in\mathscr C^\infty_0(\widehat\D)$ vanishes.
Indeed, the preceding identity gives $J_\br\widehat v=0$ almost everywhere.
By Assumption~\ref{assump:transformations_T} and \cite[Lemma~A.1]{HDPS26Laplace},
for sufficiently small $\delta>0$,
\begin{equation*}
\inf_{\br\in\mathfrak T_{\delta,\mathcal O}}
\essinf_{\widehat\x\in\widehat\D}
\operatorname{Re}J_\br(\widehat\x)
\geq C(\mathfrak T_{\mathcal O})>0.
\end{equation*}
Hence $\widehat v=0$.
Thus, for any $\widehat\ell\in\widetilde H^{-1}(\widehat\D)$, we can choose
$\widehat\varphi_n\in\mathscr C^\infty_0(\widehat\D)$ such that the associated
functionals $\widehat\ell_{n,\br}\coloneqq\widehat\ell_{\widehat\varphi_n,\br}$ converge to
$\widehat\ell$ in $\widetilde H^{-1}(\widehat\D)$.
It follows from \eqref{eq:NewtonBound} that $\widehat{\mathsf N}^{(\kappa)}_\br$
extends uniquely to a bounded operator from $\widetilde H^{-1}(\widehat\D)$
to $H^1(\widehat\D)$, with the same bound.
The extension argument in Step~4 of the proof of
\cite[Theorem~3.3]{HDPS26Laplace} yields the corresponding estimate for
$\br\in\mathfrak T_\delta$, after possibly decreasing $\delta_0$ and choosing
$0<\delta<\delta_0$.
Consequently, we obtain the following bound
\begin{equation}\label{eq:Newton_frequency_explicit}
\sup_{\br\in\mathfrak T_\delta}
\norm{\widehat{\mathsf N}^{(\kappa)}_\br}_{
\mathscr L(\widetilde H^{-1}(\widehat\D),H^1(\widehat\D))}
\leq C(\mathfrak T,\delta_0)e^{\kappa\delta A(\mathfrak T,\delta_0)}
\begin{cases}
1+\kappa^{\frac{5}{2}}+\kappa^{-\frac{1}{2}},&d=2,\\
1+\kappa^3,&d=3, 
\end{cases}
\end{equation}
for some constant $A=A(\mathfrak T,\delta_0)>0$.
This implies
\begin{equation}\label{eq:dual_newton_coarse}
\sup_{\br\in\mathfrak T_{\delta_0}}
\|\widehat{\mathsf N}_\br^{(\kappa)}\|_
{\mathscr L(\widetilde H^{-1}(\widehat\D),H^1(\widehat\D))}
\leq C(\mathfrak T,\delta_0,\kappa_0)
(1+\kappa)^{\frac{d+3}{2}}e^{\kappa\delta_0 A},
\quad\kappa\geq\kappa_0.
\end{equation}

\textbf{Step 3: Holomorphy of the shape-to-operator map.}
For the holomorphy argument, we relate the functionals
$\widehat\ell_{\widehat\varphi,\br}$ from \textbf{Step~2} to a dense family
of functionals independent of $\br$.
For $\br\in\mathfrak T_\delta$, define the weak Jacobian
$\mathsf J_\br:\widetilde H^{-1}(\widehat\D)\to
\widetilde H^{-1}(\widehat\D)$ by
\begin{equation*}
\dual{\mathsf J_\br\widehat\ell}{\widehat v}_
{\widetilde H^{-1}(\widehat\D)\times H^1(\widehat\D)}
=\dual{\widehat\ell}{\widehat v}_
{\widetilde H^{-1}(\widehat\D)\times H^1(\widehat\D),\br}.
\end{equation*}
It follows from Lemma~\ref{lem:bilinear_ext_D} and
\cite[Proposition~2.1, item~(ii)]{HDPS26Laplace} that $\mathsf J_\br$ is a bounded
isomorphism with uniformly bounded inverse and that there exists $\delta>0$
such that the maps
\begin{equation*}
\mathfrak T\rightarrow
\mathscr L\bigl(\widetilde H^{-1}(\widehat\D),\widetilde H^{-1}(\widehat\D)\bigr):
\quad\br\mapsto\mathsf J_\br,\quad\br\mapsto\mathsf J_\br^{-1}
\end{equation*}
admit bounded holomorphic extensions to $\mathfrak T_\delta$.

Set $\mathsf Q_\br\coloneqq\widehat{\mathsf N}^{(\kappa)}_\br\mathsf J_\br^{-1}$.
For fixed $f\in L^2(\widehat\D)$, let
\begin{equation*}
\dual{\widehat\ell}{\widehat v}_
{\widetilde H^{-1}(\widehat\D)\times H^1(\widehat\D)}
=\int_{\widehat\D}f\widehat v\,\dd\widehat\x,
\quad \widehat v\in H^1(\widehat\D).
\end{equation*}
The definition of the weak Jacobian gives
\begin{equation*}
\dual{\mathsf J_\br^{-1}\widehat\ell}{\widehat v}_
{\widetilde H^{-1}(\widehat\D)\times H^1(\widehat\D),\br}
=\int_{\widehat\D}\frac{f}{J_\br}\widehat vJ_\br\,\dd\widehat\x.
\end{equation*}
Since $f/J_\br\in L^2(\widehat\D)$, this identifies
$\mathsf J_\br^{-1}\widehat\ell=\widehat\ell_{f/J_\br,\br}$.
Thus, the explicit representation \eqref{eq:defExt}, as used in
\textbf{Step~2} yields
$\mathsf Q_\br\widehat\ell
=\widehat{\mathsf N}^{(\kappa)}_\br
\left(\frac{f}{J_\br}\right)$,
where the operator on the right is the $L^2(\widehat\D)$-data operator from
Theorem~\ref{thm:shape_hol_Newton_L2}. The map $\br\mapsto J_\br^{-1}$ is holomorphic in
$L^\infty(\widehat\D)$ by \cite[Lemma~A.1 and Proposition~2.1, item~(ii)]{HDPS26Laplace}.
Thus, Theorem~\ref{thm:shape_hol_Newton_L2} implies holomorphy of
$\br\mapsto\mathsf Q_\br\widehat\ell$ in $H^1(\widehat\D)$.

The functionals represented by $L^2(\widehat\D)$ densities in the standard
duality pairing are dense in $\widetilde H^{-1}(\widehat\D)$, and
\eqref{eq:Newton_frequency_explicit}, together with the uniform bound
for $\mathsf J_\br^{-1}$, gives a uniform bound for $\mathsf Q_\br$.
For any $\widehat\ell\in\widetilde H^{-1}(\widehat\D)$, choose
$f_n\in L^2(\widehat\D)$, independent of $\br$, such that the functionals defined by
\begin{equation*}
\dual{\widehat\ell_n}{\widehat v}_
{\widetilde H^{-1}(\widehat\D)\times H^1(\widehat\D)}
=\int_{\widehat\D}f_n\widehat v\,\dd\widehat\x,
\quad\widehat v\in H^1(\widehat\D),
\end{equation*}
converge to $\widehat\ell$ in $\widetilde H^{-1}(\widehat\D)$.
The uniform bound for $\mathsf Q_\br$ implies
\begin{equation*}
\sup_{\br\in\mathfrak T_\delta}
\|\mathsf Q_\br\widehat\ell_n-\mathsf Q_\br\widehat\ell\|_{H^1(\widehat\D)}
\longrightarrow0.
\end{equation*}
Since $\br\mapsto \mathsf Q_\br\widehat\ell_n$ is holomorphic for each $n$,
the limit $\br\mapsto \mathsf Q_\br\widehat\ell$ is holomorphic as well.
The argument in Step~5 of the proof of
\cite[Theorem~3.3]{HDPS26Laplace} therefore yields holomorphy of the map
$\mathfrak T_\delta\ni\br\mapsto\mathsf Q_\br\in
\mathscr L(\widetilde H^{-1}(\widehat\D),H^1(\widehat\D))$. Since
$\widehat{\mathsf N}^{(\kappa)}_\br=\mathsf Q_\br\mathsf J_\br$, the assertion follows.

\textbf{Step 4: Frequency-explicit operator norm bound.}
By Assumption~\ref{assump:transformations_T}, item~(i), we have
$\sup_{\widetilde\br\in\mathfrak T}
\|\widetilde\br\|_{L^\infty(\widehat\D)}<\infty$.
Thus, since $\D_\br$ is a bounded Lipschitz domain for all $\br\in\mathfrak T_{\delta_0}\cap
W^{1,\infty}(\widehat\D;\IR^d)$, cf.
\cite[Lemma~B.1, item~(i)]{HDPS26Laplace}, there is
$L>\delta_0+\sup_{\widetilde\br\in\mathfrak T}
\|\widetilde\br\|_{L^\infty(\widehat\D)}$ such that
\begin{equation}\label{eq:choice_B_L}
\bigcup_{\br\in\mathfrak T_{\delta_0}\cap
W^{1,\infty}(\widehat\D;\IR^d)}\overline{\D_{\br}}
\subset
B(\boldsymbol0,L)
\eqqcolon B_L.
\end{equation}
For $f\in L^2(B_L)$, define
\begin{equation*}
(\mathsf N^{(\kappa)}_{B_L} f)(\x)\coloneqq
\int_{B_L}G^{(\kappa)}(\x-\y)f(\y)\,\dd\y,
\quad\x\in B_L.
\end{equation*}
For the frequency-explicit estimate, we take our cue from
\cite[Eq.~(4.12)]{Spence2014Wavenumber}, which recalls the
free resolvent estimate of \cite[Theorems~3 and~4]{Vainberg1975ShortWave}.
For every fixed $\kappa_0>0$, there exists
$C_{\mathrm{res}}(L,\kappa_0)>0$, depending on $L$ and $\kappa_0$ only, such that
for all $f\in L^2(B_L)$ and any $\kappa\geq\kappa_0$
\begin{equation}\label{eq:dual_newton_free_H2}
\begin{aligned}
\kappa\|\mathsf N^{(\kappa)}_{B_L} f\|_{L^2(B_L)}
+\|\mathsf N^{(\kappa)}_{B_L} f\|_{H^1(B_L)}
+
&\kappa^{-1}\|\mathsf N^{(\kappa)}_{B_L} f\|_{H^2(B_L)}
\\
&\leq C_{\mathrm{res}}(L,\kappa_0)\|f\|_{L^2(B_L)}.
\end{aligned}
\end{equation}
The Fubini--Tonelli theorem and $G^{(\kappa)}(\x-\y)=G^{(\kappa)}(\y-\x)$ yield
\[
\int_{B_L}(\mathsf N^{(\kappa)}_{B_L} f)g \dd\x
=\int_{B_L}f(\mathsf N^{(\kappa)}_{B_L} g) \dd\x,
\]
for all $f,g\in L^2(B_L)$, implying that
$\mathsf N^{(\kappa)}_{B_L}:L^2(B_L)\rightarrow H^2(B_L)$ is symmetric when considered as an operator from and to $L^2(B_L)$.
Thus, using a density argument, it extends to an operator $\mathsf N^{(\kappa)}_{B_L}:\widetilde H^{-2}(B_L)\to L^2(B_L)$, and \eqref{eq:dual_newton_free_H2} yields for $\kappa\geq\kappa_0$
\[
\begin{aligned}
\|\mathsf N^{(\kappa)}_{B_L}\|_{\mathscr L(\widetilde H^{-2}(B_L),L^2(B_L))}
&=
\|\mathsf N^{(\kappa)}_{B_L}\|_{\mathscr L(L^2(B_L),H^2(B_L))}
\\
&\leq C_{\mathrm{res}}(L,\kappa_0)(1+\kappa).
\end{aligned}
\]
Interpolation yields for $\kappa\geq\kappa_0$
\begin{equation} \label{eq:freebound}
    \|\mathsf N^{(\kappa)}_{B_L}\|_{\mathscr L(\widetilde H^{-1}(B_L),H^1(B_L))}
\leq C_{\mathrm{res}}(L,\kappa_0)(1+\kappa).
\end{equation}
From the explicit representation of the Newton potential \eqref{eq:Nhatexplicit} and the identification in \textbf{Step~1}, we obtain that $\widehat{\mathsf N}^{(\kappa)}_\br \widehat\ell_{\widehat\varphi,\br}  = \tau_\br \mathsf N^{(\kappa)}_{B_L} \varphi_{\br}$ for all $\widehat\varphi\in\mathscr C_0^\infty(\widehat\D)$ and $\varphi_\br  = \widehat{\varphi} \circ \br^{-1}$.
Hence
\begin{align*}
    \norm{\widehat{\mathsf N}^{(\kappa)}_\br \widehat\ell_{\widehat\varphi,\br}}_{H^1(\widehat\D)} \leq C   \norm{ \mathsf N^{(\kappa)}_{B_L}  \varphi_\br}_{H^1(\D_\br)} \leq C   \norm{ \mathsf N^{(\kappa)}_{B_L} \varphi_\br}_{H^1(B_L)}.
\end{align*}
where $C=C(\mathfrak T,\delta_0)>0$ is independent of $\kappa$, $\br$, and $\widehat\varphi$.
Combining this with \eqref{eq:freebound} we arrive at
\begin{align}
\label{eq:freebound2}
    \norm{\widehat{\mathsf N}^{(\kappa)}_\br \widehat\ell_{\widehat\varphi,\br}}_{H^1(\widehat\D)} \leq C(L,\kappa_0) (1+\kappa) \|  \varphi_\br\|_{\widetilde{H}^{-1}(B_L)},\quad \kappa \geq \kappa_0,
\end{align}
where $\varphi_\br$ is extended by zero outside $\D_{\br}$. We can  then bound its norm as follows 
\begin{align*}
    \|\varphi_\br\|_{\widetilde{H}^{-1}(B_L)} &= \sup_{v \in H^1(B_L) \setminus \{0\}} \frac{\left\lvert\langle \varphi_\br , v \rangle_{\widetilde H^{-1}(B_L)\times H^1(B_L)}\right\rvert}{\|v\|_{H^1(B_L)}} \\& =  \sup_{v \in H^1(B_L)\setminus \{0\}} \frac{\left\lvert\langle \varphi_\br , v \rangle_{\widetilde H^{-1}(\D_{\br})\times H^1(\D_{\br})}\right\rvert}{\|v\|_{H^1(B_L)}}   \leq \|\varphi_\br\|_{\widetilde{H}^{-1}(\D_{\br})}.
\end{align*}
On the other hand, we have that 
\begin{align*}
    \|\varphi_\br\|_{\widetilde{H}^{-1}(\D_{\br})} &= \sup_{v \in H^1(\D_\br)\setminus \{0\}} \frac{\left\lvert\langle \varphi_\br , v \rangle_{\widetilde H^{-1}(\D_\br)\times H^1(\D_\br)}\right\rvert}{\|v\|_{H^1(\D_\br)}} \\
    &= \sup_{v \in H^1(\D_\br) \setminus \{0\}} \frac{\left\lvert\displaystyle\int_{\widehat\D} \widehat \varphi \tau_\br v J_\br  \ \d\widehat\x\right\rvert}{\|v\|_{H^1(\D_\br)}},
\end{align*}
Finally, using the pairing identity in \textbf{Step~1},
Lemma~\ref{lem:bilinear_ext_D}, item~\normalfont{(i)}, and the uniform bound of
$\tau_\br$ on the latter equation and combining with \eqref{eq:freebound2}, we obtain
\begin{equation*}
\norm{\widehat{\mathsf N}^{(\kappa)}_\br\widehat\ell_{\widehat\varphi,\br}}_{H^1(\widehat\D)}
\leq C_0(\mathfrak T,\delta_0,\kappa_0)(1+\kappa)
\norm{\widehat\ell_{\widehat\varphi,\br}}_{\widetilde H^{-1}(\widehat\D)}.
\end{equation*}
The density established in \textbf{Step~2} therefore gives the uniform bound
\begin{equation}\label{eq:dual_newton_real_bound}
\sup_{\br\in\mathfrak T_{\delta_0}\cap W^{1,\infty}(\widehat\D;\IR^d)}
\|\widehat{\mathsf N}_\br^{(\kappa)}\|_
{\mathscr L(\widetilde H^{-1}(\widehat\D),H^1(\widehat\D))}
\leq C_0(\mathfrak T,\delta_0,\kappa_0)(1+\kappa),
\quad\kappa\geq\kappa_0.
\end{equation}
Step~3 implies that the map
$\br\mapsto\widehat{\mathsf N}^{(\kappa)}_\br$ admits a bounded holomorphic extension to
$\mathfrak T_{\delta_0}$. Applying Corollary~\ref{cor:abstract_frequency_bound} to \eqref{eq:dual_newton_coarse} and
\eqref{eq:dual_newton_real_bound} implies \eqref{eq:Newton_explicit_refined_frequency}.
\end{proof}

\subsection{Shape Holomorphy of the Single and Double Layer Potentials}
\label{sec:layer_potentials_hol}
We prove shape holomorphy of the single and double layer potentials.
To this end, for each $\br \in \mathfrak{T}$, we set
\begin{equation*}\label{eq:single_layer_potential}
    \widehat{\mathsf{S}}^{(\kappa)}_{\br}
    =
    \tau_{\br} \,
    {\mathsf{S}}^{(\kappa)}_{\Gamma_\br}
    \tau^{-1}_{\br}
    \in
    \mathscr{L}
    \left(	
	{H}^{-\frac{1}{2}}(\normalfont\widehat{\Gamma})
	,
	H^{1}({\widehat{\text{D}}})
    \right),
\end{equation*}
where the single layer potential
${\mathsf{S}}^{(\kappa)}_{\Gamma_\br}: {H}^{-\frac{1}{2}}(\Gamma_\br) \rightarrow
H^{1}(\text{D}_\br)$ are introduced in \eqref{eq:single_layer}. 

\begin{theorem}\label{thm:holo_single_layer_pot}
Let $\mathfrak T$ be a set of domain transformations
satisfying Assumptions~\ref{assump:transformations_T}
and~\ref{assump:domain_T_2}.
Then there exists $\delta_0=\delta_0(\mathfrak T)>0$, such that,
for every $\kappa>0$ and $0<\delta<\delta_0$, the map
\begin{equation*}\label{eq:map_hol_single_potential}
	\mathcal{S}^{(\kappa)}:
	\mathfrak{T}
	\rightarrow
	\mathscr{L}
	\left(	
		{H}^{-\frac{1}{2}}(\normalfont\widehat{\Gamma})
		,
		H^{1}({\widehat{\text{D}}})
	\right):
	\br
	\mapsto
	\widehat{\mathsf{S}}^{(\kappa)}_{\br}
\end{equation*}
admits a bounded holomorphic extension to $\mathfrak T_\delta$.
Moreover, there exists $C_2=C_2(\mathfrak T,\delta_0)>0$ such that, for every
$\kappa_0>0$, there exists $C_1=C_1(\mathfrak T,\delta_0,\kappa_0)>0$
satisfying
\begin{equation}\label{eq:single_layer_refined_frequency}
\sup_{\br\in\mathfrak T_\delta}
\norm{\widehat{\mathsf S}^{(\kappa)}_\br}_{\mathscr L(H^{-\frac{1}{2}}(\widehat\Gamma),H^1(\widehat\D))}
\leq C_1(1+\kappa)
e^{\frac{4\delta}{\pi R}\left(\kappa\delta_0C_2+\frac{d+1}{2}\log(1+\kappa)\right)}
\end{equation}
for all $\kappa\geq\kappa_0$, $0<R<\delta_0$, and
$0<\delta<\min\{R,\delta_0-R\}$.
\end{theorem}

\begin{proof}
Recalling that ${\mathsf{S}}^{(\kappa)}_{\Gamma_\br} = {\mathsf{N}}^{(\kappa)}_{\br} \gamma^\prime_{\Gamma_\br}$,
for each $\br \in \mathfrak{T}$, cf.~\cite[Definition~3.1.5]{SS10}, we obtain
\begin{equation*}
	\widehat{\mathsf{S}}^{(\kappa)}_{\br}
    =
    \tau_{\br} 
    {\mathsf{N}}^{(\kappa)}_{\br} 
    \tau^{-1}_{\br}
    \tau_{\br}
    \gamma^\prime_{\Gamma_\br}
    \tau^{-1}_{\br} 
    = 
    \widehat{\mathsf{N}}^{(\kappa)}_\br
    \tau_{\br}
    \gamma^\prime_{\Gamma_\br}
    \tau^{-1}_{\br} .
\end{equation*}
Theorem~\ref{thm:shape_hol_Newton} implies that there exists $\delta_0>0$ such that the map
$\mathfrak T\ni\br\mapsto\widehat{\mathsf N}^{(\kappa)}_\br
\in\mathscr L(\widetilde H^{-1}(\widehat\D),H^1(\widehat\D))$
admits a bounded holomorphic extension to $\mathfrak T_{\delta_0}$.
Following the proof of \cite[Theorem~3.6]{HDPS26Laplace} and using
Lemmas~\ref{lem:bilinear_ext_D} and~\ref{lem:bilinear_Gamma}, we obtain the map
\begin{equation*}
\mathfrak T\rightarrow
\mathscr L(H^{-\frac{1}{2}}(\widehat\Gamma),\widetilde H^{-1}(\widehat\D)):
\br\mapsto\tau_\br\gamma'_{\Gamma_\br}\tau_\br^{-1}
\end{equation*}
admits a bounded holomorphic extension to $\mathfrak T_{\delta_0}$,
for a possibly smaller $\delta_0>0$.
Moreover, item~\normalfont{(i)} of
Lemma~\ref{lem:bilinear_Gamma} and item~\normalfont{(iii)} of
Lemma~\ref{lem:bilinear_ext_D} imply that
\begin{equation*}
\sup_{\br\in\mathfrak T_{\delta_0}}
\norm{\tau_\br\gamma'_{\Gamma_\br}\tau_\br^{-1}}_
{\mathscr L(H^{-\frac{1}{2}}(\widehat\Gamma),\widetilde H^{-1}(\widehat\D))}
\leq
\frac{C_1(\mathfrak T,\delta_0)}{C_2(\mathfrak T,\delta_0)}
\norm{\gamma_{\widehat\Gamma}}_
{\mathscr L(H^1(\widehat\D),H^{\frac{1}{2}}(\widehat\Gamma))}.
\end{equation*}
The constants $C_1$ and $C_2$ are given in
\eqref{eq:C_1_boundary} and~\eqref{eq:inf_sup_sesq_r_D}, respectively.

Finally, combining the Newton potential bound
\eqref{eq:Newton_explicit_refined_frequency} with the bound for
$\tau_\br\gamma'_{\Gamma_\br}\tau_\br^{-1}$ we obtain
\eqref{eq:single_layer_refined_frequency}.%
\end{proof}

 For each
$\br \in \mathfrak{T}$, we set
\begin{equation*}
    \widehat{\mathsf{D}}^{(\kappa)}_{\br}
    \coloneqq
    \tau_{\br} \,
    {\mathsf{D}}^{(\kappa)}_{\Gamma_\br}
    \tau^{-1}_{\br}
    \in
    \mathscr{L}
    \left(	
	{H}^{\frac{1}{2}}(\normalfont\widehat{\Gamma})
	,
	H^{1}({\widehat{\text{D}}})
    \right),
\end{equation*}
where the double layer potential $ {\mathsf{D}}^{(\kappa)}_{\Gamma_\br}$
has been introduced in \eqref{eq:double_layer}.

\begin{theorem}\label{thm:DLhlmrp}
Let $\mathfrak T$ be a set of domain transformations
satisfying Assumptions~\ref{assump:transformations_T}
and~\ref{assump:domain_T_2}.
Then there exists $\delta_0=\delta_0(\mathfrak T)>0$ such that,
for every $\kappa>0$ and $0<\delta<\delta_0$, the map
\begin{equation*}\label{eq:map_hol_D_H_1/2_H1}
	\mathcal{D}^{(\kappa)}:
	\mathfrak{T}
	\rightarrow
	\mathscr{L}
	\left(	
		{H}^{\frac{1}{2}}(\normalfont\widehat{\Gamma})
		,
		H^{1}(\widehat{\text{D}})
	\right):
	\br
	\mapsto
	\widehat{\mathsf{D}}^{(\kappa)}_{\br}
\end{equation*}
admits a bounded holomorphic extension to $\mathfrak T_\delta$.
Moreover, there exists $C_2=C_2(\mathfrak T,\delta_0)>0$ such that, for every
$\kappa_0>0$, there exists $C_1=C_1(\mathfrak T,\delta_0,\kappa_0)>0$
satisfying
\begin{equation}\label{eq:double_layer_main_frequency}
\sup_{\br\in\mathfrak T_\delta}
\norm{\widehat{\mathsf D}^{(\kappa)}_\br}_{\mathscr L(H^{\frac{1}{2}}(\widehat\Gamma),H^1(\widehat\D))}
\leq C_1(1+\kappa)^{\frac{3}{2}}
e^{\frac{4\delta}{\pi R}\left(\kappa\delta_0C_2+\frac{1}{2}\log(1+\kappa)\right)}
\end{equation}
for all $\kappa\geq\kappa_0$, $0<R<\delta_0$, and
$0<\delta<\min\{R,\delta_0-R\}$.
\end{theorem}

\begin{proof}
We proceed in several steps as follows.

\textbf{Step 1: Shape holomorphy.}
For each $\br\in\mathfrak T$, let
$\widehat\phi\in H^{\frac{1}{2}}(\widehat\Gamma)$, $\phi=\tau_\br^{-1}\widehat\phi$, and
consider the solution operator
\begin{equation*}
\widetilde{\mathsf T}_\br:
H^{\frac{1}{2}}(\Gamma_\br)\rightarrow H^1(\D_\br):
\phi\mapsto\widetilde u_{\br},
\end{equation*}
of the Laplacian boundary value problem
\begin{equation*}
-\Delta u_\br
=
0
\quad
\text{in}~\D_\br
\qquad
\text{and}
\qquad
\gamma_{\Gamma_\br}u_\br=\phi.
\end{equation*}
For each $\br \in \mathfrak{T}$, we define
\begin{align*}
\widehat{\mathsf T}_\br
&\coloneqq\tau_\br\widetilde{\mathsf T}_\br\tau_\br^{-1}\in \mathscr L\left(H^{\frac{1}{2}}(\widehat\Gamma),H^1(\widehat\D)\right),\\
\widehat{\gamma_{1,\Gamma_\br}\widetilde{\mathsf T}_\br}
&\coloneqq
\tau_\br\gamma_{1,\Gamma_\br}\widetilde{\mathsf T}_\br\tau_\br^{-1}
\in\mathscr{L}\left(H^{\frac{1}{2}}(\widehat{\Gamma}),
H^{-\frac{1}{2}}(\widehat{\Gamma})\right).
\end{align*}
Observing that $\mathcal L^{(\kappa)}\widetilde{\mathsf T}_\br\varphi
=-\kappa^2\widetilde{\mathsf T}_\br\varphi$
and inserting it into \eqref{eq:repformula}, we obtain
\begin{equation*}
    {\mathsf{D}}^{(\kappa)}_{\Gamma_\br} \phi
    =
    -\widetilde{\mathsf T}_\br \phi
    -\kappa^2\mathsf{N}^{(\kappa)}_\br
    \widetilde{\mathsf T}_\br\phi
    +
    {\mathsf{S}}^{(\kappa)}_{\Gamma_\br} \gamma_{1,\Gamma_\br} \widetilde{\mathsf T}_\br \phi,
\end{equation*}
which pulled back to the reference domain reads
\begin{equation}\label{eq:decomp_double_layer}
    \widehat{\mathsf{D}}^{(\kappa)}_\br \widehat\phi
    =
    -\widehat{\mathsf T}_\br \widehat{\phi}
    -\kappa^2\widehat{\mathsf{N}}^{(\kappa)}_{\br}
    \widehat{\mathsf T}_\br\widehat{\phi}
    +
    \widehat{\mathsf{S}}^{(\kappa)}_\br
    \widehat{\gamma_{1,\Gamma_\br}\widetilde{\mathsf T}_\br}
    \widehat{\phi}.
\end{equation}
With the shape holomorphy of $\br\mapsto\widehat{\mathsf{T}}_\br$ being shown in \cite[Theorem~3.4]{HDPS26Laplace} and the one of $\br\mapsto\widehat{\gamma_{1,\Gamma_\br}\widetilde{\mathsf T}_\br}$ being shown in the proof of \cite[Theorem~3.7]{HDPS26Laplace}, 
Theorem~\ref{thm:shape_hol_Newton} and Theorem~\ref{thm:holo_single_layer_pot} imply holomorphy of $\mathcal{D}^{(\kappa)}$. It thus remains to show the frequency-explicit operator norm bound.

\textbf{Step 2: Crude frequency-explicit operator norm bound.}
Theorem~\ref{thm:shape_hol_Newton_L2} and Theorem~\ref{thm:holo_single_layer_pot}
imply that there exists
$\delta_0, C_2>0$ such that, for every $\kappa_0>0$, there exists $C_1=C_1(\kappa_0)>0$
such that
\begin{equation}\label{eq:kNewton}
\begin{aligned}
\sup_{\br\in\mathfrak T_\delta}
\norm{\widehat{\mathsf N}^{(\kappa)}_\br}_
{\mathscr L(L^2(\widehat\D),H^1(\widehat\D))}
&\leq C_1e^{\kappa\delta C_2},\\
\sup_{\br\in\mathfrak T_\delta}
\norm{\widehat{\mathsf S}^{(\kappa)}_\br}_
{\mathscr L(H^{-\frac{1}{2}}(\widehat\Gamma),H^1(\widehat\D))}
&\leq C_1(1+\kappa)e^{\kappa\delta C_2}.
\end{aligned}
\end{equation}
for all $\kappa\geq\kappa_0$ and $0<\delta\leq\delta_0$.
With $\widehat{\mathsf T}_\br$ independent of $\kappa$,  \eqref{eq:decomp_double_layer} and
\eqref{eq:kNewton} imply that there exists $C_1=C_1(\kappa_0)>0$ such that
\begin{equation}\label{eq:double_layer_preliminary_bound}
\sup_{\br\in\mathfrak T_\delta}
\norm{\widehat{\mathsf D}^{(\kappa)}_\br}_
{\mathscr L(H^{\frac{1}{2}}(\widehat\Gamma),H^1(\widehat\D))}
\leq C_1(1+\kappa)^2e^{\kappa\delta C_2}.
\end{equation}

\textbf{Step 3: Improved frequency-explicit operator norm bound.}
Our goal is to apply Corollary~\ref{cor:abstract_frequency_bound}, for which we show that for real domain transformations in $\mathfrak T_{\delta_0}\cap W^{1,\infty}(\widehat\D;\IR^d)$ there exists 
$C_0=C_0(\mathfrak T,\delta_0,\kappa_0)>0$ such that
\begin{equation}\label{eq:double_layer_real_bound}
\sup_{\br\in\mathfrak T_{\delta_0}\cap W^{1,\infty}(\widehat\D;\IR^d)}
\norm{\widehat{\mathsf D}^{(\kappa)}_\br}_
{\mathscr L(H^{\frac{1}{2}}(\widehat\Gamma),H^1(\widehat\D))}
\leq C_0(1+\kappa)^{\frac{3}{2}},\quad\kappa\geq\kappa_0.
\end{equation}

To show \eqref{eq:double_layer_real_bound}, we exploit that real domain transformations $\br\in\mathfrak T_{\delta_0}\cap W^{1,\infty}(\widehat\D;\IR^d)$ yield bounded Lipschitz domains $\D_\br$, cf.\ \cite[Lemma~B.1, item~(i)]{HDPS26Laplace}, for which
\cite[Lemma~4.3]{Spence2014Wavenumber},
implies that there is $C=C(\mathfrak T,\delta_0,\kappa_0)>0$
such that for all $\phi\in H^{\frac{1}{2}}(\Gamma_\br)$,
\begin{equation}\label{eq:double_layer_real_L2}
\norm{\mathsf D^{(\kappa)}_{\Gamma_\br}\phi}_{L^2(\mathcal O)}
\leq C(1+\kappa)^{\frac{1}{2}}\norm{\phi}_{L^2(\Gamma_\br)},
\quad\kappa\geq\kappa_0.
\end{equation}

In the following we extend this estimate from the $L^2(\mathcal{O})$- to the $H^1(\widehat{\D})$-norm. To this end, define
\[
	\mathsf{E}_\br\phi
	=
	\begin{cases}
		-\widetilde{\mathsf T}_\br\phi,&\text{in }\D_\br,\\
		0,&\text{in }\mathcal O\setminus\overline{\D_\br}.
	\end{cases}
\]
The mapping properties of $\mathsf{D}_{\Gamma_\br}^{(\kappa)}$ and $\mathsf{E}_\br$ imply
$\mathsf D^{(\kappa)}_{\Gamma_\br}\phi,\mathsf{E}_\br\phi\in H^1(\D_\br)$ and
$\mathsf D^{(\kappa)}_{\Gamma_\br}\phi,\mathsf{E}_\br\phi\in H^1(\mathcal O\setminus\overline{\D_\br})$.
Recalling that $\gamma_{\Gamma_\br}\widetilde{\mathsf T}_\br\phi=\phi$,
it follows from the jump relations $[\mathsf D^{(\kappa)}_{\Gamma_\br}\phi-\mathsf{E}_\br\phi]=0$ on
$\Gamma_\br$, with $[\cdot]$ referring to the jump from \eqref{eq:jumps}. This implies $\mathsf D^{(\kappa)}_{\Gamma_\br}\phi-\mathsf{E}_\br\phi\in H^1(\mathcal O)$.

Now, let $\chi\in\mathscr{C}_0^\infty(\mathcal{O})$, with $0\leq\chi\leq 1$ and $\chi=1$ in a sufficiently small neighbourhood of $\cup_{\br \in \mathfrak T_{\delta_0}\cap W^{1,\infty}(\widehat\D;\IR^d)}\D_\br$ and consider
\[
z_\br=\chi^2\mathsf D^{(\kappa)}_{\Gamma_\br}\phi-\mathsf{E}_\br\phi=\chi^2(\mathsf D^{(\kappa)}_{\Gamma_\br}\phi-\mathsf{E}_\br\phi)\in H^1_0(\mathcal O).
\]
Since $\mathsf D^{(\kappa)}_{\Gamma_\br}\phi$ satisfies the homogeneous Helmholtz equation in both
$\D_\br$ and $\mathcal O\setminus\overline{\D_\br}$, the trace of
$z$ on $\partial\mathcal O$ vanishes, and the jumps of the Neumann trace of $\mathsf D^{(\kappa)}_{\Gamma_\br}$ and Dirichlet trace of $z_\br$ vanish on $\Gamma_\br$, the distributional Green's formula gives
\begin{equation*}
\int_{\mathcal O}
\bigl(\nabla \mathsf D^{(\kappa)}_{\Gamma_\br}\phi\cdot\nabla\overline z_\br
-\kappa^2\mathsf D^{(\kappa)}_{\Gamma_\br}\phi\overline z_\br\bigr)\,\dd\x
=0,
\end{equation*}
where $\nabla \mathsf D^{(\kappa)}_{\Gamma_\br}\phi$ is understood as the $L^2(\mathcal O)$ function defined in $\mathcal O\setminus \D_\br$ and $\D_\br$.
Inserting $z_\br=\chi^2\mathsf D^{(\kappa)}_{\Gamma_\br}\phi-\mathsf{E}_\br\phi$ we obtain
\begin{equation*}
\begin{aligned}
\norm{\chi\nabla \mathsf D^{(\kappa)}_{\Gamma_\br}\phi}_{L^2(\mathcal O)}^2
={}&\kappa^2\norm{\chi \mathsf D^{(\kappa)}_{\Gamma_\br}\phi}_{L^2(\mathcal O)}^2
-2\operatorname{Re}\left\{
\int_{\mathcal O}
\chi\overline{\mathsf D^{(\kappa)}_{\Gamma_\br}\phi}\,\nabla \mathsf D^{(\kappa)}_{\Gamma_\br}\phi\cdot\nabla\chi\,\dd\x\right\}\\
&+\operatorname{Re}\left\{
\int_{\D_\br}\nabla \mathsf D^{(\kappa)}_{\Gamma_\br}\phi\cdot\nabla\overline{\mathsf{E}_\br\phi}\,\dd\x\right\}
-\kappa^2\operatorname{Re}\left\{
\int_{\D_\br}\mathsf D^{(\kappa)}_{\Gamma_\br}\phi\overline{\mathsf{E}_\br\phi}\,\dd\x\right\}.
\end{aligned}
\end{equation*}
Cauchy--Schwarz and Young's inequality yield
\begin{equation*}
\begin{aligned}
&\frac{1}{2}\norm{\chi\nabla \mathsf D^{(\kappa)}_{\Gamma_\br}\phi}_{L^2(\mathcal O)}^2\\
&\qquad\leq\left(\frac{3}{2}\kappa^2
+4\norm{\nabla\chi}_{L^\infty(\mathcal O)}^2\right)
\norm{\mathsf D^{(\kappa)}_{\Gamma_\br}\phi}_{L^2(\mathcal O)}^2
+\norm{\nabla\mathsf{E}_\br\phi}_{L^2(\D_\br)}^2
+\frac{\kappa^2}{2}\norm{\mathsf{E}_\br\phi}_{L^2(\D_\br)}^2,
\end{aligned}
\end{equation*}
which implies
\begin{equation}\label{eq:double_layer_energy_bound}
\norm{\mathsf D^{(\kappa)}_{\Gamma_\br}\phi}_{H^1(\D_\br)}
\leq C_\chi(1+\kappa)
\left(\norm{\mathsf D^{(\kappa)}_{\Gamma_\br}\phi}_{L^2(\mathcal O)}
+\norm{\mathsf{E}_\br\phi}_{H^1(\D_\br)}\right).
\end{equation}
Combining \eqref{eq:double_layer_real_L2} and
\eqref{eq:double_layer_energy_bound} and observing that
\[
\norm{\mathsf{E}_\br\phi}_{H^1(\D_\br)}
\leq
\norm{\widetilde{\mathsf T}_\br}_
{\mathscr L(H^{\frac{1}{2}}(\Gamma_\br),H^1(\D_\br))}
\norm{\phi}_{H^{\frac{1}{2}}(\Gamma_\br)},
\]
yields
\begin{equation*}
\norm{\mathsf D^{(\kappa)}_{\Gamma_\br}\phi}_{H^1(\D_\br)}
\leq C(1+\kappa)^{\frac{3}{2}}
\norm{\phi}_{H^{\frac{1}{2}}(\Gamma_\br)}.
\end{equation*}
Finally, for $\widehat\phi\in H^{\frac{1}{2}}(\widehat\Gamma)$, set
$\phi=\tau_\br^{-1}\widehat\phi$.
Pullback to the reference domain yields \eqref{eq:double_layer_real_bound}.
Applying Corollary~\ref{cor:abstract_frequency_bound} to \eqref{eq:double_layer_real_bound} and
\eqref{eq:double_layer_preliminary_bound} implies
\eqref{eq:double_layer_main_frequency}, and thus the assertion.
\end{proof}

\subsection{Shape Holomorphy of the Single Layer BIO}
\label{sec:single_layer_hol}
We recall from \eqref{eq:transported_calderon_bios} the definition of
$\widehat{\mathsf{V}}^{(\kappa)}_\br$, with $\mathsf{V}^{(\kappa)}_{\Gamma_\br}$ defined in
\eqref{eq:def_single_layer}. In the following, we establish shape holomorphy of the
map
$\mathfrak{T}\ni\br\mapsto\widehat{\mathsf{V}}^{(\kappa)}_\br\in
\mathscr{L}\bigl(H^{-\frac{1}{2}}(\widehat{\Gamma}),
H^{\frac{1}{2}}(\widehat{\Gamma})\bigr)$.
\begin{theorem}\label{thm:Vhlmrp}
Let $\mathfrak T$ be a set of domain transformations
satisfying Assumptions~\ref{assump:transformations_T}
and~\ref{assump:domain_T_2}.
Then there exists $\delta_0>0$ such that,
for every $\kappa>0$ and $0<\delta<\delta_0$, the map
\begin{equation*}\label{eq:map_hol_V_H_-1/2_H_1/2}
	\mathcal{V}^{(\kappa)}:
	\mathfrak{T}
	\rightarrow
	\mathscr{L}
	\left(
		H^{-\frac{1}{2}}(\widehat\Gamma),
		H^{\frac{1}{2}}(\widehat\Gamma)
	\right):
	\br
	\mapsto
	\widehat{\mathsf V}^{(\kappa)}_\br
\end{equation*}
admits a bounded holomorphic extension to $\mathfrak T_\delta$.
Moreover, there exists $C_2=C_2(\mathfrak T,\delta_0)>0$ such that, for every
$\kappa_0>0$, there exists $C_1=C_1(\mathfrak T,\delta_0,\kappa_0)>0$
satisfying
\begin{equation}\label{eq:V_refined_frequency}
\sup_{\br\in\mathfrak T_\delta}
\norm{\widehat{\mathsf V}^{(\kappa)}_\br}_
{\mathscr L(H^{-\frac{1}{2}}(\widehat\Gamma),H^{\frac{1}{2}}(\widehat\Gamma))}
\leq C_1(1+\kappa)^{\frac{d-1}{2}}
e^{\frac{4\delta}{\pi R}\left(\kappa\delta_0C_2+2\log(1+\kappa)\right)}
\end{equation}
for all $\kappa\geq\kappa_0$, $0<R<\delta_0$, and
$0<\delta<\min\{R,\delta_0-R\}$.
\end{theorem}

\begin{proof}
By definition, we have that, for each $\br\in\mathfrak T$,
\begin{equation*}
\widehat{\mathsf V}^{(\kappa)}_\br
=\tau_\br\gamma_{\Gamma_\br}\tau_\br^{-1}
\widehat{\mathsf S}^{(\kappa)}_\br
=\gamma_{\widehat\Gamma}\widehat{\mathsf S}^{(\kappa)}_\br.
\end{equation*}
The holomorphy assertion then follows from
Theorem~\ref{thm:holo_single_layer_pot} and the continuity of the
Dirichlet trace operator
$\gamma_{\widehat\Gamma}:H^1(\widehat\D)
\rightarrow H^{\frac{1}{2}}(\widehat\Gamma)$.

To show the frequency-explicit operator norm bound \eqref{eq:V_refined_frequency}, choose $\delta_0>0$ such that Theorems~\ref{thm:shape_hol_Newton}
and~\ref{thm:holo_single_layer_pot} and
\cite[Lemma~B.1, item~(i)]{HDPS26Laplace} hold on
$\mathfrak T_{\delta_0}$. Observing that
\[
\widehat{\mathsf V}^{(\kappa)}_\br
=
\tau_\br\gamma_{\Gamma_\br}\tau_\br^{-1}
\widehat{\mathsf S}^{(\kappa)}_\br
=
\tau_\br\gamma_{\Gamma_\br}\tau_\br^{-1}
\widehat{\mathsf{N}}^{(\kappa)}_\br
\tau_{\br}
\gamma^\prime_{\Gamma_\br}
\tau^{-1}_{\br} 
\]
the frequency explicit operator norm bound
\eqref{eq:dual_newton_coarse} from
Theorem~\ref{thm:shape_hol_Newton} implies that there exist
$C_1=C_1(\mathfrak T,\delta_0,\kappa_0)>0$ and
$C_2=C_2(\mathfrak T,\delta_0)>0$ such that
\begin{equation} \label{eq:VprevboundComplex}
\sup_{\br\in\mathfrak T_{\delta_0}}
\norm{\widehat{\mathsf V}^{(\kappa)}_\br}_
{\mathscr L(H^{-\frac{1}{2}}(\widehat\Gamma),H^{\frac{1}{2}}(\widehat\Gamma))}
\leq C_1(1+\kappa)^{\frac{d+3}{2}}e^{\kappa\delta_0C_2},
\quad\kappa\geq\kappa_0.
\end{equation}
Aiming to use Corollary~\ref{cor:abstract_frequency_bound}, we further observe that for each $\br\in\mathfrak T_{\delta_0}\cap
W^{1,\infty}(\widehat\D;\IR^d)$ 
\cite[Corollary~1.2, item~(i)]{GrahamLohndorfMelenkSpence2015}
and Proposition~\ref{prop:iso_pullback} imply that there is $C_0=C_0(\mathfrak T,\delta_0,\kappa_0)>0$ such that
\begin{equation}\label{eq:single_layer_real_bound}
\sup_{\br\in\mathfrak T_{\delta_0}\cap W^{1,\infty}(\widehat\D;\IR^d)}
\norm{\widehat{\mathsf V}^{(\kappa)}_\br}_
{\mathscr L(H^{-\frac{1}{2}}(\widehat\Gamma),H^{\frac{1}{2}}(\widehat\Gamma))}
\leq C_0(1+\kappa)^{\frac{d-1}{2}},\quad\kappa\geq\kappa_0.
\end{equation}
Finally, the bound \eqref{eq:V_refined_frequency} follows from Corollary~\ref{cor:abstract_frequency_bound} using \eqref{eq:VprevboundComplex} and \eqref{eq:single_layer_real_bound}.
\end{proof}

\subsection{\texorpdfstring{Shape Holomorphy of the Double Layer BIO and Its Adjoint}{Shape Holomorphy of the Double Layer BIO and Its Adjoint}}
\label{sec:ShHolDLOp}

We recall from \eqref{eq:transported_calderon_bios} the definitions of
$\widehat{\mathsf{K}}^{(\kappa)}_\br$ and $\widehat{\mathsf{K}}^{(\kappa)\prime}_\br$, with the
corresponding BIOs defined in \eqref{eq:def_double_layer} and
\eqref{eq:def_adjoint_double_layer}. We establish shape holomorphy of the
maps $\mathfrak{T}\ni\br\mapsto\widehat{\mathsf{K}}^{(\kappa)}_\br\in
\mathscr{L}\bigl(H^{\frac{1}{2}}(\widehat{\Gamma}),
H^{\frac{1}{2}}(\widehat{\Gamma})\bigr)$ and
$\mathfrak{T}\ni\br\mapsto\widehat{\mathsf{K}}^{(\kappa)\prime}_\br\in
\mathscr{L}\bigl(H^{-\frac{1}{2}}(\widehat{\Gamma}),
H^{-\frac{1}{2}}(\widehat{\Gamma})\bigr)$.

\begin{theorem}\label{thm:KHlm}
Let $\mathfrak T$ be a set of domain transformations
satisfying Assumptions~\ref{assump:transformations_T}
and~\ref{assump:domain_T_2}.
Then there exists $\delta_0=\delta_0(\mathfrak T)>0$ such that,
for every $\kappa>0$ and $0<\delta<\delta_0$, the map
\begin{equation*}\label{eq:map_hol_K_H_+-1/2_H_+-1/2}
	\mathcal{K}^{(\kappa)}:
	\mathfrak{T}
	\rightarrow
	\mathscr{L}
	\left(
		H^{\frac{1}{2}}(\widehat\Gamma),
		H^{\frac{1}{2}}(\widehat\Gamma)
	\right):
	\br
	\mapsto
	\widehat{\mathsf K}^{(\kappa)}_\br
\end{equation*}
admits a bounded holomorphic extension to $\mathfrak T_\delta$.
Moreover, there exists $C_2=C_2(\mathfrak T,\delta_0)>0$ such that, for every
$\kappa_0>0$, there exists $C_1=C_1(\mathfrak T,\delta_0,\kappa_0)>0$
satisfying
\begin{equation}\label{eq:K_explicit_frequency}
\sup_{\br\in\mathfrak T_\delta}
\norm{\widehat{\mathsf K}^{(\kappa)}_\br}_
{\mathscr L(H^{\frac{1}{2}}(\widehat\Gamma),H^{\frac{1}{2}}(\widehat\Gamma))}
\leq C_1(1+\kappa)^{\frac{3}{2}}
e^{\frac{4\delta}{\pi R}\left(\kappa\delta_0C_2+\frac{1}{2}\log(1+\kappa)\right)}
\end{equation}
for all $\kappa\geq\kappa_0$, $0<R<\delta_0$, and
$0<\delta<\min\{R,\delta_0-R\}$.
\end{theorem}

\begin{proof}
Recalling the definition of the double layer operator,
we obtain, for each $\br\in\mathfrak T$,
\begin{equation*}
\widehat{\mathsf K}^{(\kappa)}_\br
=\tau_\br\left(\gamma_{\Gamma_\br}
\mathsf D^{(\kappa)}_{\Gamma_\br}
+\frac{1}{2}\mathsf I_{\Gamma_\br}\right)\tau_\br^{-1}
=\gamma_{\widehat\Gamma}\widehat{\mathsf D}^{(\kappa)}_\br
+\frac{1}{2}\mathsf I_{\widehat\Gamma}.
\end{equation*}
The holomorphy assertion follows from Theorem~\ref{thm:DLhlmrp}
and the continuity of the Dirichlet trace operator
$\gamma_{\widehat\Gamma}:H^1(\widehat\D)
\rightarrow H^{\frac{1}{2}}(\widehat\Gamma)$.
The frequency estimate \eqref{eq:K_explicit_frequency} follows
directly from \eqref{eq:double_layer_main_frequency}.

\end{proof}

\begin{theorem}\label{thm:KpHlm}
Let $\mathfrak T$ be a set of domain transformations
satisfying Assumptions~\ref{assump:transformations_T}
and~\ref{assump:domain_T_2}.
Then there exists $\delta_0>0$, depending only on $\mathfrak T$, such that,
for every $\kappa>0$ and $0<\delta<\delta_0$, the map
\begin{equation*}\label{eq:map_hol_Kadj}
	\mathcal{K}^{(\kappa)\prime}:
	\mathfrak{T}
	\rightarrow
	\mathscr{L}
	\left(
		H^{-\frac{1}{2}}(\widehat\Gamma),
		H^{-\frac{1}{2}}(\widehat\Gamma)
	\right):
	\br
	\mapsto
	\widehat{\mathsf K}^{(\kappa)\prime}_\br
\end{equation*}
admits a bounded holomorphic extension to $\mathfrak T_\delta$.
There exists $C_2=C_2(\mathfrak T,\delta_0)>0$ such that, for every
$\kappa_0>0$, there exists $C_1=C_1(\mathfrak T,\delta_0,\kappa_0)>0$
satisfying
\begin{equation*}
\sup_{\br\in\mathfrak T_\delta}
\norm{\widehat{\mathsf K}^{(\kappa)\prime}_\br}_
{\mathscr L(H^{-\frac{1}{2}}(\widehat\Gamma),H^{-\frac{1}{2}}(\widehat\Gamma))}
\leq C_1(1+\kappa)^{\frac{3}{2}}
e^{\frac{4\delta}{\pi R}\left(\kappa\delta_0C_2+\frac{1}{2}\log(1+\kappa)\right)}
\end{equation*}
for all $\kappa\geq\kappa_0$, $0<R<\delta_0$, and
$0<\delta<\min\{R,\delta_0-R\}$.
\end{theorem}

\begin{proof}
Recalling that double layer and adjoint double layer operator are
adjoint with respect to the duality product \eqref{eq:transport_boundary} for each $\br\in\mathfrak T$,
we obtain, for each $\br\in\mathfrak T$,
$\widehat\phi\in H^{-\frac{1}{2}}(\widehat\Gamma)$, and
$\widehat\varphi\in H^{\frac{1}{2}}(\widehat\Gamma)$, that
\begin{equation}\label{eq:dualityproofDBIO}
\dual{\widehat{\mathsf K}^{(\kappa)\prime}_\br\widehat\phi}
{\widehat\varphi}_
{H^{-\frac{1}{2}}(\widehat\Gamma)\times H^{\frac{1}{2}}(\widehat\Gamma),\br}
=\dual{\widehat\phi}
{\widehat{\mathsf K}^{(\kappa)}_\br\widehat\varphi}_
{H^{-\frac{1}{2}}(\widehat\Gamma)\times H^{\frac{1}{2}}(\widehat\Gamma),\br}.
\end{equation}
Choosing $\delta_0>0$ such that Theorem~\ref{thm:KHlm}
and Lemma~\ref{lem:bilinear_Gamma} hold implies that, for $\br\in\mathfrak T_{\delta_0}$, \eqref{eq:dualityproofDBIO} determines
$\widehat{\mathsf K}^{(\kappa)\prime}_\br\widehat\phi
\in H^{-\frac{1}{2}}(\widehat\Gamma)$ uniquely.
Theorem~\ref{thm:KHlm} and the argument of
\cite[Lemma~2.8]{HDPS26Laplace}, applied to \eqref{eq:dualityproofDBIO},
therefore yield a bounded holomorphic extension of
$\mathcal K^{(\kappa)\prime}$ in operator norm.
The continuity and inf-sup estimates
\eqref{eq:C_1_boundary} and~\eqref{eq:inf_sup_sesq_r_boundary} yield
\begin{equation*}
\norm{\widehat{\mathsf K}^{(\kappa)\prime}_\br}_
{\mathscr L(H^{-\frac{1}{2}}(\widehat\Gamma),H^{-\frac{1}{2}}(\widehat\Gamma))}
\leq
\frac{C_1(\mathfrak T,\delta_0)}{C_2(\mathfrak T,\delta_0)}
\norm{\widehat{\mathsf K}^{(\kappa)}_\br}_
{\mathscr L(H^{\frac{1}{2}}(\widehat\Gamma),H^{\frac{1}{2}}(\widehat\Gamma))},
\end{equation*}
from which the frequency explicit operator norm bound follows with
Theorem~\ref{thm:KHlm}.
\end{proof}

\subsection{Shape Holomorphy of the Hypersingular BIO}\label{sec:hypersingular_op_hol}
We recall from
\eqref{eq:transported_calderon_bios} the definition of
$\widehat{\mathsf{W}}^{(\kappa)}_\br$, with $\mathsf{W}^{(\kappa)}_{\Gamma_\br}$ defined in
\eqref{eq:def_hypersingular}. We establish shape holomorphy of the map
$\mathfrak{T}\ni\br\mapsto\widehat{\mathsf{W}}^{(\kappa)}_\br\in
\mathscr{L}\bigl(H^{\frac{1}{2}}(\widehat{\Gamma}),
H^{-\frac{1}{2}}(\widehat{\Gamma})\bigr)$.

\begin{theorem}\label{thm:WHlm}
Let $\mathfrak T$ be a set of domain transformations
satisfying Assumptions~\ref{assump:transformations_T}
and~\ref{assump:domain_T_2}.
Then there exists $\delta_0=\delta_0(\mathfrak T)>0$ such that,
for every $\kappa>0$ and $0<\delta<\delta_0$, the map
\begin{equation*}\label{eq:map_hol_W}
	\mathcal{W}^{(\kappa)}:
	\mathfrak{T}
	\rightarrow
	\mathscr{L}
	\left(
		H^{\frac{1}{2}}(\widehat\Gamma),
		H^{-\frac{1}{2}}(\widehat\Gamma)
	\right):
	\br
	\mapsto
	\widehat{\mathsf W}^{(\kappa)}_\br
\end{equation*}
admits a bounded holomorphic extension to $\mathfrak T_\delta$.
Moreover, there exists $C_2=C_2(\mathfrak T,\delta_0)>0$ such that, for every
$\kappa_0>0$, there exists $C_1=C_1(\mathfrak T,\delta_0,\kappa_0)>0$
satisfying
\begin{equation*}
\sup_{\br\in\mathfrak T_\delta}
\norm{\widehat{\mathsf W}^{(\kappa)}_\br}_
{\mathscr L(H^{\frac{1}{2}}(\widehat\Gamma),H^{-\frac{1}{2}}(\widehat\Gamma))}
\leq C_1(1+\kappa)^{\frac{d+1}{2}}
e^{\frac{4\delta}{\pi R}\left(\kappa\delta_0C_2+2\log(1+\kappa)\right)}
\end{equation*}
for all $\kappa\geq\kappa_0$, $0<R<\delta_0$, and
$0<\delta<\min\{R,\delta_0-R\}$.
\end{theorem}

\begin{proof}
We proceed in several steps.

\textbf{Step 1: Shape holomorphy.}
We recall that
$\mathsf W^{(\kappa)}_{\Gamma_\br}
=-\gamma_{1,\Gamma_\br}\mathsf D^{(\kappa)}_{\Gamma_\br}$,
and also that the double layer potential solves the homogeneous
Helmholtz equation, cf.~\eqref{eq:single_double_sol}.
The Neumann trace and the double layer potential are transported
as one composition according to
$\widehat{\mathsf W}^{(\kappa)}_\br
=-\tau_\br\gamma_{1,\Gamma_\br}
\mathsf D^{(\kappa)}_{\Gamma_\br}\tau_\br^{-1}$.
Choosing $Z_{\D_\br}=\tau_\br^{-1}Z_{\widehat\D}\tau_\br$ in
\eqref{eq:defneumanntrace}, and recalling
\eqref{eq:transport_boundary}, we obtain, for each
$\br\in\mathfrak T$ and all
$\widehat\phi,\widehat\varphi\in H^{\frac{1}{2}}(\widehat\Gamma)$,
\begingroup
\small
\begin{equation}\label{eq:Wdualityrep}
\begin{aligned}
\dual{\widehat{\mathsf W}^{(\kappa)}_\br\widehat\phi}
{\widehat\varphi}_
{H^{-\frac{1}{2}}(\widehat\Gamma)\times H^{\frac{1}{2}}(\widehat\Gamma),\br}
&=-\int_{\widehat\D}\!
(D\br)^{-\top}\widehat\nabla
\widehat{\mathsf D}^{(\kappa)}_\br\widehat\phi
\!\cdot\!(D\br)^{-\top}\widehat\nabla Z_{\widehat\D}\widehat\varphi
J_\br\,\dd\widehat\x\\
&\quad+\kappa^2\int_{\widehat\D}\widehat{\mathsf D}^{(\kappa)}_\br\widehat\phi\,
Z_{\widehat\D}\widehat\varphi
J_\br\,\dd\widehat\x.
\end{aligned}
\end{equation}
\endgroup
Lemmas~\ref{lem:shape_hol_bilinear_forms} and~\ref{lem:bilinear_Gamma},
Theorem~\ref{thm:DLhlmrp}, and the argument of
\cite[Lemma~2.8]{HDPS26Laplace}, applied to \eqref{eq:Wdualityrep},
imply that there is $\delta_0>0$ such that $\mathcal W^{(\kappa)}$
admits a bounded holomorphic extension to $\mathfrak T_{\delta_0}$
in operator norm.

\textbf{Step 2. Crude frequency-explicit operator norm bound.}
Choose $\delta_0>0$ such that
\eqref{eq:kNewton}, \eqref{eq:newton_L2_explicit_bound}, and
\cite[Lemma~B.1, item~(i)]{HDPS26Laplace} hold for all $0<\delta\leq\delta_0$.
It follows from \eqref{eq:inf_sup_sesq_r_boundary} combined with \eqref{eq:Wdualityrep} that
\begin{equation}
\label{eq:WprevBound}
\begin{aligned}
\norm{\widehat{\mathsf W}^{(\kappa)}_\br\widehat\phi}_
{H^{-\frac{1}{2}}(\widehat\Gamma)}
&\leq
\frac{\norm{Z_{\widehat\D}}_
{\mathscr L(H^{\frac{1}{2}}(\widehat\Gamma),H^1(\widehat\D))}}
{C_2(\mathfrak T,\delta_0)}
\Bigl(
C_{\mathsf a}(\mathfrak T,\delta_0)
\norm{\widehat{\mathsf D}^{(\kappa)}_\br\widehat\phi}_{H^1(\widehat\D)}\\
&\quad+\kappa^2C_{\mathsf b}(\mathfrak T,\delta_0)
\norm{\widehat{\mathsf D}^{(\kappa)}_\br\widehat\phi}_{L^2(\widehat\D)}
\Bigr).
\end{aligned}
\end{equation}
Our aim is to convert this estimate into a crude frequency-explicit operator norm bound for $\widehat{\mathsf W}^{(\kappa)}_\br\colon H^{1/2}(\widehat{\Gamma})\to H^{-1/2}(\widehat{\Gamma})$. With an operator norm bound for $\widehat{\mathsf D}^{(\kappa)}_\br\colon H^{1/2}(\widehat{\Gamma})\to H^1(\widehat{\D})$ provided in Theorem~\ref{thm:DLhlmrp}, we obtain in the following an operator norm bound for $\widehat{\mathsf D}^{(\kappa)}_\br\colon H^{1/2}(\widehat{\Gamma})\to L^2(\widehat{\D})$.

We proceed by showing bounds for the single layer and Newton potentials and combine them with \eqref{eq:decomp_double_layer}. First, notice that by
Theorem~\ref{thm:holo_single_layer_pot} we obtain
\begin{align*}
    \norm{\widehat{\mathsf{S}}^{(\kappa)}_\br
    }_{\mathscr{L}\left(H^{-\frac{1}{2}}(\widehat{\Gamma}),L^{2}(\widehat{\D})\right)}
    \leq C(\mathfrak{T},\delta_0)
    \norm{\widehat{\mathsf{N}}^{(\kappa)}_\br
    }_{\mathscr{L}\left(\widetilde{H}^{-1}(\widehat{\D}),L^{2}(\widehat{\D})\right)},
\end{align*}
then, \eqref{eq:kNewton} yields 
\begin{align*}
    \sup_{\br\in\mathfrak T_\delta}\norm{\widehat{\mathsf S}^{(\kappa)}_\br}_
    {\mathscr L(H^{-\frac{1}{2}}(\widehat\Gamma),L^2(\widehat\D))}
    &\leq C(\mathfrak T,\delta_0)
    \sup_{\br\in\mathfrak T_\delta}\norm{\widehat{\mathsf N}^{(\kappa)}_\br}_
    {\mathscr L(\widetilde H^{-1}(\widehat\D),L^2(\widehat\D))}\\
    &\lesssim C(\mathfrak T,\delta_0)
    \sup_{\br\in\mathfrak T_\delta}\norm{\widehat{\mathsf N}^{(\kappa)}_{\br}}_
    {\mathscr L(L^2(\widehat\D),H^1(\widehat\D))}\\
    &\leq C_1e^{\kappa\delta C_2}
\end{align*}
for some constants $C_1=C_1(\mathfrak T,\delta_0,\kappa_0)>0$ and
$C_2=C_2(\mathfrak T,\delta_0)>0$.
Second, \eqref{eq:newton_L2_explicit_bound} implies
\begin{equation*}
\sup_{\br\in\mathfrak T_\delta}
\norm{\widehat{\mathsf N}^{(\kappa)}_\br}_
{\mathscr L(L^2(\widehat\D),L^2(\widehat\D))}
\leq C_1\kappa^{\frac{d-3}{2}}e^{\kappa\delta C_2}.
\end{equation*}
Combining these two bounds and \cite[Theorem~3.4]{HDPS26Laplace}, with 
\eqref{eq:decomp_double_layer} yields
\begin{equation}\label{eq:DPL2frequency}
\sup_{\br\in\mathfrak T_\delta}
\norm{\widehat{\mathsf D}^{(\kappa)}_\br}_
{\mathscr L(H^{\frac{1}{2}}(\widehat\Gamma),L^2(\widehat\D))}
\leq C_1(1+\kappa)^{\frac{d+1}{2}}e^{\kappa\delta C_2}.
\end{equation}
Inserting the latter estimate and \eqref{eq:double_layer_preliminary_bound} into \cref{eq:WprevBound} implies that
\begin{equation}\label{eq:hypersingular_preliminary_bound}
\sup_{\br\in\mathfrak T_{\delta_0}}
\norm{\widehat{\mathsf W}^{(\kappa)}_\br}_
{\mathscr L(H^{\frac{1}{2}}(\widehat\Gamma),H^{-\frac{1}{2}}(\widehat\Gamma))}
\leq C_1(1+\kappa)^{\frac{d+5}{2}}e^{\kappa\delta_0 C_2},
\quad\kappa\geq\kappa_0.
\end{equation} 

\textbf{Step 3. Improved frequency-explicit operator norm bound.}
The limiting entity for the frequency-explicit bound \cref{eq:hypersingular_preliminary_bound}
is \eqref{eq:DPL2frequency}, used to estimate the second term in \eqref{eq:WprevBound}. To improve the bound we consider real-valued domain
transformations $\br\in\mathfrak T_{\delta_0}\cap
W^{1,\infty}(\widehat\D;\IR^d)$, implying that the $\D_\br$ are bounded
and Lipschitz according to \cite[Lemma~B.1, item~(i)]{HDPS26Laplace}.

Let $\br\in\mathfrak T_{\delta_0}\cap
W^{1,\infty}(\widehat\D;\IR^d)$; then from \eqref{eq:defneumanntrace}, for every $\phi,\psi\in H^{\frac{1}{2}}(\Gamma_\br)$, we have
\begin{equation}\label{eq:hypersingular_harmonic_weak_identity}
\begin{aligned}
\dual{\mathsf W^{(\kappa)}_{\Gamma_\br}\phi}{\psi}_
{H^{-\frac{1}{2}}(\Gamma_\br)\times H^{\frac{1}{2}}(\Gamma_\br)}
&=-\int_{\D_\br}\nabla \mathsf D^{(\kappa)}_{\Gamma_\br}\phi\cdot\nabla  Z_{\D_\br}\psi\,\dd\x\\
&\quad+\kappa^2\int_{\D_\br}(\mathsf D^{(\kappa)}_{\Gamma_\br}\phi)( Z_{\D_\br}\psi)\,\dd\x,
\end{aligned}
\end{equation}
Similarly to \eqref{eq:reformulatetoB}, we obtain
\begin{equation}
\label{eq:secondtermW}
\begin{aligned}
&\int_{\D_\br}(\mathsf{D}^{(\kappa)}_{\Gamma_\br}\phi)(Z_{\D_\br}\psi)\,\dd\x\\
&\hspace{2cm}=
\int_{\D_\br}\int_{\Gamma_\br}\bm\nu_\br(\y)\cdot
\nabla_\y G^{(\kappa)}(\x-\y)\phi(\y)\,\dd s_\y(Z_{\D_\br}\psi)(\x)\,\dd\x\\
&\hspace{2cm}=
\int_{\Gamma_\br}\phi(\y)\bm\nu_\br(\y)\cdot\int_{\D_\br}
\nabla_\y G^{(\kappa)}(\y-\x)(Z_{\D_\br}\psi)(\x)\,\dd\x\,\dd s_\y\\
&\hspace{2cm}=
\int_{\Gamma_\br}\phi(\gamma_{1,\Gamma_\br}\mathsf N^{(\kappa)}_\br Z_{\D_\br}\psi)\,\dd s_\y,
\end{aligned}
\end{equation}
as, for any $h\in H^1(\D_\br)$, it holds $\mathsf N^{(\kappa)}_\br h
\in H^2_{\mathrm{loc}}(\mathbb R^d)$,
$\nabla \mathsf N^{(\kappa)}_\br h\in H^1_{\mathrm{loc}}(\mathbb R^d;\IC^d)$,
and therefore $\gamma_{1,\Gamma_\br}\mathsf N^{(\kappa)}_\br h
=\bm\nu_\br\cdot\gamma_{\Gamma_\br}\nabla \mathsf N^{(\kappa)}_\br h\in L^2(\Gamma_\br)$.
As in \cite[Eq.~(108)]{LabarcaHiptmair2025}, we obtain
\begin{equation}\label{eq:hypersingular_newton_normal_identity}
\gamma_{1,\Gamma_\br}\mathsf N^{(\kappa)}_\br h
=\bm\nu_\br\cdot\left(
\gamma_{\Gamma_\br}\mathsf N^{(\kappa)}_\br\nabla h
-\mathsf V^{(\kappa)}_{\Gamma_\br}
(\bm\nu_\br\gamma_{\Gamma_\br}h)\right),
\end{equation}
with operators $\mathsf N^{(\kappa)}_\br$ and
$\mathsf V^{(\kappa)}_{\Gamma_\br}$ acting component-wise.

To estimate the latter term, let $B_L$ be as in \eqref{eq:choice_B_L}.
For $q\in L^2(\mathbb R^d)$ supported in $B_L$ and $\kappa\geq\kappa_0$,
the free resolvent estimate \eqref{eq:dual_newton_free_H2} and the
multiplicative trace inequality, cf.\ \cite[Theorem~1.5.1.10]{Grisvard}, yield
\begin{equation}\label{eq:hypersingular_newton_Dirichlet_trace}
	\norm{\gamma_{\Gamma_\br}\mathsf N_\br^{(\kappa)} q}_{L^2(\Gamma_\br)}
	\leq 
	C_{\mathrm{tr}}(\mathfrak{T},\delta_0)C_{\mathrm{res}}(L,\kappa_0)
	\kappa^{-\frac{1}{2}}\norm{q}_{L^2(\mathbb R^d)},
\end{equation}
where $C_{\mathrm{tr}}(\mathfrak{T},\delta_0)>0$ depends only on $\mathfrak{T}$
and $\delta_0$.
Moreover, \cite[Theorem~3.3]{CGL09} implies
\begin{equation}\label{eq:hypersingular_single_layer_L2_bound}
	\norm{\mathsf V^{(\kappa)}_{\Gamma_\br}}_
	{\mathscr L(L^2(\Gamma_\br),L^2(\Gamma_\br))}
	\leq 
	C_{\mathrm S}(\mathfrak T,\delta_0)\kappa^{\frac{d-3}{2}},
\end{equation}
with fixed $C_{\mathrm S}(\mathfrak T,\delta_0)>0$.
Combining \eqref{eq:hypersingular_newton_normal_identity}, \eqref{eq:hypersingular_newton_Dirichlet_trace}, and
\eqref{eq:hypersingular_single_layer_L2_bound}, we obtain
\begin{equation}\label{eq:hypersingular_newton_normal_bound}
\begin{aligned}
\norm{\gamma_{1,\Gamma_\br}\mathsf N^{(\kappa)}_\br h}_{L^2(\Gamma_\br)}
&\leq C\left(\kappa^{-\frac{1}{2}}\norm{\nabla h}_{L^2(\D_\br)}
+\kappa^{\frac{d-3}{2}}\norm{\gamma_{\Gamma_\br} h}_{L^2(\Gamma_\br)}\right)\\
&\leq C(\mathfrak T,\delta_0,\kappa_0)
\kappa^{\frac{d-3}{2}}\norm{h}_{H^1(\D_\br)}.
\end{aligned}
\end{equation}
Using the previous bound in \eqref{eq:secondtermW} gives
\begin{equation}\label{eq:hypersingular_harmonic_mass_bound}
\begin{aligned}   
\left|\int_{\D_\br}(\mathsf D^{(\kappa)}_{\Gamma_\br}\phi)(Z_{\D_\br}\psi)\,\dd\x\right|
& \leq C(\mathfrak T,\delta_0,\kappa_0)\kappa^{\frac{d-3}{2}}
\norm{\phi}_{L^2(\Gamma_\br)}
\norm{Z_{\D_\br} \psi}_{H^{1}(\D_\br)}, \\
&\leq C(\mathfrak T,\delta_0,\kappa_0)\kappa^{\frac{d-3}{2}}
\norm{\phi}_{H^{\frac{1}{2}}(\Gamma_\br)}
\norm{\psi}_{H^{\frac{1}{2}}(\Gamma_\br)}.
\end{aligned}
\end{equation}
Combining \eqref{eq:hypersingular_harmonic_weak_identity}, \eqref{eq:double_layer_real_bound}, and \eqref{eq:hypersingular_harmonic_mass_bound} yields
\[
\begin{aligned}
\norm{\mathsf W^{(\kappa)}_{\Gamma_\br}}_
{\mathscr L(H^{\frac{1}{2}}(\Gamma_\br),H^{-\frac{1}{2}}(\Gamma_\br))}
&\leq C(\mathfrak T,\delta_0,\kappa_0)
\left((1+\kappa)^{\frac{3}{2}}+\kappa^{\frac{d+1}{2}}\right)\\
&\leq
C(\mathfrak T,\delta_0,\kappa_0)(1+\kappa)^{\frac{d+1}{2}},
\end{aligned}
\]
from which pulling back to the reference domain yields, for $\kappa\geq\kappa_0$,
\begin{equation}\label{eq:hypersingular_real_bound}
	\sup_{\br\in\mathfrak T_{\delta_0}\cap W^{1,\infty}(\widehat\D;\IR^d)}
	\norm{\widehat{\mathsf W}^{(\kappa)}_\br}_
	{\mathscr L(H^{\frac{1}{2}}(\widehat\Gamma),H^{-\frac{1}{2}}(\widehat\Gamma))}
	\leq 
	C_0(1+\kappa)^{\frac{d+1}{2}}.
\end{equation}
The assertion follows from applying Corollary~\ref{cor:abstract_frequency_bound} to
\eqref{eq:hypersingular_real_bound} and
\eqref{eq:hypersingular_preliminary_bound} and using {\bf Step~1}.
\end{proof}

\begin{remark}[Comparison with known wavenumber-explicit bounds]
\label{rem:real_frequency_comparison}
We compare the orders in $\kappa$ obtained in
Theorems~\ref{thm:Vhlmrp}--\ref{thm:WHlm} for real-valued
transformations, i.e., in the limit $\delta\downarrow0$ for fixed $R$,
with those available in the literature for a fixed boundary
$\Gamma=\partial\D$ of a bounded Lipschitz domain.

For the single layer BIO, our upper bound of order $\kappa^{\frac{d-1}{2}}$
coincides with the bound in Corollary~1.2, item~(i), of \cite{GrahamLohndorfMelenkSpence2015},
which, to our knowledge, is the best bound available on Lipschitz boundaries.
For domains with connected exterior, Theorem~1.5, items~(a) and~(c), and Corollary~1.9 in
\cite{GalkowskiSpence2019} provide upper bounds of orders
$\kappa^{\frac12}\log\kappa$, $\kappa^{\frac13}\log\kappa$, and
$\kappa^{\frac13}$ for piecewise smooth, piecewise curved, and convex
$C^\infty$ curved boundaries, respectively. For piecewise smooth
boundaries in $d=2$, however, the Lipschitz bound $\kappa^{\frac12}$
above avoids the logarithmic factor.

For the double layer BIO and its adjoint, we obtain an upper bound
of order $\kappa^{\frac32}$ in
$H^{\frac12}(\Gamma)\to H^{\frac12}(\Gamma)$ and
$H^{-\frac12}(\Gamma)\to H^{-\frac12}(\Gamma)$ on Lipschitz boundaries.
For comparison, Theorem~3.5 of \cite{CGL09} gives an
$L^2(\Gamma)\to L^2(\Gamma)$ bound of order $\kappa^{\frac{d-1}{2}}$
on Lipschitz boundaries. Under additional regularity, bounds are
available for the stronger mapping
$H^{-\frac12}(\Gamma)\to H^{\frac12}(\Gamma)$: Corollary~1.2, item~(ii), in
\cite{GrahamLohndorfMelenkSpence2015} gives an upper bound of order
$\kappa^{\frac{d+1}{2}}$ for $C^2$ boundaries, while, for domains with connected exterior,
Theorem~1.5, items~(b) and~(c), and Corollary~1.9 in
\cite{GalkowskiSpence2019} provide upper bounds of orders
$\kappa^{\frac54}\log\kappa$, $\kappa^{\frac76}\log\kappa$, and
$\kappa$ for globally $C^{2,\alpha}$ piecewise smooth, globally
$C^{2,\alpha}$ piecewise curved, and convex $C^\infty$ curved
boundaries, respectively, for some $\alpha\in(0,1]$. 

For the hypersingular BIO, we obtain an upper bound of order
$\kappa^{\frac{d+1}{2}}$ in
$H^{\frac12}(\Gamma)\to H^{-\frac12}(\Gamma)$ on Lipschitz boundaries.
Theorem~4.7 in \cite{GalkowskiMarchandSpence2022} implies upper bounds
of orders $\kappa\log\kappa$ for piecewise smooth boundaries and
$\kappa$ for convex $C^\infty$ curved boundaries in the same spaces.
\end{remark}

\section{Concluding Remarks}
\label{sec:concluding_remarks}
We have established shape holomorphy of the Calder\'on projector for
the Helmholtz equation on the boundaries of bounded Lipschitz domains
in dimensions $d=2,3$, obtained from a fixed reference domain by
bijective and bi-Lipschitz transformations
(Theorem~\ref{thm:Calderon_Projector}). This extends the analysis of
\cite{HDPS26Laplace} from the Laplace to the Helmholtz case. As in
that work, the four BIOs are expressed through the Newton potential,
the Dirichlet and Neumann traces and their adjoints, together with
the corresponding solution extension. Shape holomorphy of the BIOs follows from
the corresponding properties of these operators.

In addition, we provide a complete wavenumber-explicit analysis of the
bounds for the four Helmholtz BIOs on Lipschitz boundaries. The domain of
complex holomorphic extension can be chosen independently of $\kappa>0$.
On a domain of size $\delta>0$, the bounds contain algebraic factors in
$\kappa$ multiplied by $\exp(C\kappa\delta)$, with $C$ independent of
$\kappa$ and $\delta$. Choosing $\delta$ proportional to $\kappa^{-1}$
keeps the exponential factor uniformly bounded
and retains the algebraic growth of the real-domain bounds.
In particular, for real transformations, corresponding to $\delta=0$,
we obtain polynomial upper bounds of orders
$\kappa^{\frac{d-1}{2}}$, $\kappa^{\frac32}$, and
$\kappa^{\frac{d+1}{2}}$ for the single layer, double layer and its
adjoint, and hypersingular BIOs, respectively.
These estimates recover, to our knowledge, the best available single
layer bound and provide wavenumber-explicit bounds for the other three
BIOs on Lipschitz boundaries, as discussed in Remark~\ref{rem:real_frequency_comparison}.

Following \cite[Section~4]{HDPS26Laplace}, the proved results can be applied
to interior, exterior, screen, and transmission problems with
holomorphic transported data, whenever the corresponding boundary
integral formulations are well posed. Combined with the present shape
holomorphy results, the restriction and extension operators of
\cite[Sections~4.2--4.3]{HDPS26Laplace} provide the tools needed to
establish shape holomorphy of the multiple traces formulation
\cite{HJ12}, by combining the subdomain Calder\'on projectors and
coupling the Dirichlet and Neumann traces across the interfaces.

Composition with a holomorphic parametrization of the domain
transformations yields parametric holomorphy of the transported BIOs,
with bounds explicit in the wavenumber. Together with suitable
summability assumptions on the parametrization, this regularity allows
polynomial approximation of parametric maps with dimension-independent
convergence rates, as explained in \cite{CCS15}.

Current work concerns the use of these estimates in frequency-explicit
uncertainty quantification and model order reduction for scattering by
composite objects with parametrically varying interfaces. In particular,
the combination of proper orthogonal decomposition and neural networks
in \cite{WKHH25}, together with the fully discrete analysis of
\cite{DH2026}, provides a framework for constructing surrogates of
high-dimensional parameter-to-solution maps.

Furthermore, ongoing work also extends the present analysis to Maxwell BIOs, with the
aim of establishing shape holomorphy and bounds explicit in the frequency
and the size of the holomorphic extension for electromagnetic scattering
on Lipschitz boundaries. A complementary direction concerns scattering by
periodic structures with parametrically varying geometries, including
the diffraction gratings considered in \cite{AylwinPintoSilvaOelker2025}.
In this setting, we aim to complement the reduced basis approach to
uncertainty quantification developed in that work with frequency-explicit
holomorphy estimates for the corresponding BIOs.

\appendix
\section{Proof of Theorem~\ref{thm:shape_hol_Newton_L2}}
\label{app:proof_shape_hol_Newton_L2}

We establish shape holomorphy of 
the map $\br\mapsto\widehat{\mathsf N}^{(\kappa)}_\br$
by adapting the argument of \cite[Theorem~4.11]{henriquez2024domain}
and applying \cite[Theorems~4.2 and~4.3]{henriquez2024domain}.

\textbf{Step 1: Auxiliary geometric bounds.}
We provide some auxiliary geometric bounds which will be useful later on. To simplify presentation we define
\begin{equation*}
	\boldsymbol\delta_\br\coloneqq\br(\widehat\x)-\br(\widehat\y)
	\quad
	\text{and}
	\quad
	\rho_\br
	\coloneqq
	\sqrt{\boldsymbol\delta_\br\cdot\boldsymbol\delta_\br},
\end{equation*}
for each $\br\in\mathfrak{T}_\delta$, with $\delta>0$
to be specified below, and $(\widehat\x,\widehat\y)\in (\widehat\D\times\widehat\D)^\star$. Here, $\sqrt{\cdot}$ denotes the principal branch of the square root.

As a consequence of Lemma~\ref{lem:lipschitz_preservation}, there
exist $\delta=\delta(\mathfrak{T})>0$ and $\eta=\eta(\mathfrak{T})>0$
such that
\begin{equation*}
	\operatorname{Re}
	\left\{
		\boldsymbol\delta_\br\cdot\boldsymbol\delta_\br
	\right\}
	\geq\eta\norm{\widehat\x-\widehat\y}^2
	\quad\text{and}\quad
	\operatorname{Re}\rho_\br
	\geq\sqrt\eta\norm{\widehat\x-\widehat\y}>0.
\end{equation*}
To control the imaginary part of $\rho_\br$, consider $\widetilde{\br}\in\mathfrak{T}$ with
$\norm{\br-\widetilde{\br}}_{W^{1,\infty}(\widehat\D;\IC^d)}<\delta$.
Since $\widetilde{\br}$ is real, we obtain
\begin{equation*}
	\left(
		\boldsymbol\delta_\br-\boldsymbol\delta_{\widetilde\br}
	\right)
	\cdot
	\big(
		\boldsymbol\delta_\br
       	+
		\boldsymbol\delta_{\widetilde\br}
    	\big)
	=
	\boldsymbol\delta_\br\cdot\boldsymbol\delta_\br
	-
	\norm{\boldsymbol\delta_{\widetilde\br}}^2,
\end{equation*}
and
\begin{align*}
	\snorm{\operatorname{Im}
		\big\{
		\boldsymbol\delta_\br\cdot\boldsymbol\delta_\br
		\big\}
	}
	\leq
	2\delta
	\left(
		\operatorname{Lip}(\br;\widehat\D)
       	+
		\operatorname{Lip}(\widetilde{\br};\widehat\D)
    \right)
    \norm{\widehat\x-\widehat\y},
\end{align*}
due to $\norm{\boldsymbol\delta_\br-\boldsymbol\delta_{\widetilde\br}}
\leq2\norm{\br-\widetilde\br}_{L^\infty(\widehat\D;\mathbb C^d)}
<2\delta$.
Finally, observing that $\operatorname{Im}z=2\operatorname{Re}\sqrt z
\operatorname{Im}\sqrt z$ for all $\operatorname{Re}z>0$, we obtain with Lemma~\ref{lem:lipschitz_preservation} that
\begin{align}
\label{eq:fbound}
\left|\operatorname{Im}
\rho_\br\right|
\leq\delta C_2(\mathfrak{T},\delta),
\end{align}
with $C_2(\mathfrak T,\delta)$ defined in
\eqref{eq:kernel_exponential_constant}.

\textbf{Step 2: Shape holomorphy.}
To show shape holomorphy, we define
\begin{equation*}
	\mathsf g^{(\kappa)}_{\br,\mathbb C}
	(\widehat\x,\widehat\y)
	\coloneqq
	\Phi_d^{(\kappa)}(\rho_\br)J_\br(\widehat\y),
	\quad
	(\widehat\x,\widehat\y)\in
	(\widehat\D\times\widehat\D)^\star,
\end{equation*}
for each $\br \in \mathfrak{T}_\delta$, with $\delta>0$ as in Step~1,
and where, for $\operatorname{Re}\{\rho\}>0$,
\begin{equation}\label{eq:newton_radial_helmholtz_kernel}
	\Phi_2^{(\kappa)}(\rho)\coloneqq
	\frac{\imath}{4}\mathrm H_0^{(1)}(\kappa \rho)
	\quad
	\text{and}
	\quad
	\Phi_3^{(\kappa)}(\rho)\coloneqq
	\frac{e^{\imath\kappa \rho}}{4\pi \rho}.
\end{equation}
We denote the corresponding integral operator \eqref{eq:Nhatexplicit} by $\widehat{\mathsf N}^{(\kappa)}_{\br,\mathbb C}$.
In analogy to \cite[Theorem~4.11, pp.~95--96]{henriquez2024domain} we obtain that
there exists $\delta>0$ such that the map
$
\mathfrak T\ni\br
\mapsto
\widehat{\mathsf N}^{(\kappa)}_{\br}
\in\mathscr L(L^2(\widehat\D),H^1(\widehat\D))
$
admits a bounded holomorphic extension into $\mathfrak T_\delta$.

\textbf{Step 3: Crude frequency-explicit operator norm bound.}
Following the proofs of \cite[Theorem~4.2]{henriquez2024domain} and \cite[Theorem~4.3]{henriquez2024domain}, frequency-explicit bounds in modulus of the kernel provide us with frequency-explicit bounds in operator norm. Applied to the setting of this article, observe that for almost every $\widehat\x\in\widehat\D$ with
$\widehat\x\ne\widehat\y$, it holds
\begin{equation*}\label{eq:newton_gradient_formula}
\widehat\nabla_{\widehat\x}
\mathsf g^{(\kappa)}_{\br,\mathbb C}
=J_\br(\widehat\y)D\br(\widehat\x)^\top
\boldsymbol\delta_\br(\widehat\x,\widehat\y)
\frac{(\Phi_d^{(\kappa)})'(\rho_\br(\widehat\x,\widehat\y))}
{\rho_\br(\widehat\x,\widehat\y)},
\end{equation*}
where
\begin{equation}\label{eq:radial_first_derivatives}
(\Phi_2^{(\kappa)})'(\rho)
=-\frac{\imath\kappa\mathrm H_1^{(1)}(\kappa\rho)}{4}
\quad\text{and}\quad
(\Phi_3^{(\kappa)})'(\rho)
=\frac{e^{\imath\kappa\rho}(\imath\kappa\rho-1)}{4\pi\rho^2}.
\end{equation}
For $d=3$, \eqref{eq:fbound} yields
\begin{equation}\label{eq:bound_singular_kernel_3d}
\begin{aligned}
|\mathsf g^{(\kappa)}_{\br,\mathbb C}
(\widehat\x,\widehat\y)|
&\leq C(\mathfrak T,\delta)
e^{\kappa\delta C_2(\mathfrak T,\delta)}
\norm{\widehat\x-\widehat\y}^{-1},\\
\norm{\widehat\nabla_{\widehat\x}
\mathsf g^{(\kappa)}_{\br,\mathbb C}
(\widehat\x,\widehat\y)}
&\leq C(\mathfrak T,\delta)(1+\kappa)
e^{\kappa\delta C_2(\mathfrak T,\delta)}
\norm{\widehat\x-\widehat\y}^{-2},
\end{aligned}
\end{equation}
whereas for $d=2$ the Hankel function bounds for $\operatorname{Re}\{z\}>0$
\begin{equation}\label{eq:hankel_bounds}
\begin{aligned}
\left|\mathrm H_0^{(1)}(z)\right|
&\leq
\left|e^{\imath z}\right|
\left(\frac{2}{\pi\operatorname{Re}\{z\}}\right)^{\frac{1}{2}},\\
\left|\mathrm H_1^{(1)}(z)\right|
&\leq
\left|e^{\imath z}\right|
\left(
\left(\frac{2}{\pi\operatorname{Re}\{z\}}\right)^{\frac{1}{2}}
+\frac{2}{\pi\operatorname{Re}\{z\}}
\right),
\end{aligned}
\end{equation}
cf.\ \cite[Theorem~4.10, p.~93]{henriquez2024domain} and
\cite[Theorem~4.11, p.~95]{henriquez2024domain}, \eqref{eq:fbound}, and \eqref{eq:radial_first_derivatives} yields
\begin{equation}\label{eq:bound_singular_kernel_2d}
\begin{aligned}
|\mathsf g^{(\kappa)}_{\br,\mathbb C}
(\widehat\x,\widehat\y)|
&\leq C(\mathfrak T,\delta)\kappa^{-\frac{1}{2}}
e^{\kappa\delta C_2(\mathfrak T,\delta)}
\norm{\widehat\x-\widehat\y}^{-\frac{1}{2}},\\
\norm{\widehat\nabla_{\widehat\x}
\mathsf g^{(\kappa)}_{\br,\mathbb C}
(\widehat\x,\widehat\y)}
&\leq C(\mathfrak T,\delta)(1+\kappa^{\frac{1}{2}})
e^{\kappa\delta C_2(\mathfrak T,\delta)}
\norm{\widehat\x-\widehat\y}^{-1}.
\end{aligned}
\end{equation}

Following the proof of \cite[Theorem 4.2]{henriquez2024domain},
\eqref{eq:bound_singular_kernel_3d} and \eqref{eq:bound_singular_kernel_2d}
give bounds on the operator in $L^1(\widehat\D)$ and $L^\infty(\widehat\D)$.
The Riesz--Thorin interpolation theorem then yields
\begin{equation}\label{eq:newton_L2_explicit_bound}
\norm{\widehat{\mathsf N}^{(\kappa)}_{\br,\mathbb C}}_
{\mathscr L(L^2(\widehat\D),L^2(\widehat\D))}
\leq C(\mathfrak T,\delta)\kappa^{\frac{d-3}{2}}
e^{\kappa\delta C_2(\mathfrak T,\delta)},
\end{equation}
whereas following the arguments in the proof of
\cite[Theorem~4.3]{henriquez2024domain} yields
\begin{equation*}
\norm{\widehat\nabla\widehat{\mathsf N}^{(\kappa)}_{\br,\mathbb C}}_
{\mathscr L(L^2(\widehat\D),L^2(\widehat\D;\IC^d))}
\leq C(\mathfrak T,\delta)(1+\kappa^{\frac{d-1}{2}})
e^{\kappa\delta C_2(\mathfrak T,\delta)}.
\end{equation*}
This implies
\begin{equation}\label{eq:newton_reference_preliminary_bound}
\sup_{\br\in\mathfrak T_\delta}
\norm{\widehat{\mathsf N}^{(\kappa)}_\br}_{
\mathscr L(L^2(\widehat\D),H^1(\widehat\D))}
\leq C(\mathfrak T,\delta)e^{\kappa\delta C_2(\mathfrak T,\delta)}
\begin{cases}
1+\kappa^{\frac{1}{2}}+\kappa^{-\frac{1}{2}},&d=2,\\
1+\kappa,&d=3.
\end{cases}
\end{equation}

\textbf{Step 4: Frequency-explicit operator norm bound.}
Aiming to apply Corollary~\ref{cor:abstract_frequency_bound}, we derive in the following frequency-explicit operator norm bounds for real domain transformation. This allows us to work in the physical domain and to use classical resolvent estimates.

To this end, let $\delta_0>0$ be such that Lemma~\ref{lem:lipschitz_preservation} and \cite[Lemma~B.1, item~(i)]{HDPS26Laplace} hold on $\mathfrak T_{\delta_0}$, implying that for each
$\br\in\mathfrak T_{\delta_0}\cap W^{1,\infty}(\widehat\D;\IR^d)$,
the map $\br:\widehat\D\rightarrow\D_\br$ is bi-Lipschitz and
$\D_\br$ is a bounded Lipschitz domain.
Let $L>0$ and $B_L=B(\boldsymbol0,L)$ be as in
\eqref{eq:choice_B_L} in the proof of Theorem~\ref{thm:shape_hol_Newton}.
For $f\in L^2(\mathbb R^d)$ supported in $B_L$, define
\begin{equation}\label{eq:free_newton_potential_definition}
    u(\x)\coloneqq
    \int_{\mathbb R^d}G^{(\kappa)}(\x-\y)f(\y)\,\dd\y,
    \quad \x\in\mathbb R^d.
\end{equation}
Now, for $\widehat\varphi\in L^2(\widehat\D)$, let $\varphi_\br$ be the
extension of $\tau_\br^{-1}\widehat\varphi$ by zero to $\mathbb R^d$,
and let $u_\br$ be given by
\eqref{eq:free_newton_potential_definition} with $f=\varphi_\br$.
Then $\widehat{\mathsf N}^{(\kappa)}_\br\widehat\varphi
=\tau_\br u_\br$. For every fixed $\kappa_0>0$, applying the free resolvent
estimate \eqref{eq:dual_newton_free_H2} to $f=\varphi_\br$ yields
\[
\norm{u_\br}_{H^1(\D_\br)}
\leq
\norm{u_\br}_{H^1(B_L)}
\leq C_{\mathrm{res}}(L,\kappa_0)
\norm{\varphi_\br}_{L^2(\mathbb R^d)},
\quad\kappa\geq\kappa_0.
\]
Pulling back to the reference domain and using Lemma~\ref{lem:lipschitz_preservation} and
Assumption~\ref{assump:transformations_T}, item~\normalfont{(iv)}, we obtain
\begin{equation}\label{eq:real_transported_newton_uniform_bound}
    \sup_{\br\in\mathfrak T_{\delta_0}\cap
    W^{1,\infty}(\widehat\D;\IR^d)}
    \norm{\widehat{\mathsf N}^{(\kappa)}_\br}_{\mathscr{L}(L^2(\widehat\D),H^1(\widehat\D))}
    \leq C_0,\quad\kappa\geq\kappa_0,
\end{equation}
for some $C_0=C_0(\mathfrak T,\delta_0,\kappa_0)\geq 1$.
\eqref{eq:newton_reference_preliminary_bound} and \eqref{eq:real_transported_newton_uniform_bound} combined with
Corollary~\ref{cor:abstract_frequency_bound} yield
\eqref{eq:newton_reference_improved_complex_bound}.
\hfill\proofbox

\section{Auxiliary Result for the Proof of Theorem~\ref{thm:shape_hol_Newton}}
\label{app:auxiliary_shape_hol_Newton}

We establish shape holomorphy of the operator $\mathsf B^{(\kappa)}_\br$
defined in \eqref{eq:B_operator}.

\begin{proposition}\label{prop:shape_hol_B_r_O}
Let $\mathfrak T$ be a set of domain transformations
satisfying Assumptions~\ref{assump:transformations_T}
and~\ref{assump:domain_T_2}. Then there exists $\delta_0=\delta_0(\mathfrak{T})>0$ such that, for each $\kappa>0$, the map
\begin{equation*}
\mathcal B^{(\kappa)}_{\mathcal O}:
\mathfrak T_{\mathcal O}\rightarrow
\mathscr L(H^{\frac{1}{2}}(\partial\mathcal O),H^1(\widehat\D)):
\br\mapsto\mathsf B^{(\kappa)}_\br
\end{equation*}
admits a bounded holomorphic extension into
$\mathfrak T_{\delta_0,\mathcal O}$. For every
$0<\delta\leq\delta_0$, there are $C_1=C_1(\mathfrak{T}_{\mathcal{O}},\delta_0)>0, C_2=C_2(\mathfrak{T}_{\mathcal{O}},\delta_0)>0$,  such that
\begin{equation}\label{eq:BopNorm}
\sup_{\br\in\mathfrak T_{\delta,\mathcal O}}
\norm{\mathsf B^{(\kappa)}_\br}_{\mathscr{L}(H^{\frac{1}{2}}(\partial\mathcal O),H^1(\widehat\D))}
\leq C_1
\left(1+\kappa^{\frac{d+1}{2}}\right)
e^{\kappa\delta C_2}.
\end{equation}
The precise expression for $C_2$ is given in \eqref{eq:kernel_exponential_constant} with $\mathcal O$ and
$\mathfrak T_{\mathcal O}$ in place of $\widehat\D$ and $\mathfrak T$.

Moreover, for each $\br\in\mathfrak T_{\delta_0,\mathcal O}\cap
W^{1,\infty}(\mathcal O;\IR^d)$ it holds
\begin{equation}\label{eq:B_operator_real_O_r}
(\mathsf B^{(\kappa)}_\br\phi)(\widehat\y)
=\int_{\partial\mathcal O}
J_{\mathrm S,\partial\mathcal O,\br}(\widehat\x)\,
\boldsymbol n_{\mathcal O_\br}(\br(\widehat\x))
\cdot\nabla G^{(\kappa)}
\bigl(\br(\widehat\x)-\br(\widehat\y)\bigr)
\phi(\widehat\x)\,\dd s_{\widehat\x},
\quad
\widehat\y\in\widehat\D,
\end{equation}
for all $\phi\in H^{\frac{1}{2}}(\partial\mathcal O)$.
\end{proposition}

\begin{proof}
We recall the notation for $\boldsymbol\delta_\br$ and $\rho_\br$ from Appendix~\ref{app:proof_shape_hol_Newton_L2}. The argument of the proof itself is analogous to the proof of \cite[Proposition E.1]{HDPS26Laplace}, but requires special care for the complex-valued kernel function.

\textbf{Step 1: Auxiliary geometric bounds.}
Observe that \cite[Lemma~A.2]{HDPS26Laplace} implies existence of $\delta,\eta>0$ such that
\[
\operatorname{Re}\{\boldsymbol\delta_\br\cdot\boldsymbol\delta_\br\}
\geq \eta\,\operatorname{dist}
\left(\partial\mathcal{O},\overline{\widehat{\D}}\right)^2>0,
\]
for each $\br\in\mathfrak{T}_{\delta,\mathcal{O}}$,
$\widehat{\x}\in\partial\mathcal{O}$, and
$\widehat{\y}\in\widehat{\D}$.
This implies
\begin{equation}\label{eq:rho_bounds_O}
\operatorname{Re}\{\rho_\br\}
=
\operatorname{Re}\{\sqrt{\boldsymbol{\delta}_\br\cdot\boldsymbol{\delta}_\br}\}
\geq\sqrt\eta\operatorname{dist}
\left(\partial\mathcal{O},\overline{\widehat{\D}}\right),
\end{equation}
and
\begin{equation*}
\left|\rho_\br\right|
\leq
\sup_{\br\in\mathfrak T_{\delta_0,\mathcal O}}
\sup_{(\widehat\x,\widehat\y)\in\partial\mathcal O\times\widehat\D}
\norm{\boldsymbol\delta_\br}.
\end{equation*}
The latter expression is finite due to the compactness of  $\mathfrak T_{\mathcal O}$ and implies
\begin{equation}\label{eq:Im_rho_bound_O}
\left|\operatorname{Im}
\{\rho_\br\}\right|
\leq\delta C_2(\mathfrak T_{\mathcal O},\delta_0)
\end{equation}
in analogy to \eqref{eq:fbound}. The preceding bounds imply that $\boldsymbol\delta_\br$ takes values in a compact subset
of $\{\z\in\mathbb C^d:\operatorname{Re}(\z\cdot\z)>0\}$ where
the holomorphic extension of $G^{(\kappa)}$ is given by
$G_{\mathbb C}^{(\kappa)}(\z)
=\Phi_d^{(\kappa)}((\z\cdot\z)^{\frac{1}{2}})$ with $\Phi_d^{(\kappa)}$ as in \eqref{eq:newton_radial_helmholtz_kernel}.

\textbf{Step 2: Frequency-explicit bounds for the kernel function.}
Observing that the same argument as in the proof of \cite[Lemma~2.4]{HDPS26Laplace} implies that $J_{\mathrm{S},\partial\mathcal{O},\br}\,
\boldsymbol{n}_{\mathcal{O}_\br}\circ\br$
admits a bounded holomorphic extension $\br\mapsto\boldsymbol c_{\partial\mathcal O,\br,\mathbb C}$,
the kernel of $\mathsf B^{(\kappa)}_\br$ extends for all $0<\delta\leq\delta_0$ and $\br\in\mathfrak T_{\delta,\mathcal O}$ to
\begin{equation}\label{eq:kernel_B_complex_O}
\mathsf t^{(\kappa)}_{\br,\mathbb C}(\widehat\x,\widehat\y)
\coloneqq\boldsymbol c_{\partial\mathcal O,\br,\mathbb C}(\widehat\x)
\cdot\nabla G_{\mathbb C}^{(\kappa)}
\bigl(\boldsymbol\delta_\br\bigr),
\end{equation}
with gradient in $\widehat\y$ reading
\begin{equation}\label{eq:gradient_kernel_B_complex_O}
\widehat\nabla_{\widehat\y}
\mathsf t^{(\kappa)}_{\br,\mathbb C}(\widehat\x,\widehat\y)
=-D\br(\widehat\y)^\top
D(\nabla G_{\mathbb C}^{(\kappa)})
\bigl(\boldsymbol\delta_\br\bigr)
\boldsymbol c_{\partial\mathcal O,\br,\mathbb C}
(\widehat\x).
\end{equation}
In the following we bound \eqref{eq:kernel_B_complex_O} and \eqref{eq:gradient_kernel_B_complex_O} by exploiting
\begin{align*}
\nabla G_{\mathbb C}^{(\kappa)}(\boldsymbol\delta_\br)
&=
(\Phi_d^{(\kappa)})'(\rho_\br)\frac{\boldsymbol\delta_\br}{\rho_\br}\\
D(\nabla G_{\mathbb C}^{(\kappa)})(\boldsymbol\delta_\br)
=&(\Phi_d^{(\kappa)})''(\rho_\br)
\frac{\boldsymbol\delta_\br\boldsymbol\delta_\br^\top}
{\boldsymbol\delta_\br\cdot\boldsymbol\delta_\br}
+\frac{(\Phi_d^{(\kappa)})'(\rho_\br)}{\rho_\br}
\left(I-\frac{\boldsymbol\delta_\br\boldsymbol\delta_\br^\top}
{\boldsymbol\delta_\br\cdot\boldsymbol\delta_\br}\right).\label{eq:hessian_G_O}
\end{align*}
Using \eqref{eq:radial_first_derivatives}, \eqref{eq:hankel_bounds} for $d=2$, and
\eqref{eq:rho_bounds_O}--\eqref{eq:Im_rho_bound_O} to estimate \cref{eq:kernel_B_complex_O} yields
\begin{equation}\label{eq:TNorm}
\begin{aligned}
|\mathsf t^{(\kappa)}_{\br,\mathbb C}(\widehat\x,\widehat\y)|
&\leq C(\mathfrak T_{\mathcal O},\delta_0)
(1+\kappa^{\frac{d-1}{2}})
e^{\kappa\delta C_2(\mathfrak T_{\mathcal O},\delta_0)}.
\end{aligned}
\end{equation}
To estimate \eqref{eq:gradient_kernel_B_complex_O}, observe that
differentiating \eqref{eq:newton_radial_helmholtz_kernel} and
$(\mathrm H_1^{(1)})'(z)=\mathrm H_0^{(1)}(z)
-\frac{\mathrm H_1^{(1)}(z)}{z}$ for $d=2$, c.f.\ \cite[equation~9.1.27]{AbramowitzStegun1964}, yields
\begin{equation}\label{eq:radial_second_derivatives}
\begin{aligned}
(\Phi_2^{(\kappa)})''(\rho_\br)
&=-\frac{\imath\kappa^2}{4}\mathrm H_0^{(1)}(\kappa\rho_\br)
+\frac{\imath\kappa}{4\rho_\br}\mathrm H_1^{(1)}(\kappa\rho_\br),\\
(\Phi_3^{(\kappa)})''(\rho_\br)
&=\frac{e^{\imath\kappa\rho_\br}}{4\pi}
\left(-\frac{\kappa^2}{\rho_\br}
-\frac{2\imath\kappa}{\rho_\br^2}+\frac{2}{\rho_\br^3}\right).
\end{aligned}
\end{equation}
For $d=2$, 
\eqref{eq:hankel_bounds} and
\eqref{eq:rho_bounds_O}--\eqref{eq:Im_rho_bound_O} yield
\[
\kappa^2\left|\mathrm H_0^{(1)}\!\left(
\kappa\rho_\br\right)\right|
\leq C\kappa^{\frac{3}{2}}e^{\kappa\delta
C_2(\mathfrak T_{\mathcal O},\delta_0)},\quad
\text{and}
\quad
\kappa\left|\mathrm H_1^{(1)}\!\left(
\kappa\rho_\br\right)\right|
\leq C(1+\kappa^{\frac{1}{2}})e^{\kappa\delta
C_2(\mathfrak T_{\mathcal O},\delta_0)},
\]
implying that
\[
\norm{D(\nabla G_{\mathbb C}^{(\kappa)})(\boldsymbol\delta_\br)}
\leq
C(1+\kappa^{\frac{3}{2}})e^{\kappa\delta
C_2(\mathfrak T_{\mathcal O},\delta_0)}.
\]
For $d=3$, \eqref{eq:radial_first_derivatives},
\eqref{eq:radial_second_derivatives}, and
\eqref{eq:rho_bounds_O}--\eqref{eq:Im_rho_bound_O} yield
\[
\norm{D(\nabla G_{\mathbb C}^{(\kappa)})(\boldsymbol\delta_\br)}
\leq
C(1+\kappa^2)e^{\kappa\delta
C_2(\mathfrak T_{\mathcal O},\delta_0)}.
\]
In summary, we obtain from \eqref{eq:gradient_kernel_B_complex_O} that
\begin{equation}\label{eq:devTNorm}
\begin{aligned}
\norm{\widehat\nabla_{\widehat\y}
\mathsf t^{(\kappa)}_{\br,\mathbb C}(\widehat\x,\widehat\y)}
&\leq C_3(\mathfrak T_{\mathcal O},\delta_0)
(1+\kappa^{\frac{d+1}{2}})
e^{\kappa\delta C_2(\mathfrak T_{\mathcal O},\delta_0)}.
\end{aligned}
\end{equation}

\textbf{Step 3: Shape holomorphy of $\mathsf{B}_\br^{(\kappa)}$.}
Define
\begin{equation*}
(\mathsf B^{(\kappa)}_{\br,\mathbb C}\phi)(\widehat\y)
\coloneqq\int_{\partial\mathcal O}
\mathsf t^{(\kappa)}_{\br,\mathbb C}(\widehat\x,\widehat\y)
\phi(\widehat\x)\,\dd s_{\widehat\x},
\end{equation*}
and observe that $\mathsf B^{(\kappa)}_{\br,\mathbb C}$ is an extension of $\mathsf B^{(\kappa)}_{\br}$ whose kernel is a non-singular function and holomorphic with respect to $\br$. Moreover, differentiation under the integral and the Cauchy-Schwarz inequality imply
\begin{equation*}
\norm{\mathsf B^{(\kappa)}_{\br,\mathbb C}\phi}_{H^1(\widehat\D)}
\leq C_4(\mathfrak T_{\mathcal O},\delta_0)
\Bigl(\norm{\mathsf t^{(\kappa)}_{\br,\mathbb C}}_
{L^\infty(\partial\mathcal O\times\widehat\D)}
+\norm{\widehat\nabla_{\widehat\y}\mathsf t^{(\kappa)}_{\br,\mathbb C}}_
{L^\infty(\partial\mathcal O\times\widehat\D;\mathbb C^d)}
\Bigr)
\norm{\phi}_{H^{\frac{1}{2}}(\partial\mathcal O)}.
\end{equation*}
Substituting \eqref{eq:TNorm} and
\eqref{eq:devTNorm} yields
\eqref{eq:BopNorm}.

For $\br\in\mathfrak T_{\delta,\mathcal O}\cap
W^{1,\infty}(\mathcal O;\IR^d)$, \eqref{eq:B_operator_real_O_r} follows from
\eqref{eq:kernel_B_complex_O} and the identity
$\boldsymbol c_{\partial\mathcal O,\br,\mathbb C}
=J_{\mathrm S,\partial\mathcal O,\br}\,
\boldsymbol n_{\mathcal O_\br}\circ\br$ from \textbf{Step~2}.
For $\br\in\mathfrak T_{\mathcal O}$, we have $\br\equiv\operatorname{Id}$
in a vicinity of $\partial\mathcal O$, so this agrees with \eqref{eq:B_operator}.
\end{proof}

\section{Auxiliary Complex Extension Bound}
\label{app:abstract_complex_bound}

The results below is inspired by the proof of
\cite[Theorem~3.1]{LafontaineSpenceWunsch2021}, which in turn uses the three-line theorem
in a rectangle \cite[Lemma~D.1]{DyatlovZworski2019}.

\begin{lemma}
Let $F$ be a complex Banach space, let $\rho>0$, and let
$\mathcal A:\mathfrak T_\rho\rightarrow F$ be holomorphic.
\label{lem:abstract_complex_bound}
Assume that $0<M_0\leq M_1$ and
\begin{equation}\label{eq:abstract_real_complex_bounds}
\begin{aligned}
\norm{\mathcal A(\br)}_F&\leq M_0,
&&\br\in\mathfrak T_\rho\cap W^{1,\infty}(\widehat\D;\IR^d),\\
\norm{\mathcal A(\br)}_F&\leq M_1,
&&\br\in\mathfrak T_\rho.
\end{aligned}
\end{equation}
Then, for any $R,\delta>0$ satisfying $\delta<R<\rho-\delta$, one has
\begin{equation}\label{eq:abstract_exponential_bound}
\sup_{\br\in\mathfrak T_\delta}\norm{\mathcal A(\br)}_F
\leq M_0e^{\frac{4\delta}{\pi R}\log\left(\frac{M_1}{M_0}\right)}.
\end{equation}
The same assertions hold with $\mathcal O$ and $\mathfrak T_{\mathcal O}$
in place of $\widehat\D$ and $\mathfrak T$, respectively.
\end{lemma}

\begin{proof}
\textbf{Step 1: Reduction to one complex variable.}
By definition of $\mathfrak T_\delta$, for every $\br\in\mathfrak T_\delta$ there is $\widetilde{\br}\in\mathfrak T$ such that
$\norm{\br-\widetilde{\br}}_{W^{1,\infty}(\widehat\D;\IC^d)}<\delta$. Set $t\coloneqq\norm{\operatorname{Im}\{\br\}}_{W^{1,\infty}(\widehat\D;\IR^d)}$.
If $t=0$, the first estimate in \eqref{eq:abstract_real_complex_bounds} gives $\norm{\mathcal A(\br)}_F\leq M_0$.
Otherwise, set $\br_{\mathbb R}\coloneqq\operatorname{Re}\{\br\}$ and
$\boldsymbol{h}\coloneqq\frac{\operatorname{Im}\{\br\}}{t}$, and consider the map
$z\mapsto\br_{\mathbb{R}}+z\boldsymbol{h}$.
Here $\boldsymbol{h}$ is real and $\norm{\boldsymbol{h}}_{W^{1,\infty}(\widehat\D;\IC^d)}=1$. The map coincides with $\br$ at
$z=\imath t$ and is real for real $z$.
Moreover, $\widetilde{\br}$ being real implies
$t\leq\norm{\br-\widetilde{\br}}_{W^{1,\infty}(\widehat\D;\IC^d)}<\delta$ and
$\norm{\br_{\mathbb R}+z\boldsymbol{h}-\widetilde{\br}}_{W^{1,\infty}(\widehat\D;\IC^d)}
\leq\norm{\br-\widetilde{\br}}_{W^{1,\infty}(\widehat\D;\IC^d)}+|z|
<\delta+|z|$ for $z\in\mathbb C$. Thus, by definition of $\mathfrak T_\rho$,
the map $z\mapsto\mathcal A(\br_{\mathbb R}+z\boldsymbol{h})$
is holomorphic in the open disc $|z|<\rho-\delta$ containing $|z|\leq R<\rho-\delta$.
Thus, both estimates in \eqref{eq:abstract_real_complex_bounds} apply
for $\br_{\mathbb R}+x\boldsymbol{h} \in \mathfrak T_\rho$, $x\in[-R,R]$.

\textbf{Step 2: Holomorphic coordinate transform.}
Taking the logarithm with respect to the principal branch, consider the function
\[
w(z)\coloneqq-\frac{2\imath}{\pi}\log\frac{R+z}{R-z},
\]
on $D_R^+\coloneqq\{z\in\mathbb C:|z|<R,\ \operatorname{Im}z>0\}$.
The map $z\mapsto\frac{R+z}{R-z}$ maps $D_R^+$ onto the first quadrant of the complex plane. Hence $w$ is holomorphic, and for $z=x+\imath y\in D_R^+$ we have
\[
\operatorname{Re}\{w(z)\}
=\frac2\pi\arctan\frac{2Ry}{R^2-x^2-y^2}.
\]
Thus, the real part of $w$ takes values in $(0,1)$ on $D_R^+$, tends to zero on $(-R,R)$ and to 1 on the upper semicircle away from $\pm R$.

\textbf{Step 3: Boundary estimates.}
Observing that $\log\frac{M_1}{M_0}\geq0$, we define the holomorphic function
\begin{equation}\label{eq:H}
H(z)\coloneqq M_0^{-1}\mathcal A(\br_{\mathbb R}+z\boldsymbol{h})e^{-\log\left(\frac{M_1}{M_0}\right)w(z)},
\quad
z \in D^+_R,
\end{equation}
which satisfies $\norm{H(z)}_F\leq\frac{M_1}{M_0}$ in $D_R^+$.
Moreover, the modulus of $e^{-\log\left(\frac{M_1}{M_0}\right)w(z)}$ tends to $1$ on $(-R,R)$ and to
$\frac{M_0}{M_1}$ on the upper semicircle, away from $\pm R$.
Thus, for every $\zeta\in\partial D_R^+\setminus\{-R,R\}$,
\eqref{eq:abstract_real_complex_bounds} yields
\[
\limsup_{D_R^+\ni z\to\zeta}\norm{H(z)}_F\leq1.
\]
At $\pm R$, continuity of $z\mapsto\mathcal A(\br_{\mathbb R}+z\boldsymbol{h})$ and
$e^{-\log\left(\frac{M_1}{M_0}\right)\operatorname{Re}\{w(z)\}}\leq1$ imply
\[
\limsup_{D_R^+\ni z\to\pm R}\norm{H(z)}_F
\leq M_0^{-1}\norm{\mathcal A(\br_{\mathbb R}\pm R\boldsymbol{h})}_F\leq1.
\]

\textbf{Step 4: Application of the maximum modulus principle.}
To apply the scalar maximum modulus principle, consider any bounded linear functional $\lambda\in F'$ with $\norm{\lambda}_{F'}\leq1$.
The considerations in Step 3 imply that
\[
\limsup_{D_R^+\ni z\to\zeta}|\lambda(H(z))|
\leq\norm{\lambda}_{F'}
\limsup_{D_R^+\ni z\to\zeta}\norm{H(z)}_F
\leq 1
\]
for every $\zeta\in\partial D_R^+$. Since $z\mapsto\lambda(H(z))$ is holomorphic on $D_R^+$, the scalar maximum modulus principle
\cite[Theorem~VI.1.4]{Con78} yields
$|\lambda(H(z))|\leq 1$ for all $z\in D_R^+$.
This implies
\[
\norm{H(z)}_F
=\sup_{\substack{\lambda\in F'\\\norm{\lambda}_{F'}\leq1}}
|\lambda(H(z))|
\leq1,\quad z\in D_R^+.
\]
At $z=\imath t$, we have $\br_{\mathbb R}+\imath t\boldsymbol{h}=\br$ and $w(\imath t)=\frac{4}{\pi}\arctan\frac{t}{R}$.
Using $\norm{H(\imath t)}_F\leq1$ in \eqref{eq:H} gives $\norm{\mathcal A(\br)}_F\leq M_0e^{\log\left(\frac{M_1}{M_0}\right)w(\imath t)}$.
Since $t<\delta$ and $\frac{M_1}{M_0}\geq1$, we obtain
\begin{equation*}
\sup_{\br\in\mathfrak T_\delta}\norm{\mathcal A(\br)}_F
\leq M_0\left(\frac{M_1}{M_0}\right)^{\frac{4}{\pi}\arctan\frac{\delta}{R}},
\end{equation*}
from which $\arctan\frac{\delta}{R}\leq\frac{\delta}{R}$ yields
\eqref{eq:abstract_exponential_bound}.
\end{proof}

\begin{corollary}\label{cor:abstract_frequency_bound}
Let $\kappa_0>0$ and $0<R<\rho$.
Suppose that, for each $\kappa\geq\kappa_0$, $\mathcal A_\kappa$ satisfies the assumptions of
Lemma~\ref{lem:abstract_complex_bound} on
$\mathfrak T_\rho$, with
\[
M_0=C_0(1+\kappa)^\alpha
\quad
\text{and}
\quad
M_1=C_1(1+\kappa)^\beta e^{\kappa\rho K},
\]
where $0<C_0\leq C_1$, $\alpha,\beta\in\mathbb R$ with
$\alpha\leq\beta$, and $K\geq0$.
Assume that $\rho$, $R$, $C_0$, $C_1$, $\alpha$, $\beta$, and $K$
are independent of $\kappa$ and $\delta$.
Then, for all $\kappa\geq\kappa_0$ and $0<\delta<\min\{R,\rho-R\}$,
\begin{equation}\label{eq:abstract_frequency_bound}
\sup_{\br\in\mathfrak T_\delta}\norm{\mathcal A_\kappa(\br)}_F
\leq C_1(1+\kappa)^\alpha
 e^{\frac{4\delta}{\pi R}
 \left(\kappa\rho K+(\beta-\alpha)\log(1+\kappa)\right)}.
\end{equation}
\end{corollary}

\begin{proof}
$C_0\leq C_1$ implies that the first estimate in \eqref{eq:abstract_real_complex_bounds}
also holds with $\widetilde M_0\coloneqq C_1(1+\kappa)^\alpha$.
Applying \eqref{eq:abstract_exponential_bound} with $\widetilde M_0$ in place of $M_0$ and using $\log(M_1/\widetilde M_0)=\kappa\rho K+(\beta-\alpha)\log(1+\kappa)$ yields
\eqref{eq:abstract_frequency_bound}.
\end{proof}

\section*{Declaration on the Use of Generative Artificial Intelligence}

Generative AI was used to assist in reviewing mathematical arguments,
searching the literature, in particular for the discussion of wavenumber-explicit
bounds in Remark~\ref{rem:real_frequency_comparison}, and revising the wording of the manuscript.
This included assistance in organizing and presenting the results in
Appendix~\ref{app:abstract_complex_bound} and with Step~3 of the proof of
Theorem~\ref{thm:DLhlmrp}. The authors have checked the mathematical arguments
and assume responsibility for all content.

\bibliographystyle{siamplain}
\bibliography{references}

\end{document}